\documentclass{amsart}
\usepackage{amssymb}
\usepackage{amsmath}
\usepackage[hyperref]{degt}
\usepackage{hyperref}
\def\itemrefform#1{$(#1)$}

\addtheorem{Algorithm}{algorithm}

\def\2{\color{red}}
\def\3{\color{blue}}
\def\4{\color{magenta}}
\def\Dg:{{\bf Dg:\enspace}\ignorespaces}

\let\LE=\Lambda
\let\QF=\Phi
\let\TC=\Psi
\let\QC=\Theta
\let\SC=\Delta

\def\KK-{\pdfstr{K3}{$K3$}-\penalty0\ignorespaces}

\def\onefam#1{\multispan2 \hss$\dim=#1$\hss}
\def\upssrm#1{^{\scriptscriptstyle\mathrm{#1}}}
\def\0{}
\def\A{\upssrm{A}}

\def\K{\upssrm{K}}
\def\Kummer{\frak{K}}
\def\dual{^\vee}
\def\F{\Bbb F}

\let\graph\Gamma
\let\pencil\Pi
\let\fiber\Sigma
\let\Weyl\Delta
\def\fp{\Weyl^{\sharp}}
\def\barfp{\bar\Weyl^{\sharp}}
\def\Fano{\Cal{F}}
\def\Cone{\Cal{C}}
\def\Fn{\operatorname{Fn}}
\def\extended{^{\mbox{\tiny ex}}}
\def\extended{^\text{\rm ex}}
\def\Fex{\Fn\extended}
\def\rt{\operatorname{\frak{rt}}}
\def\root{\operatorname{root}}
\def\sat{\operatorname{sat}}
\def\satex{\sat\extended}
\def\girth{\operatorname{girth}}
\def\val{\operatorname{val}}
\def\Star{\operatorname{star}}
\def\spp{\operatorname{sp}}
\def\NS{\operatorname{NS}}
\def\Pic{\operatorname{Pic}}

\def\Nef{\operatorname{Nef}}
\def\hp{h}
\def\pdp{p}
\def\kappal{\mathfrak{l}}

\def\graphl{\operatorname{Fn}}
\def\graphe{\graphl\extended}
\def\graphex{\graph\extended}
\def\PPP{\mathbb P}
\def\KKK{\mathbb K}
\def\bS{\bar S}
\def\qgen{spanned }
\def\sstar{star}
\def\algcit{Algorithm}
\def\corcit{Corollary}
\def\lemcit{Lemma}
\def\remcit{Remark}
\def\propcit{Proposition}

\def\thmcit{Theorem}

\def\resp.{resp\PERIOD}
\let\geq\ge
\let\leq\le

\let\nodal\alpha
\let\M=M
\let\o=o

\let\ex=e
\let\iso=p

\def\basis#1{\{#1\}}
\def\pos#1{P_{#1}}
\def\eclass#1{[\![#1]\!]}

\def\CK{\Cal K}

\def\GK{\frak G}
\def\bL{\bold{L}}

\def\CC{\Cal C}

\def\CS{\Cal S}
\def\fo{\frak o}

\makeatletter
\def\newconfig#1#2{\expandafter\gdef\csname config=#1\endcsname{#2}}

\def\theconfig#1{\hyperlink{config=#1}{\csname config=#1\endcsname}}
\def\confighook#1#2{\smash{\raise9pt\hbox{\hypertarget{config=#1}{}}}%
 \protected@write\@auxout{}{\string\newconfig{#1}{#2}}%
 \csname config=#1\endcsname}
\pdef\aconfig#1#2{\ifmmode\confighook{#1}{#2}\else$\confighook{#1}{#2}$\fi}
\pdef\config#1{\ifmmode\theconfig{#1}\else$\theconfig{#1}$\fi}
\makeatother

\def\tpencil{\smash{\tilde\pencil}}

\makeatletter
\let\HyPsd@CatcodeWarning\@gobble
\makeatother

\let\fiber=\Sigma
\let\graph\Gamma
\let\pencil\Pi
\def\Fano{\Cal{F}}
\def\CK{\Cal K}
\def\CS{\Cal S}
\def\CV{\Cal V}
\def\bCS{\bar\CS}
\def\F{\Bbb F}
\def\GK{\frak G}
\def\val{\operatorname{val}}

\def\sat{\operatorname{sat}}
\let\nodal=\alpha

\def\bv{\bold v}
\def\gram{\frak{m}}
\def\perm(#1){\hyperlink{@-perm}{#1'}}

\let\gcup\sqcup
\def\Sym{\operatorname{Sym}}
\def\num{\operatorname{nc}}
\def\Num{\operatorname{NC}}
\def\true{\text{true}}
\def\false{\text{false}}
\def\smax{_\text{max}}
\def\smin{_\text{min}}
\def\maxlist{\Cal G\smax}

\def\Mlines{M_{\text{lines}}}
\def\Msing{M_{\text{sing}}}
\def\Sat_#1{\hyperlink{@-Sat}{\operatorname{Sat}_#1}}
\def\pat{\hyperlink{@-pat}{\operatorname{pat}}}
\def\gpat_#1{\hyperlink{@-scount}{\pi_{#1}}}
\def\Pat{\hyperref[eq.pat.Gamma]{\operatorname{Pat}}}
\def\range_#1{\hyperlink{@-range}{\operatorname{range}_{#1}}}
\def\Sec_#1{\hyperref[eq.Sec]{\operatorname{Sec}_{#1}}}
\def\bnd_#1{\hyperref[eq.sec.bound]{b_{#1}}}
\def\sym#1{\hyperlink{@-sym}{[#1]}}
\def\Com{\hyperlink{@-com}{C}}
\def\Coms{\hyperlink{@-com}{C_*}}
\def\tpencil{\smash{\tilde\pencil}}

\makeatletter
\def\anchor#1{\vbox to\z@{\vss\hbox{\raise9pt\vbox to\z@{\hypertarget{@-#1}{}}}}}
\makeatother

\title{Counting Lines with Vinberg's algorithm}

\author{Alex Degtyarev}

\address{%
	Department of Mathematics\\
	Bilkent University\\
	06800 Ankara, TURKEY}

\email{
	degt@fen.bilkent.edu.tr}

\author{S\l awomir Rams}
\address{Institute of Mathematics, Jagiellonian University,
	ul. {\L}ojasiewicza 6,  30-348 Krak\'{o}w, Poland}
\email{slawomir.rams@uj.edu.pl}

\thanks{%
	A.D. was partially supported by the T\"{U}B\DOTaccent{I}TAK grant 118F413.
    S.R. was   supported by the National Science Centre, Poland, Opus  grant
	no.\ 2017/25/B/ST1/00853.}

\dedicatory{To the memory of\, $\boxed{\text{\!Ernest Borisovich Vinberg}}$}



\begin{document}

	\begin{abstract}
	 We combine classical Vinberg's algorithms \cite{Vinberg:polyhedron} with the
	 lattice-theoretic/arithmetic approach
	 from \cite{degt:lines} to
	 give a method of classifying large line configurations on
	 complex quasi-polarized \KK-surfaces. We apply our method to
	 classify all complex \KK-octic surfaces
	 with at worst Du Val singularities and at least $32$ lines.
     The upper bound on the number of lines is $36$, as in the smooth
     case, with at most $32$ lines if the singular locus is non-empty.
	\end{abstract}
	
	\maketitle

\section{Introduction}\label{S.intro}

In the last decade,
a substantial progress towards the complete
understanding of configurations of lines on \emph{smooth polarized} \KK-surfaces
 has been made (see \autoref{s.sart} below).
Unfortunately,  the methods that led to
that developement
do not yield optimal results (if any) when the surfaces in question have singularities.
The main aim of this paper is to address this problem.
More precisely, we combine classical Vinberg's algorithms \cite{Vinberg:polyhedron} with the
lattice-theoretic/arithmetic approach
from \cite{degt:lines} to create a uniform framework for the study of
configurations of lines on complex projective  $K3$ surfaces
with at worst
Du Val (\latin{aka} simple or $\bA$--$\bD$--$\bE$) singularities.

Among other things,
we discover a new phenomenon unthinkable in the realm
of smooth polarized \KK-surfaces: a surface with a larger N\'{e}ron--Severi
lattice~$N'$ may have fewer lines than that with a smaller lattice
$N\subset N'\subset H_2(X)$; we discuss this phenomenon in
\autoref{s.extendedversusplain}.
In particular, unlike the smooth case, it is no longer sufficient to
confine oneself to surfaces~$X$ with $\NS(X)$ \qgen by
lines
(and the quasi-polarization, see \autoref{conv.lines}).
It \emph{is}, however, sufficient to assume that $\NS(X)$ is \qgen by
lines \emph{and exceptional divisors}.
Hence,
in most statements, we make an assumption on the number of
lines, \ie, the size of the \emph{plain} Fano graph
$\graphl(X,\hp)$, \latin{aka} the dual adjacency graph
of lines on~$X$, see~\eqref{eq-fanograph-def},
but we classify \emph{extended} Fano graphs $\graphe(X,\hp)$,
\ie, bi-colored
graphs of both lines and exceptional divisors,
see~\eqref{eq-def-exFano-surface}.
Obviously, these extended
graphs also give us more detailed information about the lines and
singular points of
the original projective surface $X_{2d}\subset\Cp{d+1}$
itself.

To test our approach, we classify \KK-octics with many lines and
at worst
Du Val
singularities.
The principal results on octics are stated in Theorems~\ref{thm-main}
and~\ref{th.special},
where we distinguish \emph{triquadrics} (ideal theoretical intersections
of three quadrics) \vs. \emph{special octics} (see \autoref{def.triquadric}).
Recall that a complex \KK-surface is called \emph{singular} if
its
Picard rank is maximal: $\rho(X)=20$;
when mapped to a projective space, such surfaces are projectively rigid.
In general,
when speaking about an ``$s$-parameter family'', we always mean
the dimension $20-\rho(X)$ modulo the projective group.

\theorem[see \autoref{S.Proof}]\label{thm-main}
Let $X_{8} \subset \Cp{5}$ be
a degree $8$ \KK-surface  with
at worst Du Val singularities.
Then $X_8$ contains at most $36$
\rom(at most~$32$
if  $\Sing X_8 \neq \emptyset$\rom) lines. Moreover,
if $X_8$ contains at least $32$ lines, then it is one of the surfaces
listed in
\autoref{tab.main}.
Thus, if $X_8$ has $32$ lines and $\Sing X_8\ne\varnothing$, then $X_8$ is a
triquadric\rom: there are two connected $1$-parameter families and four
singular surfaces.
\endtheorem

We   observe 
(see \autoref{th.Kummer}) that there are \emph{exactly} two disjoint
$3$-parameter
families of \emph{Kummer octics}
with
a distinguished set of $16$
Kummer divisors
mapped to lines; simple as it is,
we could not find this statement in the literature.
Each family consists of triquadrics, and each contains octics with
$32$ lines, \viz. \config{QC.32.K} (generic in one family)
and \config{QC.32.4} (rigid in the other one) in \autoref{tab.main}.
See \autoref{s.Kummer} for details.

Special octics constitute a codimension~$1$ family in the space of all
octics and are subject to a stronger bound.

\theorem[see \autoref{proof.special}]\label{th.special}
There are two connected families of special octics
with at worst Du Val singularities
and at least  $30$ lines\rom:
a $2$-parameter
family of smooth surfaces with $33$ lines \rom(see
$\config{TC.33}$ in \autoref{tab.main}\rom)
and a $1$-parameter family of surfaces with
three $\bA_1$-type points and
$30$ lines \rom(see $\config{TC.30.3}$ in
\autoref{tab.mainsmall}\rom).
\endtheorem

\table{\rm
	\caption{{\KK-octics with at least $32$ lines (see \autoref{thm-main} and
			\autoref{s.classification})}}\label{tab.main}
	\def\(#1){_{#1}}
	\def\MMi#1{\expandafter\MMii#1\endMM}
	\def\MMii#1,#2\endMM{$#1,#2$}
	\def\-{\rlap{$\mdag$}}
	\def\*{\rlap{$\mstar$}}
	\def\r{\rlap{$\mreal$}}
	\def\={\relax\afterassignment\ul\count0=}
	\def\ul{\underline{\the\count0}}
	\def\u#1{\mathbf{#1}}
	\hbox to\hsize{\hss\vbox{\halign{\strut\quad
				$#$\hss\quad&\hss$#$\hss\quad&\hss$#$\hss\quad&
				\hss$#$\hss\quad&\hss$#$\hss\quad&\hss$#$\hss\quad&$#$\hss&
				\quad#\hss\cr
				\noalign{\hrule\vspace{2pt}}
				\graphex&\!\!\Sing X_8\!\!&\ls|\Aut\graphex|&\!{\det}\!&(r,c)&
				\!\!\ls|\Aut X_8|\!\!&\,\,T:=\NS(X)^\perp&\,\,Ref\cr
				\noalign{\vspace{1pt}\hrule\vspace{2pt}}
				\aconfig{QC.36'}{\QC_{36}'}
				&      &  64&  32&(1,0)&  16&[4,0,8]        &\ref{lem.quad}, \ref{lem.a3p}\cr
				\aconfig{QC.36''}{\QC_{36}''}
				&      & 576&  36&(1,0)&  72&[6,0,6]        &\ref{lem.quad}, \ref{lem.a3p}\cr
				\aconfig{QC.34'}{\QC_{34}'}
				&      &  96&    &\onefam1  &\bU\(2)\oplus[12]&\ref{lem.quad}, \ref{lem.a3p}\cr
				\aconfig{QC.33}{\QC_{33}}
				&      & 192&  80&(1,0)&  24&[8,4,12]       &\ref{prop.Kummer-}, \ref{lem.a3p}\cr
				\aconfig{TC.33}{\TC_{33}}
				&      &6912&  36&\onefam2  &\bU^2\(3)      &\ref{lem.a2p}\cr
				\aconfig{QC.32}{\QC_{32}^{}}
				&      &  96&  60&(1,0)&  24&[4,2,16]       &\ref{lem.a3p}\cr
				\aconfig{QC.32'}{\QC_{32}'}
				&      & 384&    &\onefam2  &\bU\(2)\oplus[-4]\oplus[4]&\ref{lem.quad}\cr
				\aconfig{QC.32''}{\QC_{32}''}
				&      & 512&    &\onefam2  &\bU\(2)\oplus\bU\(4)&\ref{lem.quad}\cr
				\aconfig{QC.32'''}{\QC_{32}'''}
				&      & 768&  36&\onefam2  &\bU^2\(3)       &\ref{lem.quad}\cr
				\aconfig{QC.32.K}{\QC_{32}\K}
				&      &23040& 32&\onefam3  &\bU^2\(2)\oplus[-4]&\ref{prop.Kummer.64}, \ref{lem.quad}\cr
				\aconfig{QC.32.4}{\QC_{32}^4}
				&4\bA_1& 256&  16&(1,0)& 128&[4,0,4]        &\ref{prop.Kummer.256}, \ref{lem.quad}, \ref{lem.a3p}\cr
				\aconfig{QC.32.2'}{\QC_{32}^{2\prime}}
				&2\bA_1&  48&  27&(1,0)&  24&[6,3,6]        &\ref{lem.quad}\cr
				\aconfig{QC.32.2''}{\QC_{32}^{2\prime\prime}}
				&2\bA_1&  64&  36&(1,0)&  16&[6,0,6]        &\ref{lem.quad}\cr
				\aconfig{QC.32.2'''}{\QC_{32}^{2\prime\prime\prime}}\!
				&2\bA_1& 128&    &\onefam1  &\bU\(2)\oplus[8]&\ref{lem.quad}, \ref{lem.a3p}\cr
				\aconfig{QC.32.1'}{\QC_{32}^{1\prime}}
				& \bA_1&  16&  32&(1,0)&   8&[6,2,6]        &\ref{lem.quad}, \ref{lem.a3p}\cr
				\aconfig{QC.32.1''}{\QC_{32}^{1\prime\prime}}
				& \bA_1&  96&  32&\onefam1  &\bU\(3)\oplus[6]&\ref{lem.quad}\cr
				\noalign{\vspace{1pt}\hrule\vspace{2pt}}
				\crcr}%
}\hss}%
\quad
In the tables, we abbreviate $\bU_n:=\bU(n)=[0,n,0]$, \cf.
\autoref{s.notation}\hfill
}
\endtable

\table{\rm
	\caption{{Other \KK-octics with many lines (see
			\autoref{s.classification})}}\label{tab.mainsmall}
	\def\(#1){_{#1}}
	\def\MMi#1{\expandafter\MMii#1\endMM}
	\def\MMii#1,#2\endMM{$#1,#2$}
	\def\-{\rlap{$\mdag$}}
	\def\*{\rlap{$\mstar$}}
	\def\r{\rlap{$\mreal$}}
	\def\2#1{\!{\le}2#1\!}
	\def\={\relax\afterassignment\ul\count0=}
	\def\ul{\underline{\the\count0}}
	\def\u#1{\mathbf{#1}}
	\hbox to\hsize{\hss\vbox{\halign{\strut\quad
				$#$\hss\quad&\hss$#$\hss\quad&\hss$#$\hss\quad&
				\hss$#$\hss\quad&\hss$#$\hss\quad&\hss$#$\hss\quad&$#$\hss&
				\quad#\hss\cr
				\noalign{\hrule\vspace{2pt}}
				\graphex&\Sing X_8&\ls|\Aut\graphex|&\det&(r,c)&
				\!\!\ls|\Aut X_8|\!\!&T:=\NS(X)^\perp&Ref\cr
				\noalign{\vspace{1pt}\hrule\vspace{2pt}}	
				\aconfig{QF.30'}{\QF_{30}'}
				&      & 240& 140&(1,1)&  60&[12,2,12]      &\ref{lem.a4p}\cr
				\aconfig{QF.30''}{\QF_{30}''}
				&      &  40& 135&(1,1)&  10&[12,3,12]      &\ref{lem.a4p}\cr
				\aconfig{QF.30.5''}{\QF_{30}^{\prime\prime5}}
				&5\bA_1&  40&  15&(1,0)&  20&[2,1,8]        &\ref{lem.a4p}\cr
				\aconfig{TC.30.3}{\TC_{30}^3}
				&3\bA_1& 864&    &\onefam1  &\bU\(3)\oplus[6]&\ref{lem.a2p}\cr
				\aconfig{SC.28'}{\SC_{28}'}
				&      & 576&  44&(1,0)&  96&[4,2,12]       &\ref{prop.Kummer.256}\cr
				\noalign{\vspace{1pt}\hrule\vspace{2pt}}
				\aconfig{QC.28.2+1}{\QC_{28}^{2,1}}&2\bA_2,\bA_1
				&  16&  16&(1,0)&  16&[2,0,8]        &\ref{lem.quad}\cr
				\aconfig{QC.25.1+4}{\QC_{25}^{1,4}}&\bA_2,4\bA_1
				&   2&  16&(1,0)&   2&[2,0,8]        &\ref{lem.a3p}\cr
				\noalign{\vspace{1pt}\hrule}
				\crcr}}\hss}
}
\endtable

The first column of Tables~\ref{tab.main} and~\ref{tab.mainsmall} refers to the
extended Fano graphs of octics; the other entries are explained in
\autoref{s.classification} below.
For obvious reasons, we do not try to depict these graphs;
the precise descriptions are
available electronically
(in the form of a \texttt{GRAPE}
\cite{GRAPE:nauty,GRAPE:paper,GRAPE} records)
in~\cite{DR:octic.graphs}.
Most graphs do not have a known ``name'' and, therefore,
\latin{de facto} these descriptions are their definitions.

Thus,
we obtain
a complete picture of the large line configurations on a class of varieties
(\viz. complete intersections of quadrics)
that
have been a
subject of research since the XIX-th century.
As long as \KK-surfaces are concerned,
it is well understood that
the larger the integer $d$, the smaller the maximal Fano graphs
of smooth complex \KK-surfaces of degree $2d$ are (see \cite{degt:lines}).
Moreover,
if $d$ is
sufficiently large,
no Fano graph of a
degree $2d$ \KK-surface can be hyperbolic (\cf. \autoref{S.taxonomy}).
Thus,
our approach can as well be applied
to classify large Fano graphs of quasi-polarized \KK-surfaces $(X,\hp)$ with
$\hp^2 > 8$.  On the other hand, the case of sextics (\cf. \cite{Degtyarev.Rams.sextics}),
quartics (\cf. \cite{Degtyarev.Rams.triangular}), and,
especially, double planes ramified at sextic curves would require a
considerably more thorough
treatment of the
configurations
containing
a triangle or a quadrangle.

As yet another justification of our interest in the problems above,
we recall that
\KK-surfaces are $2$-dimensional hyperk\"{a}hler varieties.
We hope that our new algorithms/\penalty0methods
developed in dimension~$2$,
apart from being of interest on their own, may contribute to a better understanding of the
higher-dimensional case.

\subsection{Classification of \KK-octics with many lines}\label{s.classification}
We use the  girth (\cf. \autoref{s.misc}) to subdivide
\emph{plain} Fano graphs of quasi-polarized \KK-surfaces into several
classes (see \autoref{S.taxonomy}) and
obtain a more refined classification/bound for each class.
A Fano graph $\graph:=\graphl(X,\hp)$ is called
\roster*
\item
\emph{triangular} (the $\TC_*$-series),
if $\girth(\graph)=3$; such graphs
appear only as Fano graphs of  special octics (see \autoref{lem.2stars}),
\item
\emph{quadrangular} (the $\QC_*$-series),
if $\girth(\graph)=4$,
\item
\emph{pentagonal} (the $\QF_*$-series),
if $\girth(\graph)=5$,
\item
\emph{astral} (the $\SC_*$-series),
if $\girth(\graph)\ge6$ and $\graph$ has a vertex
of valency $\ge4$.
\endroster
All other graphs are \emph{locally elliptic} (the $\LE_*$-series),
\ie, one has $\val v\le3$ for each
vertex $v\in\graph$
(and we still assume $\girth(\graph)\ge6$ to exclude a few trivial
cases).

Principal properties of \KK-octics with large  line configurations are
collected in Tables~\ref{tab.main} and~\ref{tab.mainsmall}, where,
inevitably, we have to restate
some results of~\cite{degt:lines} concerning smooth octics.
The first column refers
to the isomorphism classes of the extended
 Fano graphs introduced
 elsewhere in the paper; the subscript always stands for the number of
 lines.
 Then, for each ~graph  $\graphex:=\graphe(X,\hp)$, we list,
\roster*
\item
the number and types of the singular points of the corresponding octics;
  \item
 the order $\ls|\Aut\graphex|$ of the full automorphism group of~$\graphex$,
 \item
 the transcendental lattice
 $T:=\NS(X)^\perp\subset H_2(X;\Z)$  of generic
  \KK-surfaces $X$ with $\graphe(X,\hp)
 \cong\graphex$;
 this lattice determines a family of \emph{abstract} \KK-surfaces
 (see \autoref{s.graph2T} for the discussion of the uniqueness),
 \item
 references to the parts of the text where the graph appears.
 \endroster
 For the rigid configurations ($\rank T=2$), we list, in addition,
 \roster*
 \item
 the determinant $\det T$,
 \item
 the numbers
 $(r,c)$ of, respectively, real projective isomorphism classes and
 pairs of complex conjugate projective isomorphism
 classes of octics
$(X,\hp)$
with  $\graphe(X,\hp)\cong\graphex$,
 \item
 the order $\ls|\Aut X_8|$ of the
 group of \emph{projective} automorphisms of~$X_8$.
 \endroster
  Each of the nine non-rigid configurations~$\graph$ listed in
the tables  
  is realized by a single connected
 equilinear deformation family~${\mathcal M}(\graph)$; we indicate the dimension
 $\dim\bigl({\mathcal M}(\graph)/\!\PGL(6,\C)\bigr)=\rank T-2$ and the minimum of the discriminants
 of the singular \KK-surfaces in ${\mathcal M}(\graph)$
 (whenever we know its value).

The proof of \autoref{thm-main} is based on the study of various types of Fano graphs of \KK-octics.
In particular,
\autoref{thm-main}
implies
that each triquadric with at least $32$ lines contains a quadrangle (\ie, an $\tA_3$ configuration of lines) and its singularities
(if any) are of type $\bA_1$.
As
part of the proof, we obtain
a classification
 of maximal pentagonal configurations (see  \autoref{s.pent}), maximal astral configurations  (see  \autoref{s.astral})
 and
 examples of large  line configurations on octics with $\bA_2$-singularities (the entries $\QC_{28}^{2,1}$, $\QC_{25}^{1,4}$
 in  \autoref{tab.mainsmall}). Combined with Nikulin's theory~\cite{Nikulin:forms}
 (\cf. \autoref{S.Proof}), this yields the following extra bounds.

\remark[see \autoref{S.Proof}]\label{thm-main-degenerate}
Let $X_{8} \subset \Cp{5}$ be a degree $8$ \KK-surface  with
at worst Du Val singularities
and $\graph$ its Fano graph. Then\rom:
\roster
\item
if $\graph$ is pentagonal, then $\ls|\graph|\le30$ and the maximum is
attained at the three singular surfaces with the extended Fano graphs $\QF_*$ in
\autoref{tab.mainsmall}\rom;
\item
if $\graph$ is astral, then $\ls|\graph|\le28$ and the maximum is attained at
a unique singular Kummer octic with the extended Fano graph \config{SC.28'} in
\autoref{tab.mainsmall}.
\endroster
\endremark
For completeness,
the extremal locally elliptic graphs ($\LE_*$ of size~$24$ or~$25$) are
described
in \autoref{s.elliptic}: they are not listed in \autoref{tab.mainsmall} as they
are realized by too many \KK-octics with non-isomorphic transcendental
lattices.


\subsection{Contents of the paper}\label{s.contents}

The paper splits into two parts. The first one, \viz.
\autoref{S.Fano}--\autoref{S.taxonomy} and \autoref{S.basicalg},
describes a general strategy for
the classification of large configurations of lines on projective models of
complex  \KK-surfaces with at worst Du Val singularities.
The second part, \autoref{S.extreme}--\autoref{S.Proof} and
Appendices~\ref{S.Msigma},~\ref{S:largepencils},
demonstrates the effectiveness of our approach in the case of
\KK-octics.

In \autoref{S.Fano} we
lay a
theoretical foundation for
the computation
of the configuration of lines on a quasi-polarized \KK-surface $(X, \hp)$
in terms of the lattice $\NS(X) \ni \hp$.
Unlike the smooth case~\cite{degt:lines,DIS},   
in order to avoid overcounting (see \autoref{rem.vample}),
we have to use Vinberg's algorithm~\cite{Vinberg:polyhedron} and
compute two layers of the fundamental polyhedron (\cf. \eqref{eq-def-exFano-lattice}).
The result depends on the choice of a Weyl chamber, and we discuss the extent
to which it is well defined (see Lemmas \ref{lem.Weyl.independent}
and~\ref{lem.compatible}).
Then, in \autoref{s.linesonK3}, we recall the geometry
of the  nef cone of $X$
and relate our abstract construction to
the geometric set of lines on~$X$ (\autoref{thm-each-chamber-ok}).
The section concludes with a discussion of
Saint-Donat's conditions~\cite{Saint-Donat} for various
degenerations of the quasi-polarization
(see \autoref{s.admissible}, \autoref{s.hyperbolic}).

In \autoref{S.graphs} we study polarized lattices generated by lines and,
thus, constructed from graphs. The principal innovation here is the concept
of extensibility (\autoref{def.ext}) which is to replace the admissibility
conditions used in the smooth case. A simple criterion is given by
\autoref{lem.ext}, which also asserts that, in spite of the ambiguity in the
choice of a Weyl chamber, each extensible graph has a well-defined
saturation.
The condition for the geometric realizability of a graph is given by
Theorems~\ref{th.K3} and~\ref{th.K3.ext}; we also state a version for special
octics (\autoref{th.octics}, our primary concern) and hyperelliptic
polarizations (\autoref{th.hyperelliptic}, very similar).
Finally, we address the new phenomenon mentioned at the beginning (see
\autoref{wrn.sat} and a detailed discussion in
\autoref{s.extendedversusplain}) and explain how it affects our proof
strategy.

In \autoref{S.taxonomy} we
recall
(after~\cite{degt:lines})
the taxonomy of hyperbolic graphs
and combinatorial counterparts of elliptic pencils on \KK-surfaces.
In \autoref{s.approach} we
explain our approach to the classification of large
configurations
of lines.
After a thorough examination of the local properties
specific to Fano graphs of octics (and proving
\autoref{th.special}) in \autoref{S.extreme},
this general approach, mostly computer aided, is illustrated in human readable form
in \autoref{S.Kummer}, on the example of (almost) Kummer octics. In
\autoref{s.quartics} we also announce a few new results
(mostly examples) concerning spatial quartics.

In \autoref{S.triquadrics} we state the results of the computation in the
form of a number of bounds for various types of graphs. These statements are
used in \autoref{S.Proof} to prove \autoref{thm-main}.
The computation leading to \autoref{S.triquadrics} is heavily computer aided
(we used \GAP~\cite{GAP4}).
The algorithms and a few technical
details are outlined
in the appendices: after a brief introduction in \autoref{S.app},
the ``general'' procedures and tests (applicable to other similar problems
as well)
are explained in
\autoref{S.basicalg}, and they are applied to the two
particular tasks at hand, see
\eqref{eq.bounds}, in Appendices~\ref{S.Msigma} and \ref{S:largepencils}.

\remark	\label{remark-no-code}
Although the computation leading to \autoref{thm-main}
is
computer aided, we maintain the standards of exposition
set in other papers on rational curves on \KK-surfaces
(\eg, \cite{DIS,degt:lines,degt:singular.K3,degt:supersingular,rams.schuett,rams.schuett:char3})
and neither present nor refer to the \GAP\ scripts in the text.
We explain and motivate this standpoint in \autoref{S.app};
among other reasons,
any discussion of purely mathematical results with the help of programs
that can be executed by a particular computer algebra system will soon be
rendered obsolete by the constant evolution of scientific software.
Still, all \GAP\ scripts used in this paper can be
found as \cite{degt.Rams8.anc}. 
For the reader willing to repeat the computation with the help of~\cite{degt.Rams8.anc}
the instructions on the usage of our ``high level'' functions
are found in \verb+./code/hh/_readme.txt+;
the output is in \verb+./FOUND.txt+.
\endremark


\subsection{History of the problem}\label{s.sart}
Configurations of rational curves on surfaces
have been a
subject of intensive research
since the very beginning of algebraic geometry.
The
question what the maximal number of lines on  surfaces in a given family is
(\eg, on smooth  degree $s$ hypersurfaces in $\Cp{3}$ for a fixed integer $s > 2$)
 has a long history.
There are quite a few approaches to this question.

In the case of a hypersurface in  $\Cp{3}(\KKK)$, one can study the geometry
of the so-called \emph{flecnodal divisor}, \ie,
the locus of fourfold contact of lines with the surface in question.
This idea goes back to Salmon and Clebsch~\cite{Clebsch} and yields
the bound of at most $s(11s-24)$ lines on a smooth complex projective surface
of degree $s>2$.
Combined with certain  properties of  fibrations,
it was used in Segre's
proof~\cite{Segre} of the fact that $64$ is the maximal number of lines on
a smooth quartic in $\Cp3$. (A
gap in this proof was recently bridged
in~\cite{rams.schuett}.)
The flecnodal divisor  appears also in the proof of the sharp bound
on the number of lines on quartics in  $\PPP^{3}(\KKK)$
for algebraically closed fields~$\KKK$ such that $\mbox{char}(\KKK) \neq 2,3$
(see~\cite{rams.schuett}), the sharp bound for complex affine quartics
(see~\cite{Rams.Gonzales}), and the best known bounds for  surfaces of degree $s>4$
(see \cite{Bauer.Rams,rams.schuett:quintics}).

One can use  Bogomolov--Miyaoka--Yau's orbibundle inequality to obtain
bounds on the number of lines (more generally,
rational curves of a
bounded degree) on smooth complex
\KK-surfaces of degree $2d>4$ (see \cite{Miyaoka:lines}).

One can try to find an appropriate hyperplane section of $X$ and count the lines
that meet each of its components separately. This approach, combined with
the study of elliptic fibrations and Segre's surfaces of principal lines,
yields the sharp upper bound for projective quartics
when $\mbox{char}(\KKK) =3$ (see \cite{rams.schuett:char3}). It is
 also used in the arguments in \cite{rams.schuett,Segre,Veniani}.

In the case of a rational surface $X$
embedded \via\ a
linear system $|\hp|$, one can try to classify all solutions $v\in H_2(X)$
to the equation $v\cdot\hp =1$.  In general, this method fails for surfaces of non-negative Kodaira dimension.

There
is a very elegant approach of Elkies~\cite{Elkies},
which is
quite
efficient when the number of lines on a surface is large in
comparison with its Picard number.

Finally, one
can try to classify the potential sublattices
generated by the classes of lines  in the N\'{e}ron-Severi lattice
(\resp. second homology) of the surface in question.
In the presence of  Torelli-type theorems
this method not only leads  to sharp bounds but also provides
examples of surfaces with large line configurations.
In the case of \KK-surfaces,
this approach was pioneered in \cite{DIS}. It
gave the complete classification of smooth complex quartics with
more than $52$
lines, the sharp  bounds for $\KKK=\R$, and a bound for
$\KKK=\Q$.
Its refinements led to the sharp bound of at most $60$ lines on smooth quartics
when  $\mbox{char}(\KKK) =2$ (see~\cite{degt:supersingular}),
sharp bounds for supersingular quartics when $\mbox{char}(\KKK) =2,3$
(see~\cite{degt:supersingular}), and
sharp bounds for (smooth minimal) complex  degree $2d$ \KK-surfaces for $d>2$
or $d=1$
(see~\cite{degt:lines,degt:sextics}).
Besides, large configurations of lines are classified in
\cite{degt:supersingular,degt:lines,degt:sextics}, too.
Further generalizations of this approach resulted in sharp bounds on the
number of rational curves of a given degree  on smooth high-degree \KK-\ and Enriques surfaces
and a classification of maximal configurations (see \cite{rams.schuett:24,rams.schuett:12}).

 In contrast, until recently hardly anything was  known about line configurations on projective surfaces with non-empty singular locus.
The case of complex cubic surfaces was the only one with a complete classification
(see \cite{Bruce:cubics} for a modern exposition), whereas
for complex quartic surfaces with singularities there were
but a few partial results
(see \cite{degt:singular.K3,Jessop:quartic})
and bounds that
did not seem sharp (see \cite{Rams.Gonzales,Veniani}).
In particular, neither for
spatial complex surfaces
of degree $s > 3$ nor for complex
\KK-surfaces of
degree~$2d$ was it known
whether the maximal number of lines can be attained
by a surface with non-empty singular locus.

Recall that a projective complete intersection \KK-surface is of degree $4$, $6$ or $8$.
The configurations of lines on octic models of Kummer surfaces have a long history
(see \cite[\S\,10]{Dolgachev:book})
and a complete classification of large line configurations on  smooth \KK-octics
is found in \cite{degt:lines}.
In the present paper, we complete this picture
in the \KK-case, whereas our approach
sheds no light on the
line configurations
on ruled octic surfaces in $\Cp{5}$ (for the general classification of projective octics see
 \cite[\remcit~1.7]{Buium:surfaces}  and \cite[\S\,4.2]{Ionescu:3}).
In particular, we show that complex \KK-octics with
more than $32$
lines are always smooth.

Another application of the general theory
developed here can be found in
our recent paper \cite{Degtyarev.Rams.triangular}, where the sharp upper
bound of at most 52 lines on a complex  \KK-quartic with non-empty singular locus is
obtained.

\subsection{Common notation and conventions}\label{s.notation}
We work over the field $\C$.
Every elliptic fibration is assumed to have a section.
Otherwise, we speak of a genus-one fibration (\ie, a base-point free elliptic
pencil).

The N\'{e}ron-Severi group (\resp. Picard number) of the surface~$X$ is
denoted by $\NS(X)$ (\resp.  $\rho(X)$).
By an abuse of the language, a \emph{line} in the minimal resolution of
singularities $X\to X_{2d}$ of a quasi-polarized $K3$-surface
$X_{2d}\subset\Cp{d+1}$ is the strict transform of a straight line
in~$X_{2d}$.

When working with $K3$-surfaces, it is common to assume that $\NS(X)$ is the
``minimal'' lattice containing all relevant information, which, in our case,
consists of the quasi-polarization, lines, and, occasionally, exceptional
divisors (\ie, smooth rational curves contracted in~$X_{2d}$. Therefore, we
adopt the following convention.

\convention\label{conv.lines}
We say that the lattice $\NS(X)$ is
\emph{\qgen by lines} (\resp. lines and exceptional divisors) if it is
a finite index extension
of its sublattice generated by the classes of
lines (\resp. lines and exceptional divisors) on~$X$
\emph{and the quasi-polarization~$\hp$}.
This convention applies as well
to abstract polarized hyperbolic lattices (with a distinguished Weyl chamber
if only lines are to be considered), see \autoref{s.polarized}
and~\eqref{eq-def-Fano} for the definition of lines in this case.
In the case of special octics (\cf.
Theorems~\ref{th.hyperelliptic},~\ref{th.octics}), apart from~$h$,
the distinguished
$3$-isotropic class is also added
automatically, as part of the convention, to
the generating set.
\endconvention

As in \cite{degt:lines},
we use the following notation for common integral lattices:
\roster*
\item
$\bA_p$, $p\ge1$, $\bD_q$, $q\ge4$, $\bE_6$, $\bE_7$, $\bE_8$
are the \emph{negative definite} root
lattices generated by the indecomposable root systems of the same name
(see~\cite{Bourbaki:Lie});
\item
$[a]:=\Z u$ is the  lattice of rank~$1$ given by the condition $u^2=a$;
\item
$[a,b,c]:=\Z u+\Z v$, $u^2=a$, $u\cdot v=b$, $v^2=c$, is a lattice of
rank~$2$; when it is positive definite, we assume that $0<a\le c$ and
$0\le2b\le a$: then, $u$ is a shortest vector, $v$ is a next shortest one,
and the triple $(a,b,c)$ is unique;
\item
$\bU:=[0,1,0]$ is 
the unimodular even lattice of rank~$2$;
\item $L(n)$, denotes the lattice obtained by the scaling of
a given lattice~$L$ by a fixed integer  $n\in\Z$;
\item
$L\dual:=\Hom(L,\Z)$ denotes the dual group;
if $L$ is nondegenerate, there is a natural inclusion
$L\dual\subset L\otimes\Q$,
equipping $L\dual$ with a $\Q$-valued bilinear form,
\item
if $L$ is non-degenerate, $\discr L:=L\dual\!/L$ is the
\emph{discriminant group}, see~\cite{Nikulin:forms};
it inherits from~$L\dual$ a $\Q/2\Z$-valued quadratic form~$q$,
\item
the  inertia indices of the quadratic form
$L \otimes \R$ are denoted by $\sigma_{\pm,0}(L)$.
\endroster
In general, we maintain the standard notation for various objects
associated to a lattice (the determinant, the discriminant group, \etc.)
---see, \eg,~\cite{Conway.Sloane}, \cite{Nikulin:forms}.

\subsection{Acknowledgements}
We are grateful to Matthias Sch\"{u}tt for
a number of inspiring discussions on \KK-surfaces during our visits to
the \emph{Leibniz Universit\"{a}t}, Hannover.
This paper was conceived during Alex Degtyarev's research stay at the
\emph{Max-Planck-Institut f\"{u}r Mathematik}, Bonn, and his
short visit to
the Jagiellonian University, Krak\'{o}w; we extend our gratitude to these
institutions for their hospitality and support.
We would also like to thank the anonymous referee of the paper for
the careful reading of the manuscript and suggesting a number of improvements
of the exposition.

\section{Lattices and Fano graphs}\label{S.Fano}

All lattices considered in this paper are even:
$v^2\in2\Z$ for each $v\in L$.

\subsection{Root lattices\pdfstr{}{ \rm(see \cite{Bourbaki:Lie})}}\label{s.root}

A \emph{root lattice} is a negative definite lattice~$R$ generated by
\emph{roots}, \ie, vectors $r\in R$ of square $(-2)$. Given a negative
definite lattice~$S$, we denote by
\[*
\root_0S:=\bigl\{r\in S\bigm|r^2=-2\bigr\}
\]
the set of roots in~$S$;
then the sublattice $\rt S$ generated by $\root_0S$ is a root lattice.

Let~$R$ be a root lattice and $\Weyl$ a Weyl chamber for
(the group generated by reflections in) $R$. We denote by $\basis\Weyl$
the set of the ``outward'' roots orthogonal to the walls of~$\Weyl$.
A subset $B\subset R$ is of the form $\basis\Weyl$ if and only if it is
a ``standard'' Dynkin basis for~$R$. Recall also that a Weyl chamber~$\Weyl$
gives rise to a partition
\[*
\root_0R=\pos\Weyl\cup(-\pos\Weyl),\qquad
\pos\Weyl\cap(-\pos\Weyl)=\varnothing,
\]
with the set $\pos\Weyl$
of \emph{positive roots} \emph{closed}:
\[*
\text{if $u,v\in \pos\Weyl$ and $u+v$ is a root, then also $u+v\in \pos\Weyl$}.
\]
The positive roots $r\in \pos\Weyl$ are the
linear combinations $\sum n_ee$, $e\in\basis\Weyl$, with all $n_e\in\N$, whereas
the \emph{negative} roots $r\in{-\pos\Weyl}$ have all coefficients in $-\N$.
Conversely, any partition
\[ \label{eq-partition-roots}
\root_0R=P\cup(-P),\quad
P\cap(-P)=\varnothing,\quad
\text{$P$ is closed},
\]
defines a unique Weyl chamber $\Weyl$ such that $P=\pos\Weyl$: the
elements of~$\basis\Weyl$
are those positive roots $e\in P$ that are \emph{indecomposable}, \ie, cannot
be represented as a sum of two or more positive roots.
Needless to say that any partition as above has the form
\[*
P:=\bigl\{r\in\root_0R\bigm|\ell(r)>0\bigr\},
\]
where $\ell\:R\to\R$ is a linear functional \emph{generic} in the sense that
$\ell(r)\ne0$ for each root $r\in\root_0R$.

\subsection{Polarized lattices}\label{s.polarized}
A nondegenerate lattice~$S$ is called \emph{hyperbolic} if $\Gs_+S=1$. A
\emph{polarized lattice} $S\ni h$ is a hyperbolic lattice~$S$ equipped with a
distinguished vector $h$ of positive square; the square $h^2$ is called the
\emph{degree} of the polarization and $S$ is said to be $h^2$-polarized.
We often use the following consequence of the requirement that $\Gs_+S=1$:
\[
\text{for any pair $u,v\in S$ one has
 $\det\operatorname{Gram}(\Z h+\Z u+\Z v)\ge0$};
\label{eq.hyperbolic}
\]
moreover, the determinant is~$0$ if and only if $h,u,v$ are linearly
dependent.

Given a polarized lattice $S\ni h$, we
consider the (obviously finite) sets
\[*
\root_n(S,h):=\bigl\{r\in S\bigm|\text{$r^2=-2$, $r\cdot h=n$}\bigr\},
\quad n\in\N,
\]
and denote by $\rt(S,h)\subset h^\perp\subset S$ the sublattice generated by
$\root_0(S,h)$.
If the polarization~$h$ is understood, it is omitted from
 the notation.
We define the
\emph{positive cone}  as
\[*
\Cone^+(S,h):=\bigl\{v\in S\otimes\R\bigm|v^2{}>{}0,\ v\cdot h{}>{}0\bigl\}.
\]
We follow \cite{Vinberg:polyhedron} and call every connected
component of 
\[*
\Cone^+(S,h) \sminus \bigcup\limits_{r^2 = -2} r^{\perp}
\]
a \emph{fundamental polyhedron}. As in the case of root lattices
(see \eqref{eq-partition-roots})
each fundamental polyhedron corresponds to a partition of the set of all roots
into positive and negative ones.

Given a lattice~$S$ (negative definite or polarized), a sublattice
$F\subset S$ (negative definite or hyperbolic containing the polarization),
and a Weyl chamber $\Weyl$ for $\rt S$, we denote by $\Weyl|_F$
the Weyl chamber for $\rt F$ defined by the partition of $\root_0F$ into
$\pm P_\Weyl\cap F$.
Note that we do \emph{not} assert that
$\basis{\Weyl|_F}=\basis\Weyl\cap F$.

A fixed Weyl chamber~$\Weyl$
for $\rt(S,h)$ gives rise
to a distinguished
fundamental polyhedron $\fp$ for the group generated by
reflections of~$S$.
Here, we regard $\Weyl$ (\resp. $\fp$) as
a polyhedral (\resp. locally  polyhedral) subcone of the
cone $\Cone^+(S,h)$
and we require that $h\in\barfp\subset\bar\Weyl$ (as usual, $\bar{\ }$ standing for the closure),
 \ie, by an abuse
of  notation,
$\Weyl$
also denotes the cone
\[*
\bigl\{v\in\Cone^+(S,h)\bigm|
 \text{$v\cdot e{{}>{}}0$ for all $e\in\basis\Weyl$}\bigr\}.
\]
In
the sequel we will say that the Weyl chamber $\Weyl \subset \rt(S,h)$ \emph{extends} to the distinguished
fundamental polyhedron $\fp \subset S$ (\resp.  $\fp$ \emph{restricts} to $\Weyl$).

\remark
We
are mainly interested in the sets of walls $\basis\Weyl$ and $\basis\fp$
(see below), thus treating $\Weyl$ and $\fp$ as combinatorial objects.
(With a certain abuse of the language we also identify a wall with the
``outward'' root defining this wall.)
However, when necessary, we follow the tradition and regard them as
\emph{open} subsets of $\Cone^+(S,h)$; their closures $\bar\Weyl$ and
$\barfp$ are referred to as the \emph{closed} Weyl chamber and fundamental
polyhedron, respectively. Note that any one of the seven sets
$\Weyl$, $\bar\Weyl$, $\basis\Weyl$, $\pos\Weyl$, $\fp$, $\barfp$, $\basis\fp$
determines the six others.
We express this relation by using the same letter $\Weyl$.
\endremark

The set $\basis\fp$ of (the ``outward'' roots orthogonal to)
the walls of~$\fp$ can be found by {\bf Vinberg's
algorithm}~\cite{Vinberg:polyhedron}:
$\basis\fp=\bigcup_{n\ge0}\basis\fp_n$, where $\basis\fp_0:=\basis\Weyl$ and
the other sets are defined recursively:
\[*
\basis\fp_n:=\bigl\{r\in\root_n(S,h)\bigm|
 \text{$r\cdot e\ge0$ for all $e\in\basis\fp_k$, $0\le k<n$}\bigr\}.
\]
Denoting $\basis\fp_+:=\bigcup_{n>0}\basis\fp_n$, it is immediate that
\[
\text{$v\cdot e\ge0$ for each $e\in\pos\Weyl$ and $v\in\barfp$ or $v\in\basis\fp_+$};
\label{eq.Delta.P}
\]
in particular, two distinct vectors $v_1,v_2\in\barfp\cup\basis\fp_+$ are
never \emph{separated by a root},
\ie, $(v_1\cdot e)(v_2\cdot e)\ge0$ for each $e\in\root_0(S,h)$,
and
\[
\text{$u\cdot v\ge0$ for any two distinct vectors $u,v\in\basis\fp$}.
\label{eq.Delta.Delta}
\]
The last assertion follows directly from the construction unless
$u\cdot h=v\cdot h>0$. In the latter case, by~\eqref{eq.hyperbolic}, the
only alternative is $u\cdot v=-1$, and then $u$ and $v$
would be separated by the root $(u-v)\in\root_0(S,h)$.

The (plain) \emph{Fano graph} of a polarized lattice $(S,h)$ with a distinguished
Weyl chamber $\Weyl$ for $\rt(S,h)$ is defined as the set of vertices
\[ \label{eq-def-Fano}
\Fn_\Weyl(S,h):=\basis\fp_1=\big\{l\in\root_1(S,h)\bigm|
 \text{$l\cdot\ex\ge0$ for all $\ex\in\basis\Weyl$}\bigr\},
\]
with two vertices $l_1\ne l_2$ connected by an edge of multiplicity
$l_1\cdot l_2$, \cf.~\eqref{eq.Delta.Delta}.
On a few occasions (usually, as the ultimate result of the computation),
we also use the bi-colored \emph{extended Fano graph}
\[ \label{eq-def-exFano-lattice}
\Fex_\Weyl(S,h):=\basis\fp_1\cup\basis\Weyl,
\]
with the same convention about the multiplicities of the edges and vertices~$v$ colored
according to the value $v\cdot h\in\{0,1\}$.
Most of the time, it is the plain graphs that are used in the algorithms,
whereas their bi-colored extended counterparts play a crucial r\^{o}le at the
end of the proof and in the statements. We discuss the relation between the
two categories in \autoref{s.extendedversusplain} below.

The vertices of $\Fn_\Weyl(S,h)$ are called \emph{lines}, whereas the vectors
$\ex\in\basis\Weyl$ are called \emph{exceptional divisors}.
Due to~\eqref{eq.Delta.P}, one can also define lines as vectors
$l\in\root_1(S,h)$ such that $l\cdot e\ge0$ for all $e\in\pos\Weyl$.
As a consequence, we have the following lemma.

\lemma\label{lem.sublattice}
Let $S\ni h$ be a polarized lattice, $\Weyl$ a Weyl
chamber for $\rt (S,h)$, and
$F$ the primitive hull of the sublattice of~$S$ generated
by~$h$ and
all
$l\in\Fn_\Weyl(S,h)$. Then, there is a canonical inclusion
$\Fn_\Weyl(S,h)\subset\Fn_{\Weyl|_F}(F,h)$.
\done
\endlemma

Thus, as long as we are interested in \emph{maximizing} the number
of
lines, it suffices to consider polarized lattices $S\ni h$
\qgen by lines.

At first sight, the graphs $\Fn_\Weyl$,
$\Fex_\Weyl$
depend on the choice of the Weyl chamber $\Weyl$. However, we have the
following simple observation.

\lemma\label{lem.Weyl.independent}
For any pair $\Weyl'$, $\Weyl''$ of Weyl chambers for $\rt(S,h)$, there are
canonical isomorphisms
$\Fn_{\Weyl'}(S,h)=\Fn_{\Weyl''}(S,h)$ and
$\Fex_{\Weyl'}(S,h)=\Fex_{\Weyl''}(S,h)$.
\endlemma

\proof
Let $R:=\rt(S,h)$. Since the Weyl group~$W$ of~$R$ acts simply transitively
on the set of Weyl chambers, $\Weyl'$ and~$\Weyl''$ are related by a unique
element $\Gs\in W$. As a product of reflections, $\Gs$ admits a canonical
(identical on $R^\perp\!\ni h$) extension to an automorphism of $S\ni h$, which
induces isomorphisms of the Fano graphs.
\endproof

We conclude
with a lemma playing a crucial r\^{o}le in our algorithms.
A Weyl chamber~$\Weyl$ is called \emph{compatible} with a subset
$\graph\subset\root_1(S,h)$ if $\graph\subset\Fn_\Weyl(S,h)$.
A root $r\in\root_0(S,h)$ is called \emph{separating}
(with respect to a subset $\graph$ as above) if there is
a pair of vectors $u,v\in\graph$ such that
$r\cdot u>0$ and $r\cdot v<0$.

\lemma\label{lem.compatible}
A subset
$\graph\subset\root_1(S,h)$ admits a compatible Weyl chamber if and
only if there is no separating root $r\in\root_0(S,h)$. In this case, each
Weyl chamber $\Weyl'$ for the lattice
$S':=(\Z h+\Z\graph)^\perp\subset S$ is a face of a
\emph{unique} Weyl chamber~$\Weyl$ for $S$ compatible with~$\graph$.
\endlemma

\proof
The necessity is given by~\eqref{eq.Delta.P}, and for the sufficiency and
uniqueness we observe that the set $\pos\Weyl$
of positive roots (see \autoref{s.root})
is the set
\[*
\pos\Weyl
 =\bigl\{r\in\root_0(S,h)\bigm|r\cdot l_*+\Ge\Gf(r)>0\bigr\},\quad
\text{where $l_*:=\sum_{l\in\graph}l$},\quad 0<\Ge\ll1,
\]
and $\Gf\:S'\to\R$ is any functional positive on~$P_{\Weyl'}$
(and extended by~$0$ on the span of~$h$ and~$\graph$).
We choose~$\Ge$ so small that $\ls|\Ge\Gf(r)|<1$ for all $r\in\root_0S$;
then, the functional $\ell\:r\mapsto r\cdot l_*+\Ge\Gf(r)$
is generic and defines a partition (see \autoref{s.root}). Indeed, if
$\ell(r)=0$, then $r\cdot l_*=\Gf(r)=0$ (since $r\cdot l_*\in\Z$). The former
implies that $r\cdot l=0$ for all $l\in\graph$
(as all products $r\cdot l$ are assumed to be of the same sign);
since
also $r\cdot h=0$, we conclude that $r\in S'$, contradicting to $\Gf(r)=0$.
\endproof

\algorithm\label{alg.ext}
Given a polarized lattice $(S,h)$, subset $\Gamma$, and
Weyl chamber $\Delta'$ as in Lemma~\ref{lem.compatible},
the
unique compatible Weyl chamber~$\Weyl$ can be constructed by
Vinberg's algorithm: $\basis\Weyl=\bigcup_{n\ge0}\basis\Weyl_n$, where
$\basis\Weyl_0:=\basis{\Weyl'}$,
$l_*:=\sum_{l\in\graph}l$, and
\[*
\basis\Weyl_n=
\bigl\{\ex\in\root_0(S,h)\sminus S'\bigm|
 \text{$\ex\cdot l_*=n$ and $\ex\cdot r\ge0$ for all $r\in\basis\Weyl_k$, $k<n$}\bigr\}.
\]
Since the set
$\root_0(S,h)$
is finite, the
algorithm terminates.
(In practice, since the elements of~$\basis\Weyl$ are linearly
independent, one can as well terminate the algorithm as soon as
$\rank\rt(S,h)$
vectors have been collected.)
Note also that, once $\basis{\Weyl'}$ is known, we do not need to use the
original functional~$\Gf$ used in the proof of \autoref{lem.compatible}
or adjust constant~$\Ge$ in \latin{loc.\ cit.}
\endalgorithm


\subsection{Lines on projective \KK-surfaces} \label{s.linesonK3}
Recall that, since a \KK-surface~$X$ is simply connected, the map
$D\mapsto[D]$ factors to an isomorphism $\Pic(X)=\NS(X)$; therefore, we
freely identify classes of divisors on~$X$ with their images in $\NS(X)$.
Given an integer $d>0$, a \emph{degree $2d$ quasi-polarized \KK-surface} is
a pair $(X,\hp)$, where $X$ is a (minimal) \KK-surface and $\hp\in\NS(X)$ a
big and nef class
such that
\begin{equation} \label{eq-b-p-free}
\mbox{$\hp^2 = 2d$ and the system $\ls|\hp|$ is base-point free}.
\end{equation}

\remark \label{rem.bpfree}
Usually,
a quasi-polarized \KK-surface is defined as a pair $(X,\hp)$,
 where  $\hp$ is a big and nef line bundle on the \KK-surface $X$. Since we
 work mostly in $\NS(X)$,  we prefer to consider $\hp$ as  a class  in $\NS(X)$.

Some authors allow
the linear system $|\hp|$  of a quasi-polarization
to have fixed components.
However, since
we are
interested in the geometry of the projective
surface $f_{\hp}(X)$,
the base-point freeness is a natural assumption.

\endremark

By \eqref{eq-b-p-free}, the system $|\hp|$ defines a morphism
$f_{\hp}\:X\to\Cp{d+1}$. It is well know (see \cite[p.~615]{Saint-Donat})
that
the restricted map $f_\hp\:X\onto f_\hp(X)$
is either birational
or of degree $2$ (\latin{aka} hyperelliptic).
In the former case, we say that the quasi-polarization~$\hp$ of~$X$ is
\emph{birational}, or that $(X,\hp)$ is \emph{birationally} quasi-polarized.
In this case,
the image $X_{2d} := f_{\hp}(X)$ is a surface of degree~$2d$ with
at worst isolated singularities, all of which are Du Val,
and $f_{\hp}\: X \rightarrow X_{2d}$ is the minimal resolution of the  singularities of $X_{2d}$
(see \cite[\thmcit~6.1]{Saint-Donat}).
In the latter case (with $f_\hp$ of degree~$2$), we say that $(X,\hp)$ is
\emph{hyperelliptic}. In this case, the image
$f_{\hp}(X)$ is $\Cp{2}$, a scroll or a cone (see \cite[\propcit~5.7]{Saint-Donat}).

Recall that, by the Hodge Index Theorem, the intersection form on $V:=\NS(X) \otimes \RR$ is a non-degenerate quadratic form of index $(1, \rho(X)-1)$,
so that
the set $\{ x \in V \, | \, x^2 > 0 \}$ consists of two components, one of which
(denoted by ${\mathcal C}_X$)
contains all ample classes on $X$.
Apart from the positive cone ${\mathcal C}_X$,
the vector space $V$ contains also the \emph{nef cone} $\Nef(X)$ (\resp.
\emph{ample cone} $\Nef(X)\sminus\partial\Nef(X)$),
\ie, the set of classes in $\NS(X) \otimes \RR$
that have non-negative intersection with all curves  on $X$ (\resp. the cone spanned by the ample classes).
Moreover, by
the Riemann--Roch theorem,
\begin{equation} \label{eq-rr-for-roots}
\mbox{ every root in $\NS(X)$ is either effective or anti-effective. }
\end{equation}
In the former case we speak of an effective root.

By the adjunction formula, for each irreducible curve
$C \subset X$ with $C^{2} = -2$, we have $p_a(C)=0$, so that $C$ is smooth rational. Such curves are called $(-2)$-curves
(in particular $(-2)$-curves are always assumed to be irreducible,
\ie, they define  nodal  classes  in  $\mbox{NS}(X)$).  
We have the following well-known fact, that we will use in the sequel (see \cite[\S\,8.1]{Huybrechts}).
\lemma  \label{lemma-nakai}
The following statements hold.
\roster
\item
Let $\Ga\in\Nef(X)$. Then $\Ga \in \partial\Nef(X)$ \rom(\ie, $\Ga$ is not
ample\rom)
if and only if
either $\Ga^2=0$ or $\Ga \cdot C=0$ for a smooth rational curve $C \subset X$.
\item
Every $(-2)$-curve $C \subset X$ defines a codimension $1$ wall of $\Nef(X)$.
\item
For every $(-2)$-curve $C \subset X$ there exists a nef class $\alpha$ such that $C$ is the only smooth rational  curve in $\alpha^{\perp}$.
\done
\endroster
\endlemma

In particular,
the above lemma shows that there is a {\sl one-to-one correspondence between the set of codimension $1$  walls of the  nef cone
and the set of smooth rational curves in $X$}.
Moreover,
by \cite[\corcit~8.2.11]{Huybrechts},
the
cone
$\Nef(X) \cap {\mathcal C}_X$
is a
closed fundamental polyhedron
for the action of the
group generated by reflections
on the  cone ${\mathcal C}_X$,
\cf. $\barfp$ in \autoref{s.polarized}.

The irreducible curves contracted by the map $f_{\hp}$  are  called
\emph{exceptional divisors}.
Obviously, all exceptional divisors
are $(-2)$-curves, so
that
each exceptional divisor defines a root in  $\rt(\NS(X),\hp)$.
Indeed, the Grauert
	 contractibility criterion implies that
each exceptional divisor $C$ has negative self-intersection. Since
$p_a(C) \geq 0$,  the adjunction formula yields the equality $C^2 = -2$.
Furthermore,
the dual adjacency graph of the exceptional
divisors is a union of simply laced Dynkin diagrams,
\cf. \cite[\propcit~5.7 and \thmcit~6.1]{Saint-Donat}.

Given a quasi-polarized \KK-surface $(X, \hp)$, we call  a $(-2)$-curve $C \subset X$ a
\emph{line} on $(X,\hp)$ if $C\cdot\hp = 1$.
This definition is justified by the fact that, whenever
the quasi-polarization
is birational and  $C$ is an irreducible curve, we have
\[
\text{$C\cdot\hp = 1$ and $C^2 = -2$ if and only if
 $f_{\hp}(C)$ is a degree one curve on $X_{2d}$}.
\label{eq-lines-octic}
\]
We follow \cite{degt:lines}
 and interpret the set
\[ \label{eq-fanograph-def}
 \graphl(X,\hp) := \bigl\{
 \text{$(-2)$-curves $C\subset X$  with $C\cdot\hp = 1$}\bigr\}
\]
as a graph without loops,
where a pair of vertices $v, w\in  \graphl(X,\hp)$ is connected by a
($v \cdot w$)-fold edge.
We call this set
the \emph{\rom(plain\rom) Fano graph} of the surface
$(X,\hp)$.
As in the case of lattices, we can also consider the
bi-colored \emph{extended Fano graph}
\[  \label{eq-def-exFano-surface}
 \graphe(X,\hp) := \bigl\{
 \text{$(-2)$-curves $C\subset X$  with $C\cdot\hp \le 1$}\bigr\},
\]
with the vertices colored according to their projective degree
$C\cdot h$.
The relations between the two classes of graphs are discussed in
\autoref{s.extendedversusplain} below.

Consider the definite root lattice $\rt(\NS(X),\hp)$ and define the positive roots  as the  effective ones,
see \eqref{eq-rr-for-roots}.
We denote by $\Weyl_X$ the fundamental Weyl chamber given by this
choice of positive roots.
Then, by
\autoref{lemma-nakai},
 the nef cone $\Nef(X)$ is the distinguished
closed fundamental polyhedron $\barfp_X  \subset \NS(X)$
extending
the Weyl chamber $\Weyl_X$.
In view of \autoref{lemma-nakai},
identifying $(-2)$-curves with their classes in $\NS(X)$,
we conclude that
\[*
\graphl(X,\hp)=\Fn_{\Delta_X}(\NS(X),\hp),\quad
\graphe(X,\hp)=\Fex_{\Delta_X}(\NS(X),\hp).
\]
Then, using \autoref{lem.Weyl.independent}, we arrive at the
following statement.

\theorem \label{thm-each-chamber-ok}
Let $(X,\hp)$ be a
quasi-polarized
\KK-surface.
Then, for \emph{any} choice of the Weyl chamber $\Weyl$ for
$\rt(\NS(X),\hp)$, there are canonical isomorphisms
\[*
\graphl(X,\hp)=\Fn_{\Delta}(\NS(X),\hp),\quad
\graphe(X,\hp)=\Fex_{\Delta}(\NS(X),\hp).\done
\]
\endtheorem

\remark \label{rem.vample}
If $\hp$ is ample,
the condition on intersection with roots
in~\eqref{eq-def-Fano} is void
and finding lines reduces to solving the equations
\[  \label{eq.vample}
v^2= -2\quad  \mbox{and}\quad v\cdot\hp=1.
\]
Otherwise, \eqref{eq.vample} results in an overcount: \eg, if
$C$ is a line, and $E$ is a $(-2)$-curve such that
	$E\cdot C=1$ and $E \cdot \hp =0$, then both $C$ and $C+E$ satisfy~\eqref{eq.vample}.
\endremark


\subsection{Admissible lattices\pdfstr{}{ \rm(see \cite{Nikulin:Weil,Saint-Donat})}}\label{s.admissible}
Let $d$ be a fixed positive integer and
let $S \ni h$ be a polarized lattice of degree $2d$.
A vector $\iso\in S$ is called \emph{$m$-isotropic},
where $m = 1, 2, 3$,  if 
$$
 \iso^2=0\quad \mbox{and}\quad \iso\cdot h=m  \, .
$$
Let ~$\Weyl$ be a fixed Weyl chamber for $\rt(S,h)$.
It
extends to a unique fundamental polyhedron
$\fp\subset\Cone^+(S,h)$
such that $h\in\barfp$ (see \autoref{s.polarized}).
In order to realize the lattice $S\ni h$ as the N\'{e}ron--Severi
lattice of a quasi-polarized \KK-surface
$(X,h)$ with $\barfp=\Nef(X)$
and control the geometry of the map
$f_h\:X\to\Cp{d+1}$,
we need to impose
certain extra conditions on~$S\ni h$:
\roster
\item\label{e.h=1}
there is no $1$-isotropic vector $\iso\in S$;
\item\label{e.h=2}
for $h^2\ge4$,
there is no $2$-isotropic vector $\iso\in S$;
\item\label{e.h=3}
for $h^2=8$, there is no $3$-isotropic vector $\iso\in S$,
\endroster
As  we will explain below,  for $h^2 = 8$, there  is also  the
condition
\roster[\lastitem]
\item\label{2h}
$h\notin2S$
(so that $f_h$ does not
factor through the Veronese embedding, see \autoref{lemma-simple} below);
\endroster
however, it holds automatically whenever $\Fn_\Weyl(S,h)\ne\varnothing$.

\definition
A polarized lattice $S\ni h$ is called \emph{$m$-admissible}
(or just \emph{admissible}, if the parameter $m=1,2,3$ is understood), if
it satisfies conditions~\iref{e.h=1}--\itemrefform{m} above,
with $h^2\ge4$ if $m=2$ and $h^2=8$ if $m=3$.
\enddefinition

The admissibility of a polarized lattice is easily established in the algorithm
computing the Fano graph
(see \autoref{ss.extensibility} below or
\cite[\algcit~2.5]{degt:lines}).
By definition,  the
notion of $m$-admissibility is independent of the choice of
a Weyl chamber~$\Weyl$ for $\rt(S,h)$ (\ie,
the choice  of a fundamental polyhedron $\fp$). Still we
have the following observation, that will be useful in the sequel.

\observation \label{obs-imposing}
Let $S\ni h$ be an $(m-1)$-admissible lattice  with  $h^2$ in the range of applicability of the corresponding condition. If $S$ fails to be
  $m$-admissible, then there exists  an $m$-isotropic vector  $\iso\in\barfp$.
\endobservation

\proof
Let $\iso$ be an $m$-isotropic vector in $S$.
Since the Weyl group of $\rt(S,h)$ is finite,
$\bar\Weyl\subset\bar\Cone^+(S,h)$ is a fundamental domain of the (extended
to~$\bar\Cone^+$) action of this group.
Hence,
after applying a sequence of reflections (which all
preserve~$h$), we can assume that $\iso\in\bar\Weyl$.
(Geometrically, if $S=\NS(X)$ and $\Weyl=\Weyl_X$, this procedure
corresponds to the passage to the moving part of the linear
system~$\ls|p|$.)

We assert that then
immediately $\iso\in\barfp$. Indeed, let $r\in\basis\fp_n$, $n\ge1$, be a
wall
such that $\iso\cdot r=-a\le-1$. By~\eqref{eq.hyperbolic}, this implies
\[
h^2\le\frac{2m(m-an)}{a^2},
\label{eq.a=e.r}
\]
and the cases to exclude are
\roster*
\item
none, if $m=1$, as the inequality would imply $h^2\le0$;
\item
$h^2=4$, \cf.~\iref{e.h=2}, and $n=a=1$ if $m=2$;
\item
$h^2=8$, \cf.~\iref{e.h=3}, and $n=a=1$ if $m=3$ (here $n \neq 2$ results from $h^2=8$).
\endroster
In the last two cases,
we have $(\iso-r)^2=0$ and $(\iso-r)\cdot h=m-1$,
\ie, $(S,h)$ is not $(m-1)$-admissible and the extra restriction makes no
sense.
\endproof

We will mainly consider $2$-admissible lattices, (\ie, those corresponding
to the birational
projective models of \KK-surfaces, see \autoref{lemma-simple} below).

\lemma\label{lem.simply-laced}
Let $S\ni h$ be a polarized lattice, $\Weyl$ a
Weyl chamber,
$\ex_1,\ex_2\in\basis\Weyl$ distinct exceptional divisors, and
$l_1,l_2\in\Fn_\Weyl(S,h)$ distinct lines. Then\rom:
\roster
\item
if $S\ni h$ is $2$-admissible, then $l_1\cdot l_2\in\{0,1\}$\rom;
\item
if $S\ni h$ is $1$-admissible, then $l_1\cdot\ex_1\in\{0,1\}$\rom;
\item
for all lattices, $\ex_1\cdot\ex_2\in\{0,1\}$.
\endroster
\endlemma

\proof
By~\eqref{eq.hyperbolic} and~\eqref{eq.Delta.Delta}, we have
$0\le l_1\cdot l_2\le2(h^2+1)/h^2$ and, since $h^2\ge4$, the maximum is
$l_1\cdot l_2=2$. However, in this case $\iso:=l_1+l_2$ is a $2$-isotropic vector
and $(S,h)$ is not $2$-admissible.
Similarly, $0\le l_1\cdot\ex_1\le2$ and, in the case $l_1\cdot\ex_1=2$, the
vector $\iso:=l_1+\ex_1$ is $1$-isotropic. The last assertion is a
well-known property of root systems with all roots of square~$(-2)$.
\endproof

To justify our interest in  $m$-isotropic vectors,
let us consider
the polarized lattice $\NS(X) \ni \hp$ for a
\KK-surface  $X$
and a big and nef class $\hp \in \NS(X)$.
 As before, we fix the
Weyl chamber $\Weyl_X$  given by the
effective roots, \cf. \eqref{eq-rr-for-roots}.
Recall that
the closure of its extension to $\bar\Cone_X$ is the nef cone $\Nef(X)$.

Let us assume that  $m \in \{1, 2, 3\}$, and (for $m > 1$)  the lattice  $\NS(X)$ is $(m-1)$-admissible. We claim that
\begin{equation}  \label{eq-m-isotropic}
\NS(X) \mbox{ is $m$-admissible iff } E\cdot\hp > m \,
 \mbox{ for every elliptic curve } E \subset X \, .
\end{equation}
Indeed, assume that $\NS(X)$ fails to be  $m$-admissible.
Then, by \autoref{obs-imposing}, the lattice $\NS(X)$ contains a nef
$m$-isotropic class~$p$. It is well-known that
$\ls|p|$ is a base-point free elliptic pencil on~$X$, and one can take for~$E$ any smooth
fiber of~$\ls|p|$. The converse statement is obvious: take
$p=[E]$.

The above discussion,
 combined with  \cite{Saint-Donat},
 yields the following well-known lemma,
which we state for the reader's convenience.

\lemma \label{lemma-simple}  \label{remark-sd}
Let $X$ be a \KK-surface and let $\hp \in \NS(X)$ be big and nef.
  Then\rom:
\roster
\item\label{1-admissible}
the linear system $|\hp|$ is base-point-free 
 if and only if the polarized lattice  $\NS(X)  \ni \hp$
 is $1$-admissible\rom;
\item\label{4.birational}
assuming that $\hp^2\ge4$, $\hp^2 \neq 8$, and
 $\NS(X)   \ni \hp$ is  $1$-admissible, the map  $f_\hp$ is birational
\rom(onto its image\rom)
if and only if
  $\NS(X)   \ni \hp$ is $2$-admissible\rom;
\item\label{8.birational}
assuming that $\hp^2=8$ and $\NS(X) \ni \hp$ is
$1$-admissible, the map  $f_\hp$ is birational
if and only if
$\NS(X)  \ni \hp$ is  $2$-admissible and $\hp\notin 2\NS(X)$\rom;
\item \label{triquadric}
assuming that $\hp^2=8$ and $f_\hp$ is birational,
\cf.~\iref{8.birational},
the image $X_8 := f_\hp(X)$ is an
intersection of three quadrics if and only if
$\NS(X) \ni \hp$ is $3$-admissible.
 \endroster
\endlemma
\proof
By \cite{Nikulin:Weil} (see also \cite[the proof of \propcit~8.1]{Saint-Donat}) the linear system $|\hp|$ is base-point free   if and only if the surface $X$ contains no
 irreducible
 genus-one curve $E$ such that $E \cdot \hp=1$. Thus, Statement~\iref{1-admissible}
 follows from
 \eqref{eq-m-isotropic}.

Statements~\iref{4.birational} and~\iref{8.birational} (\resp. Statement~\iref{triquadric}) follow from \cite[\thmcit~5.2]{Saint-Donat} (\resp.  \cite[\thmcit~7.2]{Saint-Donat})
and  \eqref{eq-m-isotropic}.
\endproof


\definition \label{def.triquadric}
Let $(X, \hp)$ be a degree $8$ birationally quasi-polarized \KK-surface.
We say that $(X, \hp)$ (or the image $X_8 :=f_{\hp}(X)$) is a
	\emph{triquadric} (\resp. \emph{a special octic})
if it satisfies the equivalent conditions of
\autoref{lemma-simple}\iref{triquadric} (\resp. otherwise).
\enddefinition


\subsection{Hyperelliptic models and special octics\pdfstr{}{ \rm(see~\cite{degt:lines})}}\label{s.hyperbolic}
Given a quasi-polarized
\KK -surface $(X,\hp)$
it is natural to consider the situation when the map $f_h$
degenerates (\ie, the hyperelliptic case)
or its image is a special octic.
In lattice-theoretic terms one  has to study
a polarized lattice  $S\ni h$, with a
fixed Weyl chamber~$\Weyl$ for $\rt(S,h)$,
that
 fails to be $m$-admissible (see  \autoref{remark-sd}).
Moreover,
in view of
 \autoref{obs-imposing},
such a degeneration amounts to
 \emph{imposing} the
existence of an $m$-isotropic vector $\iso\in\barfp$ while
assuming that $S\ni h$ is $(m-1)$-admissible for $m=2,3$.
We have the following lemma,
slightly different from the smooth case. Recall that,
for $m=3$,
our principal interest is the case $h^2=8$; nevertheless,
we state the lemma in full generality
for the sake of completeness.

\lemma\label{lem.e}
Assume that the polarized lattice $S\ni h$ is $(m-1)$-admissible and
$\iso\in\barfp$ is an $m$-isotropic vector, $m=2,3$.
Then\rom:
\roster
\item\label{P1xP1}
if $h^2=4$ and $m=2$, there are at most two $2$-isotropic
vectors $\iso,\bar\iso\in\barfp$ and, if there are two such vectors
\rom(\ie,  $\iso \neq\bar\iso$\rom), then $\iso+\bar\iso=h$\rom;
\item\label{hyperelliptic}
otherwise, $\iso$ is the unique $m$-isotropic vector in $\barfp$.
\endroster
Furthermore, if either $m=2$ and $h^2\ge4$ or $m=3$ and
$h^2\ge6$, then
one has $l\cdot\iso\in\{0,1\}$ for each line
$l\in\Fn_\Weyl(S,h)$ and $\ex\cdot\iso\in\{0,1\}$ for each exceptional
divisor $\ex\in\basis\Weyl$.
\endlemma

\proof
Let $\bar\iso$ be an $m$-isotropic vector,
 $\bar\iso  \neq \iso$.
Then, by \eqref{eq.hyperbolic},
$0<\iso\cdot\bar\iso\le2m^2\!/h^2$,
implying, with two exceptions, that $\iso\cdot\bar\iso=1$
and, hence, $\iso$ and $\bar\iso$
are separated by the root $(\iso-\bar\iso)$, so
that $\bar\iso \notin \barfp$, see \autoref{s.polarized}.

Exceptionally, one may have $\iso\cdot\bar\iso=2$ for $h^2=4$, $m=2$ or $h^2=8$,
$m=3$. In the latter case, $(h-\iso-\bar\iso)$ is a $2$-isotropic vector,
so that this case is excluded and the proof of
Statement~\eqref{hyperelliptic} is complete.

For $h^2=4$, $m=2$,
the linear dependence given by~\eqref{eq.hyperbolic} is
 $\iso+\bar\iso=h$ and both
vectors are in~$\barfp$ \emph{unless they are separated by a root}.
This
proves Statement~\eqref{P1xP1}.

For the last statement,
one has $0\le l\cdot\iso$ by~\eqref{eq.Delta.P}.
If $m=2$, then~\eqref{eq.a=e.r} implies that either
$l\cdot\iso\le1$ or  $h^2=4$ and $l\cdot\iso=2$.
However, in this exceptional case, $(h-\iso-l)$ is a $1$-isotropic vector,
contradicting to the assumption that $S$ should be $1$-admissible.
For $m=3$ and $l\cdot\iso=2$
(resp. $l\cdot\iso=3$) we get  $h^2 \leq 6$
(resp. $h^2 = 4$) from ~\eqref{eq.Delta.P}.
The exceptional case $l\cdot\iso=2$ and $h^2 = 6$ yields  the $2$-isotropic vector
$(h-\iso-l)$, so it is excluded.

Similarly, the inequality  $\ex\cdot\iso\le1$ follows
from~\eqref{eq.a=e.r}.
\endproof

The unique vector $\iso\in\barfp$ (or pair $\iso,\bar\iso\in\barfp$)
given by \autoref{lem.e} constitutes an extra
structure on the lattice.
For example, there is a canonical partition
\[
\Fn_\Weyl(S,h)=C_0\cup C_1,\qquad
C_i:=\bigl\{l\bigm|l\cdot\iso=i\bigr\},\quad i=0,1.
\label{eq.partition}
\]
(In case~\iref{P1xP1}, the numbering of~$C_i$
depends on the choice of one of the two vectors.)
Conversely, if we assume that $S$ is rationally generated
by~$h$, $\iso$, and all lines (\cf. \autoref{lem.sublattice} and remark thereafter),
the class $\iso\in\barfp$ (or the pair $\iso,\bar\iso\in\barfp$),
if it exists, is uniquely
recovered from the above partition.

\remark\label{rem.partitioned}
Although graph-theoretically the two notions are identical, we will distinguish
between \emph{partitioned} (\via\ $v\mapsto v\cdot\iso$) and \emph{bi-colored}
(\via\ $v\mapsto v\cdot h$) graphs, reserving the two words as
separate terms.
Occasionally, we even consider bi-colored partitioned graphs, so that each
\emph{part}~$C_i$, $i=0,1$, is bi-colored.
\endremark

\remark
Case~\iref{P1xP1} of \autoref{lem.e} is of limited interest when
dealing with
line configurations on \KK-surfaces. Indeed, let us consider
a quasi-polarized \KK-surface $(X,\hp)$ of degree~$4$
and assume that
$\Nef(X)$ contains two $2$-isotropic vectors.
Each nef $2$-isotropic  class defines a genus-one fibration on~$X$.
Thus, $X$ is endowed with two 
genus-one fibrations
and, by~\eqref{eq.partition},
each line on~$X$ is a component of a fiber of
exactly one of them. 
This gives
us an
obvious
upper bound of $48$ lines
(recall \cite[Proposition III.11.4]{Barth.Peters.VanDeVen}) and the equality $e(X)=24$),
which is sharp for smooth models (see~\cite{degt:singular.K3}),
extends to all models, and fails to be sharp (by at least two)
when $\hp$ is not  ample.

The other cases with $m=2$ are not interesting either:
essentially, we are speaking about a single genus-one fibration. 
\endremark

\subsection{Miscellaneous definitions}\label{s.misc}
For
the reader's convenience,
we collect below
a few common definitions used in the sequel.

A graph $\graph$ is called \emph{discrete} if it has no edges.

The
\emph{girth} $\girth(\graph)$ of a graph~$\graph$ is the length of a
shortest cycle in~$\graph$
(with the convention that the girth of a forest is~$\infty$),
and the \emph{independence number} $\nodal(\graph)$
is the cardinality of the largest independent vertex set (\latin{aka}
discrete induced subgraph).
Given a bi-colored graph~$\graph'$,
we denote by $\spp_c\graph'$
the \emph{plain} subgraph of~$\graph'$ induced by the vertices of~$\graph$
of the chosen color $c\in\{0,1\}$.

\definition  \label{def.pseudo}
Let $S\ni h$ be a $1$-admissible lattice and $\Weyl$ a distinguished Weyl
chamber.
A \emph{pseudo-vertex} of the Fano graph $\graph:=\Fn_\Weyl(S,h)$
is either
an exceptional divisor
$\ex\in\basis\Weyl$ or a $2$- or $3$-isotropic vector $\iso\in\fp$.

Given a subgraph $\pencil\subset\graph$, we define
the \emph{support} of a
\hbox{(pseudo-)}\penalty0vertex $v$ of~$\graph$
(relative to the subgraph $\pencil$)
as the set
\[*
\supp_\pencil v=\|v\|_\pencil:=\bigl\{l\in\pencil\bigm|l\cdot v=1\bigr\}
\]
(with se second notation to be used in crowded formulas).

 As usual, if the subgraph $\pencil$ is understood, the subscript
 in $\supp_\pencil v$
 is omitted. The support of a
 (pseudo-)vertex relative to the whole graph is called the \emph{\sstar}
 \[*
 \Star(v):= \supp_\graph v.
 \]
\enddefinition

We use these notions in the discussion of the local structure of Fano
graphs (\autoref{s.star}) and in various algorithms,
where we \emph{identify} an extra vertex of a graph extension with its support
(see \eg~\autoref{s.ext}).

\section{Abstract graphs}\label{S.graphs}

This section is devoted to the question
whether a given graph is (isomorphic to)  the
Fano graph of a quasi-polarized \KK-surface.
We start with a discussion of
various properties of polarized lattices
spanned
by lines.

\subsection{Extensible graphs}\label{s.extension}
Given
a graph~$\graph$, let $\Z\graph$ be the lattice freely generated by the
vertices $v\in\graph$, $v^2=-2$, with $u\cdot v=n$ whenever $u$ and~$v$ are
connected by an $n$-fold edge.
Then, we associate to~$\graph$ the polarized lattice
\[ \label{eq:fanolattice}
\Fano_{2d}(\graph):=(\Z\graph+\Z h)/\ker,\qquad h^2=2d>0,\quad
\text{$h\cdot v=1$ for $v\in\graph$},
\]
where $\ker:=\ker(\Z\graph+\Z h)=(\Z\graph+\Z h)^\perp$
stands for the kernel of the bilinear form.
Usually, the even integer $2d$ is fixed in advance and omitted from the
notation.
We speak of {\sl $2d$-polarized} graph $\graph$ when it is needed to avoid ambiguity.

In a similar manner, we define the lattice
$\Fano_{2d}(\graph')$ for a  bi-colored graph
$\graph'$;
the  last condition in  \eqref{eq:fanolattice} is replaced by the equality
$h \cdot v = c(v)$,  where
 $c(v) \in \{0,1\}$ is  the color of the vertex $v \in \graph'$.

We treat vertices of~$\graph$ as vectors in $\Fano(\graph)$ and freely use
mixed terminology:
\roster*
\item
if $v_1$, $v_2$ are adjacent in~$\graph$, we say that they \emph{intersect} or
$v_1\cdot v_2=1$;
\item
otherwise, we say that $v_1$, $v_2$ are \emph{disjoint} or $v_1\cdot v_2=0$;
\item
vertices are (linearly) independent if so are the corresponding vectors.
\endroster
We also apply to~$\graph$ other lattice theoretic terminology; \eg, we speak
about the \emph{rank}
$\rank\graph:=\rank_{2d}\graph$ or \emph{inertia indices}
$\Gs_\pm\graph:=\Gs_{2d,\pm}\graph$,
referring to the lattice $\Fano(\graph)$. Note, though, that we only consider
graphs with $\Gs_+\graph=1$, \ie, \emph{we always assume that $\Fano(\graph)$ is
hyperbolic}.

\remark\label{rem.hyperbolic}
The assumption that $\Fano_{2d}(\graph)$ should be hyperbolic, which is always
checked first (\cf., \eg, \autoref{lem.Sylvester} below), is already quite
strong. The inertia indices of $\Z\graph$ are the subject of the so-called
spectral graph theory, and almost no general results
(beyond the classical fact that Dynkin diagrams, elliptic or affine,
are precisely the graphs with $\Gs_+=0$) are known. Even if $\Z\graph$ is
hyperbolic, $\Gs_+\Fano_{2d}(\graph)$ depends on~$d$ and eventually
becomes~$2$, see~\cite{degt:lines}.
\endremark

Occasionally, we also pick an isotropic subgroup
$\CK\subset\discr\Fano(\graph)$ and consider the extension
$\Fano(\graph,\CK):=\Fano_{2d}(\graph,\CK)$ of $\Fano(\graph)$ by~$\CK$;
in this notation,
$\Fano(\graph)$ is an abbreviation for $\Fano(\graph,0)$.
(Recall
that the isomorphism classes of
even finite index extensions of a non-degenerate even lattice~$L$
are in a one-to-one correspondence with
the $q$-isotropic subgroups of $\discr L$, see~\cite{Nikulin:forms}.)

\definition\label{def.ext}
As in \autoref{s.hyperbolic}, a Weyl chamber~$\Weyl$ for
$\rt(\Fano(\graph,\CK),h)$ is called \emph{compatible} with~$\graph$
(assuming a fixed kernel~$\CK$) if $\graph\subset\Fn_\Weyl\Fano(\graph,\CK)$.
The graph~$\graph$ or pair $(\graph,\CK)$ is called \emph{extensible}
(in degree~$2d$) if it admits a compatible Weyl chamber.
\enddefinition

The concept of extensibility
replaces the requirement $\root_0\Fano(\graph,\CK)=\varnothing$ used in the
smooth case to rule out the majority of ``bad'' graphs. It is effective due
to a simple criterion given by the next statement and heredity given by
\autoref{cor.ext}.



\lemma[a corollary of \autoref{lem.compatible}]\label{lem.ext}
A pair $(\graph,\CK)$ is extensible if and only if
the lattice $\Fano(\graph,K)$ has no
separating roots with respect to $\graph$.
If this is the case, there is a \emph{unique} Weyl chamber~$\Weyl$ for
$\rt(\Fano(\graph,\CK),h)$ compatible with~$\graph$.
\done
\endlemma



Given a pair $(\graph,\CK)$ and an induced subgraph $\graph'\subset\graph$,
denote
\[*
\CK|_{\graph'}:=\bigl(\Fano(\graph')\otimes\Q\bigr)\cap\Fano(\graph,\CK)\bmod\Fano(\graph')
 \subset\discr\Fano(\graph'),
\]
so that $\Fano(\graph',\CK|_{\graph'})\subset\Fano(\graph,\CK)$ is the
primitive sublattice rationally generated by~$h$ and the vertices
of~$\graph'$.
We say that a pair $(\graph',\CK')$ is \emph{subordinate} to $(\graph,\CK)$,
denoted $(\graph',\CK')\prec(\graph,\CK)$,
if $\graph'\subset\graph$ is an induced subgraph
and $\CK'\subset\CK|_{\graph'}$.

\corollary[of \autoref{lem.ext}]\label{cor.ext}
Extensibility is a hereditary property\rom:
if $(\graph,\CK)$ is extensible, then
so is any subordinate pair $(\graph',\CK')\prec(\graph,\CK)$.
\done
\endcorollary

As another consequence, an extensible pair $(\graph,\CK)$ has a well-defined
\emph{saturation} and \emph{extended saturation}
\[*
\sat_{2d}(\graph,\CK):=\Fn_\Weyl\Fano_{2d}(\graph,\CK),\qquad
\satex_{2d}(\graph,\CK):=\Fex_\Weyl\Fano_{2d}(\graph,\CK),
\]
where $\Weyl\subset\rt\Fano(\graph,\CK)$ is the
unique Weyl chamber compatible
with~$\graph$.
As with the other notation, we usually omit the subscript~$2d$.

\warning\label{wrn.sat}
\latin{A priori}, unlike the case of smooth models,
we cannot claim that the inequality
$(\graph',\CK')\prec(\graph,\CK)$ implies the inclusion
$\sat(\graph',\CK')\subset\sat(\graph,\CK)$, as the smaller lattice may
contain extra lines that are inhibited by some extra exceptional divisors in
the larger one.
As a simple example, we have
\[*
\sat(\tA_2\oplus\bA_2)=2\tA_2\not\subset
\tA_2\oplus\bA_3=\sat(\tA_2\oplus\bA_3).
\]
Indeed, denoting by $\kappa:=e_1'+e_2'+e_3'$ the fundamental cycle of the
first summand $\tA_2$ and by $e_i$, the basis elements in the second
summand, we have that
\roster*
\item
$l:=\kappa-e_1-e_2$ must be a line in $\sat(\tA_2\oplus\bA_2)$, whereas
\item
$e:=\kappa-e_1-e_2-e_3$ is an exceptional divisor in $\sat(\tA_2\oplus\bA_3)$;
\endroster
even though $l$ is still in the second lattice, it is not a line as
$l\cdot e<0$.
Even more striking
\autoref{rem.kum.sing} below shows that the number of
lines may decrease when the lattice is extended.
This renders the study of surfaces with singular points more difficult.
For example,
any validity criteria based on the absence
of extra lines of a certain kind (\eg, the type of the full graph)
are \emph{not} hereditary and
must be avoided.
Besides, \emph{until the very end of the computation}, we cannot try to
reduce the
overcounting by using any sorting based on the full
configuration $\sat(\graph,\CK)$: only the lines explicitly contained
in~$\graph$ can be taken into account. We give more details of these
technical problems in Appendices~\ref{S.basicalg}, \ref{S.Msigma},~\ref{S:largepencils}.
\endwarning

\subsection{Admissible and geometric graphs}\label{s.aggraphs}

We call a graph ~$\graph$ (\resp. pair $(\graph,\CK)$)  \emph{$m$-admissible}
(in a fixed degree~$2d$)
if it is extensible and
the lattice $\Fano_{2d}(\graph)$ (\resp. $\Fano_{2d}(\graph,\CK)$) is
$m$-admissible, where $m = 1, 2, 3$.
Clearly, this
property is hereditary,
\cf. \autoref{cor.ext}.

Obviously, if  a given configuration can be realized as the Fano graph of a
\KK-surface $(X,\hp)$,
there must exist
a primitive isometric embedding of $\Fano(\graph,\CK)$  into  the \KK-lattice
 $H_2(X;\Z)$.
That is why we introduce the following notion.

\definition\label{def.madmissible}
Let $d\ge1$ be a fixed integer and $m \in \{ 1, 2, 3\}$.
\roster
\item
A $2d$-polarized hyperbolic lattice $S\ni h$ is called \emph{$m$-geometric}
if it is $m$-admissible and there exists a primitive isometric embedding
\[*
S\into\bL:=2\bE_8\oplus3\bU.
\]
\item
Given a \emph{plain}
graph $\graph$, an isotropic subgroup
$\CK\subset\discr\Fano_{2d}(\graph)$
is an \emph{$m$-geometric kernel} (in degree~$2d$) if the lattice
$\Fano_{2d}(\graph,\CK)$ is $m$-geometric. We put
\[*
\GK^{m}_{2d}(\graph)  := \bigl\{ \CK\subset\discr\Fano_{2d}(\graph)\bigm|
 \text{$\CK$ is m-geometric} \bigr\}.
\]
\item
A \emph{plain} graph $\graph$ is \emph{$m$-subgeometric} (in degree~$2d$)
if $\GK^{m}(\graph)$  is
non-empty.
\item\label{graph.geometric}
An $m$-subgeometric graph $\graph$ is called  \emph{$m$-geometric} (in degree~$2d$) if
\[*
\graph\cong\sat_{2d}(\graph',\CK)
\]
for a graph $\graph'\subset\graph$ and kernel $\CK \in \GK^{m}_{2d}(\graph')$.
\item
A \emph{bi-colored} graph~$\graph'$ is
\emph{$m$-geometric} (in degree~$2d$) if
\[*
\graph'=\Fex_\Weyl(S,h)
\]
for some $m$-geometric $2d$-polarized lattice $S\ni h$.
\endroster
Whenever this leads to no ambiguity, we omit the prefix ``$m$'' and/or
degree~$2d$ and speak about (sub-)geometric graphs \etc.
\enddefinition

\remark\label{rem.heredity}
The fact that
$\Fano(\graph,\CK)$ is  $m$-geometric
is \emph{not} a hereditary property of pairs $(\graph,\CK)$
as one can loose the primitivity by passing to a smaller subgroup
$\CK'\subset\CK$.
In contrast, the \emph{existence} of an $m$-geometric kernel is a
hereditary property of graphs: if $\CK\in\GK^m(\graph)$, then
$\CK|_{\graph'}\in\GK^m(\graph')$ for any induced subgraph
$\graph'\subset\graph$.
\endremark

Nikulin's results  \cite{Nikulin:forms} give us precise criteria
for
a lattice
to admit a primitive embedding
into the \KK-lattice $\bL$ (\cf. \cite[\thmcit~3.2]{DIS}).
Thus,
$m$-(sub-) geometric graphs are easily detected
using \GAP~\cite{GAP4} (see~\autoref{ss.geometric} below).

As an immediate consequence of \autoref{def.madmissible}
and general theory of  lattice-polarized \KK-surfaces
(Nikulin~\cite{Nikulin:finite.groups},
Saint-Donat~\cite{Saint-Donat}; \cf. also
\cite[\thmcit~3.11]{DIS}
and \cite[\thmcit~7.3]{degt:singular.K3}),
we have the following theorem.

\theorem\label{th.K3}
Let
$d\ge2$ be a fixed integer. A graph $\graph$ is $2$-geometric in
degree~$2d$ if and only if $\graph\cong\graphl(X,\hp)$ for a
degree~$2d$ birationally quasi-polarized \KK-surface $(X,\hp)$
such that $\NS(X)$ is \qgen by lines.
If $2d=8$, then $\graph$ is $3$-geometric if and only if $(X,\hp)$ above can
be chosen to be a triquadric.
\done
\endtheorem

\autoref{th.K3} is a ``classical'' statement dealing with lines only.
However, in view of \autoref{wrn.sat}, it may fail to describe \emph{all}
configurations of lines:
unlike the smooth case,
we cannot
assert
that $\graph\subset\graphl(X,\hp)$
if the extension $\NS(X)\supset\Fano(\graph,\CK)$ is
of positive corank.
We address this phenomenon in the next section; for now, we make a precise
statement, taking into account the exceptional divisors.

\theorem\label{th.K3.ext}
A bi-colored graph~$\graph'$
is $2$-geometric in
degree~$2d$ if and only if $\graph'\cong\graphe(X,\hp)$ for a
degree~$2d$ birationally quasi-polarized \KK-surface $(X,\hp)$.
It is $3$-geometric in degree~$8$ if and only if $(X,\hp)$ can
be chosen to be a triquadric.
\done
\endtheorem

\remark
As in the smooth case  (\cf. \cite[\thmcit~2.2]{degt:lines})
the general theory
gives us
the dimension of the family of \KK-surfaces that realize a given
graph~$\graph$:
modulo the projective group, it equals $20-\rank_{2d}\graph$.

Note that,
for  $d \neq 4$, the geometric interpretation
of the fact that a bi-colored graph is $3$-geometric can be derived from
\cite[\thmcit~7.2]{Saint-Donat}. We omit the discussion of the general case
to maintain our exposition compact.
\endremark


\subsection{Extended vs\. plain Fano graphs}\label{s.extendedversusplain}
At this point,
we need to emphasize the fact that we thoroughly
distinguish between two categories:
\roster*
\item
plain graphs---typically, the input for $\Fano(-)$ and output of~$\Fn(-)$, and
\item
bi-colored graphs---the output of $\Fex(-)$ (even if all colors are~$1$).
\endroster
(Occasionally, as in \autoref{s.partitioned} below, we also consider
partitioned graphs, both plain and bi-colored, \cf.
\autoref{rem.partitioned}.) Apart from the fact that we are interested in
large configurations of \emph{lines} (plain graphs) but strive to provide
more complete geometric information (bi-colored extended graphs),
there is a deeper reason that makes the
presence of both categories essential for the paper:
\roster*
\item
for plain graphs (regarded as graphs of lines),
we have a simple hereditary extensibility criterion
given by \autoref{lem.ext},
whereas
\item
for extended graphs, we have a complete geometric realizability criterion given by
\autoref{th.K3.ext} (\vs. \autoref{th.K3}, which needs extra
assumptions).
\endroster
The difficulty is that \autoref{th.K3.ext} applies to \emph{saturated}
bi-colored graphs only, but, unless $\spp_0\graph'=\varnothing$, we have no
means of detecting whether a given bi-colored graph $\graph'$ is saturable
and, hence, almost no means of constructing saturated bi-colored graphs other than
saturating (\via\ $\satex$) \emph{plain} ones.
For this reason, we do not introduce the notion of
subgeometric or even admissible bi-colored graph.

One may try to relate the two categories by functors
\[*
\graph'\mapsto\spp_1\graph'\ \text{(bi-colored $\to$ plain)}\quad\text{and}\quad
 \graph\mapsto\graph\extended\ \text{(plain $\to$ bi-colored)}
\]
(for the latter, see \autoref{eq-graphex} below), but this does not
work quite well: the latter, if properly defined, turns out to be multi-valued,
whereas the former does not
need to take geometric graphs to geometric (see \autoref{rem.kum.sing}
below; in general, even if $\graph'$ itself is geometric, we can only assert that
$\spp_1\graph'$ is \emph{subgeometric}).

\remark\label{rem.geometric.term}
In other words,
unless $\NS(X)$ is \qgen by lines, the plain graph $\graphl(X,h)$ does not need to be geometric in the sense of
\autoref{def.madmissible}. Still our
definition of the term
``geometric'' for plain graphs
is both consistent with the existing literature on the smooth case
and convenient 
 (it is
used in most algorithms below). Consequently, whenever we speak of a geometric plain graph, we mean \emph{geometric in the sense of \autoref{def.madmissible}}. Otherwise (when  $\graph \cong \graphl(X,h)$ for a \KK-surface $(X,h)$), we informally use the term
	``geometrically realizable''.
\endremark

\convention\label{eq-graphex}
Recall that, given an
$m$-geometric plain graph $\graph$
and a kernel $\CK \in \GK^{m}_{2d}(\graph)$, one can
 define the bi-colored \emph{extended graph} $\satex_{2d}(\graph,\CK)$.
However, this extended graph is not uniquely determined by $\graph$
alone (see, \eg,
 the extended Fano graphs $\config{QF.30.5''}$ and $\config{QF.30''}$ in
 \autoref{s.pent}),
 and we are using the vague notation $\graphex$ to refer to \emph{any} of
 these extensions (usually, when listing the exceptions).
 \endconvention

In view of these subtleties, the classification of large Fano graphs,
\emph{even plain},
is more involved than in the smooth case: it is no longer enough to work with
abstract graphs and lattices of the form $\Fano(-)$.
Here, by ``large'' we still mean graphs with
many \emph{lines}, \ie, we fix a threshold~$M$ and try to list all
geometrically realizable plain  (resp. geometric bi-colored) Fano graphs
$\graph^*=\Fn^*(X,\hp)$
satisfying $\ls|\graph|\ge M$
(\resp. $\ls|\spp_1\graphex|\ge M$).
To localize the problem, we break the proofs into two stages: the first
(the bulk of the proof) runs
essentially as in the smooth case, whereas the second, new one boils down to
a more thorough analysis of but a few extremal graphs.

\subsubsection{Stage~$1$\rom: lattices \qgen by lines}\label{ss.stage.1}
We
list all geometric
 plain
graphs~$\graph$ satisfying $\ls|\graph|\ge M$.
The algorithms are essentially those used in the smooth case, with a
few minor modifications (\cf. \autoref{wrn.sat}) and \autoref{lem.ext}
replacing the requirement that
$\root_0\Fano(\graph,\CK)=\varnothing$.
The result of this ``classical'' part of the proof is the sharp upper bound on the number of
lines (see \autoref{lem.sublattice} and remark thereafter)
and some (but not all) geometrically realizable large  configurations of lines.

\subsubsection{Stage~$2$\rom: the general case}\label{ss.stage.2}
To complete the \emph{classification}, we still need to
consider
the polarized lattices $S\ni h$ that fail to be  \qgen by lines. (For
example, the configuration of $16$ pairwise disjoint lines on the
elements of family~$\Kummer_{64}$ would be missing after Stage~$1$, see
\autoref{rem.kum.sing} below.) To this end, we need to consider all geometric
extensions $\bS\supset S$ of all geometric lattices $S:=\Fano(\graph,\CK)$,
$\ls|\graph|\ge M$,
found in Stage~$1$.
We may assume that
\roster
\item\label{stage.2.rank}
$\rank\bS>\rank S$ and, in particular, $\rank S\le19$;
\item\label{stage.2.lines}
$\bar\graph:=\Fn(\bS,h)\subset\graph$ (as otherwise
either $\Fano(\bar\graph)\supset S$ is not geometric or it is a proper
overlattice that can be used instead of~$S$) and $\ls|\bar\graph|\ge M$;
\item\label{stage.2.gens}
$\bS$ is rationally generated over~$S$ by a number  of
exceptional divisors~$e_i$.
\endroster
In other words, for each sufficiently
large subgraph $\bar\graph\subset\graph$ (to play
the r\^{o}le of $\Fn(\bS,h)$) we add to~$S$, one-by-one, a number of
exceptional divisors and analyze the result. We use essentially the same
algorithms (see \autoref{s.extension.algorithm} below) as for adding lines
(\cf. the proof of
\autoref{prop.Kummer.64}), with \autoref{lem.ext} replaced by more
general \autoref{lem.compatible}.
Note also that by an ``exceptional divisor'' we
merely mean an
extra vector~$e$ satisfying $e^2=-2$, $e\cdot h=0$,
and certain formal geometric conditions (\cf. Lemmas~\ref{lem.star}
and~\ref{lem.cycle} below for octics); we do \emph{not} insist that~$e$
remain
an exceptional divisor in $\Fex(\bS)$, mainly because
we do not know a usable criterion to control the behaviour of exceptional divisors.

\remark\label{rem.generating.set}
In practice, instead of $\bar\graph$ we use an even smaller subgraph
of~$\graph$ known to generate~$S$ over~$\Q$: this reduces the number of
choices for~$e$. Then, for each ``good'' vector~$e$ found,
we analyze all geometric finite index
extensions of the lattice $S+\Z e$ and check conditions~\iref{stage.2.rank},
\iref{stage.2.lines} above.
\endremark

\subsection{Partitioned graphs}\label{s.partitioned}
Let $\graph=C_0\cup C_1$, $C_0\cap C_1=\varnothing$, be a
partitioned graph,
\cf.~\eqref{eq.partition} and \autoref{rem.partitioned}.
For such a graph,
we redefine the lattice $\Fano_{2d}(\graph)$ in order to
\emph{impose} the existence of an $m$-isotropic vector~$\iso$:
\[ \label{eq-part-graph-lattice}
\Fano_{2d}^m(\graph):=(\Z\graph+\Z h+\Z\iso)/\ker,\qquad
h^2=2d>0,\quad\iso^2=0,\quad\iso\cdot h=m,
\]
where
the other intersections are
\[*
\text{$h\cdot v=1$ for $v\in\graph$},\quad
\text{$\iso\cdot v=i$ for $v\in C_i$, $i=0,1$}.
\]
As above, the degree~$2d$ is usually fixed and omitted from the notation.

Since we want to use partitioned graphs to study configurations
of rational curves in hyperelliptic
models/special octics, we always assume $\Fano^m(\graph)$ hyperbolic.

Obviously, all notions
introduced in  \autoref{s.aggraphs}
(\eg, extension $\Fano^m(\graph,\CK)$,
compatible Weyl chamber with $\iso\in\barfp$,
$\iso$-extensibility,
to name a few)
can be generalized
to partitioned graphs and the lattices that they define;
details are left to the reader.
A
root $r\in\root_0\Fano^m(\graph,\CK)$ is called \emph{$\iso$-separating}
if $r\cdot\iso>0$ and $r\cdot v<0$ for some vertex $v\in\graph$.
Then, \autoref{lem.ext} (and its proof) extends almost literally.

\lemma\label{lem.ext.partitioned}
A partitioned pair $(\graph,\CK)$ is $\iso$-extensible if and only if
the lattice
$\Fano^m(\graph,\CK)$ has neither separating nor $\iso$-separating roots.
If this is the case, there is a \emph{unique}
Weyl chamber for
$\rt(\Fano^m(\graph,\CK),h)$ compatible with~$\graph$.
\done
\endlemma

Crucial are the concepts of an $m'$-admissible and $m'$-(sub-)geometric
partitioned graph: in \autoref{def.madmissible}, we {\em always refer to
the lattice $\Fano^m$, $m=m'+1$.} Observe the relation between the two parameters:
when \emph{imposing} the presence of an $m$-isotropic vector, we always assume that
\emph{this $m$ is minimal possible}!

For partitioned graphs, \autoref{th.K3} is recast as follows
(where, by definition, ``\qgen by lines'' includes lines, $\hp$, \emph{and the
distinguished isotropic class}).

\theorem\label{th.hyperelliptic}
A partitioned graph $\graph$ is $1$-geometric in
degree~$2d\ge4$ if and only if
$\graph\cong\graphl(X,\hp)$ for a hyperelliptic
degree~$2d$ quasi-polarized \KK-surface $(X,\hp)$
such that $\NS(X)$ is \qgen by lines.
\done
\endtheorem

\theorem\label{th.octics}
A partitioned graph $\graph$ is $2$-geometric in
degree~$8$ if and only if
$\graph\cong\graphl(X,\hp)$ for a special
octic $(X,\hp)$
such that $\NS(X)$ is \qgen by lines.
\done
\endtheorem

The more precise
and general counterparts of Theorems~\ref{th.hyperelliptic}
and~\ref{th.octics}
in terms of bi-colored partitioned graphs
(\cf. \autoref{th.K3.ext} \vs. \ref{th.K3}) are left to the reader.

\section{The taxonomy of hyperbolic  graphs}\label{S.taxonomy}

In
this section we introduce necessary tools to study
the configurations of lines on certain \KK-surfaces.
Our approach is a generalization of  \cite[\S\,4]{degt:lines}.

Given a graph $\graph$, we can consider the lattice $\ZZ \graph$ and distinguish  three cases:
\begin{enumerate}
	\item
	elliptic --- if both $\sigma_+(\ZZ \graph)$ and  $\sigma_0(\ZZ \graph)$ vanish, 
	\item
	parabolic ---  if $\sigma_+(\ZZ \graph) = 0$ and  $\sigma_0(\ZZ \graph) >0$,
	\item
	hyperbolic --- if  $\sigma_+(\ZZ \graph) = 1$.
\end{enumerate}
In view of the Hodge Index Theorem, we never consider graphs with
$\sigma_+(\ZZ \graph) > 1$.

\subsection{Parabolic graphs}\label{s.parabolic}
A connected elliptic (\resp. parabolic) graph is called a
(simply laced) \emph{Dynkin diagram}
(\resp. \emph{affine Dynkin diagram}).
A detailed account of the properties of such diagrams (in particular, their classification)
can be found in  \cite{Bourbaki:Lie}.
When $\graph$ is parabolic or elliptic, we use $\mu(\graph) := \rank(\Z\graph /\ker)$ to denote its
\emph{Milnor number}.
Given a connected parabolic graph $\fiber$,
one can use the primitive positive annihilator of the
lattice $\Z\fiber$
to define \emph{the fundamental cycle}
\[
\kappal_\fiber := \sum_{v \in \fiber} n_v v,
\label{eq.kappa}
\]
   where all $n_v$ are positive
 (see \cite[\lemcit~I.2.12]{Barth.Peters.VanDeVen},
 \cite[\S\,3.1]{degt:lines}).
Since each parabolic/\penalty0 elliptic
graph $\graph$ is a disjoint union of simply laced Dynkin diagrams
and affine Dynkin diagrams, with at least  one affine component in the parabolic case, such a graph can be described by
a formal sum of the  $\bA$-$\bD$-$\bE$ types of its components.

Following \cite[\S\,4]{degt:lines}, we introduce an order on  the set of
isomorphism classes of
connected parabolic graphs (\ie, affine Dynkin diagrams):
for graphs
$\Sigma'$, $\Sigma''$ with $\mu(\Sigma')=\mu(\Sigma'')$ we let
$$\tA_{\mu} < \tD_{\mu} < \tE_{\mu};$$
otherwise,
\[*
\Sigma' <  \Sigma''  \quad \mbox{if and only if}  \quad  \mu(\Sigma') <  \mu(\Sigma'').
\]
Finally,
since $\Gs_0(\fiber)=1$ for each affine Dynkin diagram~$\fiber$, we have
the equality
\[
\ls|\pencil|=\mu(\pencil)+\#\{\text{parabolic components of $\pencil$}\}.
\label{eq.Pi.bound}
\]
that holds for
any parabolic graph~$\pencil$
(\cf. \cite[p.~614]{degt:lines}).


 \subsection{The type of a hyperbolic graph} \label{subs.type}
Now, we can discuss the hyperbolic case, which is of our primary
interest.
 It is well-known that any hyperbolic graph $\graph$ contains
 a connected parabolic subgraph, so we can make the following definition.
(It is difficult to find a reference for this classical assertion but, essentially, it
is part of the standard classification of Dynkin diagrams: it is
straightforward that simply laced Dynkin diagrams are precisely the connected
graphs \emph{not} containing one of the affine Dynkin diagrams.)

\definition \label{def-type-graph}
Let $\graph$ be a hyperbolic or parabolic graph.
 A \emph{minimal fiber} $\fiber \subset \graph$ is a  connected parabolic
 subgraph of $\graph$
 that is minimal with respect to the
order ``$<$'' introduced in
 \autoref{s.parabolic}. All minimal fibers are of the same type, and this
 common type is called the \emph{type} of~$\graph$. Alternatively, $\graph$
 is called a \emph{$\fiber$-graph}.

This taxonomy is applied to plain graphs only (mostly, plain Fano
graphs).
\enddefinition

Recall
that we are only interested in graphs~$\graph$ such that
$\Gs_+\Fano_{2d}(\graph)=1$ for some fixed $d>0$. Under this
assumption, if a graph has a hyperbolic component, all its other components
are elliptic. Hence,
given a hyperbolic graph  $\graph$ with a connected parabolic
 subgraph $\fiber$, we can consider the maximal parabolic subgraph
 of  $\graph$ that contains $\fiber$:
 \[*
 \pencil:=\fiber\cup\bigl\{v\in\graph\bigm|
 \text{$v\cdot l=0$ for all $l\in\fiber$}\bigr\}.
 \]
This subgraph is called the \emph{pencil} containing $\fiber$
(\cf. \autoref{s.pencilK3} below). We define the sets
\[*
\sec l:=\bigl\{v\in\graph\sminus\fiber\bigm|v\cdot l=1\bigr\},\ l\in\fiber,
\quad\text{and}\quad\sec\fiber:=\bigcup_{l\in\fiber}\sec l
\]
of \emph{sections} of~$\pencil$. All parabolic
components of $\pencil$ are called its \emph{fibers}.
We define the \emph{multiplicity} of a section $w\in\sec\fiber$
as
$\sum n_v$, see~\eqref{eq.kappa}, where the summation runs over all vertices
$v \in \fiber$
that intersect $w$.
A section of multiplicity~$1$ is called \emph{simple}.
It is obvious geometrically and easily follows from the assumption
$\Gs_+\Fano_{2d}(\graph)=1$ that the multiplicity of a section
is independent of
the choice of a parabolic component $\fiber\subset\pencil$;
thus, a section is adjacent to at least one vertex in each such
component.

It may appear that
the assumption that a graph is of a given type is relatively
weak. Below we discuss a simple example to show its importance/usefulness.

\example \label{exampleD4}
Let $\graph$ be a  \smash{$\tD_4$}-graph and
let $a_j$, $j=1, \ldots,5$, be the components of the \smash{$\tD_4$}-fiber~$\fiber$.
Since $\graph$ contains no triangles (\ie, subgraphs \smash{$\tA_2$}),
distinct sections that meet the same component  $a_j$ must be disjoint.
Similarly, from the absence of
quadrangles (\smash{$\tA_3$}) or
pentagons (by definition, \smash{$\tA_4  < \tD_4$}, see \autoref{s.parabolic}),
we conclude that all  sections in $\sec \fiber$ are pairwise disjoint.
Thus, the incidence matrix of
the union $\fiber\cup\sec\fiber$
is determined  by the five integers
$\ls|\sec a_j|$, $j=1,\ldots,5$.
\endexample


\subsection{Elliptic pencils on \KK-surfaces}\label{s.pencilK3}
Given a
quasi-polarized
\KK-surface $(X, \hp)$ such that the Fano graph
$\graphl(X,\hp)$ is either elliptic or
parabolic,
it is easy to see
(see, \eg, \cite{degt:lines}) geometrically that
$\ls|\graphl(X,\hp)|\le24$.
Thus, we have
\begin{equation} \label{eq-fano-must-be-hyperbolic}
\mbox{the graph $\graphl(X,\hp)$ of a \KK-surface with at least 25 lines is hyperbolic.}
\end{equation}
In particular,
the existence of many lines on a \KK-surface
implies the existence of
a genus-one fibration
(given by the linear system $\ls|\kappal_\fiber|$, see~\eqref{eq.kappa},
where
$\fiber\subset\graphl(X, \hp)$ is a connected parabolic subgraph)
with at least one reducible fiber consisting entirely of lines.
On the other hand,
if $X\to\Cp1$ is a genus-one fibration
on a \KK-surface~$X$, the dual adjacency
graph $\tpencil$ of the components of its reducible fibers is a union of
Dynkin diagrams $\fiber_s$ (elliptic or affine),
$\tpencil=\fiber_0\cup\ldots\cup\fiber_k$, and
we have
\[
\mu(\tpencil)\le18,\qquad
\ls|\tpencil|\le24,\qquad
\nodal(\tpencil)\le16,
\label{eq.pencil}
\]
where the first  (\resp. second, \resp. third) inequality follows from
$\sigma_-\operatorname{NS}(X)\le19$
(\resp. Euler
characteristic count, \resp. \cite{Nikulin:Kummer};
recall that $\Ga$ stands for the independence number).
Obviously, the projective degree of a curve is a function
$\deg\:\tpencil\to\NN$, $v\mapsto v\cdot\hp$, and the
\emph{linear} fiber components constitute the subgraph
\[*
\pencil:=\bigl\{v\in\tpencil\bigm|\deg v=1\bigr\},
\]
which we call a \emph{combinatorial pencil}.

Assuming that $\pencil\supset\fiber_0=:\fiber$ is a $\fiber$-graph, the
degree function has the following properties:
\[
\deg(\fiber_s)=\deg(\fiber),\quad
\deg\big|_\fiber\equiv1,\quad
\deg\big|_{\fiber_s}\not\equiv1\text{ unless }\fiber_s\ge\fiber,
\label{eq.degree}
\]
where $0\le s\le k$ and $\deg(\fiber_s):=\sum_{a\in\fiber_s}n_a\deg a$ is the
total degree of a fiber (see \eqref{eq.kappa} for the definition of
the multiplicities~$n_a$).

All configurations of singular fibers of elliptic fibrations on
complex \KK-surfaces were enumerated in \cite{Shimada:ellipticK3}.
Given a genus-one fibration
on an algebraic  \KK-surface,
its Jacobian fibration
has singular fibers of the same type and it
also is a \KK-surface
(see \cite[\S\,11.4]{Huybrechts}).
(In fact,
a combinatorial pencil in the graph of lines of a quasi-polarized
\KK-surface with sufficiently many lines does have a simple section;
the particular case of octics is
discussed in \autoref{lem.Jacobian-moved} below.)
Thus, all combinatorial pencils that
we are to consider
can be derived from the list in \cite{Shimada:ellipticK3},
with the extra constraints given by~\eqref{eq.degree} taken into
account.

\subsection{The approach} \label{s.approach}
For the reader's convenience, before passing to the technical details,
we sketch
the reasoning  that
will lead us
to the proof of \autoref{thm-main}.
As
explained in \autoref{S.intro},
similar
strategy should yield sharp
upper bounds on the number of lines on
quasi-polarized degree-$2d$
\KK-surfaces for $d>4$.

In view of \eqref{eq-fano-must-be-hyperbolic}, we  focus our attention
on hyperbolic graphs.
For Stage~$1$ (see \autoref{ss.stage.1}), we
fix a \emph{type}~$\fiber$ and consider $\fiber$-graphs in the form
\[*
\graph\supset\pencil\supset\fiber,
\]
where $\fiber\subset\graph$ is a distinguished \emph{fiber} (any
fixed
subgraph of~$\graph$ within the chosen isomorphism class) and
$\pencil$ is the pencil containing~$\fiber$.
Formally, both $\pencil$ and $\sec\fiber$ depend on~$\graph$, but we omit
this fact in the notation since, in fact, we usually construct $\graph$
\emph{starting} from a given pencil~$\pencil$ or collection of sections
$\sec\fiber$.

We fix a threshold~$M$ and
prove that, with a few exceptions
that are
listed
explicitly, the
inequality
\[
\ls|\graph|<\M
\label{eq.goal}
\]
holds for any geometric $\fiber$-graph~$\graph$. To justify this claim,
we
choose another
pair of thresholds $\M_\fiber$, $\M_\pencil$ such that
$\M_\fiber+\M_\pencil\le\M+1$
and show that, with
the same exceptions, the following inequalities hold
\[
\alignedat2
&\text{$\ls|\graph|<\M$ whenever $\ls|\sec\fiber|\ge\M_\fiber$},&\quad&
\text{see \autoref{S.Msigma}},
\ \text{and}\\
&\text{$\ls|\graph|<\M$ whenever $\ls|\pencil|\ge\M_\pencil$},&&
 \text{see \autoref{S:largepencils}}.
\endalignedat
\label{eq.bounds}
\]
Since, obviously,
\[*
\graph=\pencil\cup\sec\fiber\quad\text{(disjoint union)},
\]
the two inequalities imply~\eqref{eq.goal}.
Roughly speaking, the choice of the type $\fiber$ imposes certain extra constraints
on the Gram matrix of the sections (\cf. \autoref{exampleD4}). Thus, adding  many sections to $\fiber$ rules out (almost) all configurations of singular fibers
allowed by \cite{Shimada:ellipticK3}. On the other hand, once we choose
many fiber components, the resulting lattice has high rank and cannot
accommodate many sections.

Thus, typically, we start with a \emph{sufficiently large}
reasonably ``standard'' graph~$\graph_0$ and extend it by adding a number of extra
vertices. Each extension $\graph\supset\graph_0$ that is not
subgeometric is disregarded immediately, \cf. \autoref{rem.heredity}.
Crucial is the fact that the rank
of a (sub-)geometric graph does not exceed~$20$ --- the maximal Picard rank of
a \KK-surface.
Hence,
once an extension $\graph\supset\graph_0$ of this maximal
rank has been obtained (usually, after adding but a few extra vertices),
there remains to apply Nikulin's theory and study the
finite index extensions of the lattice $\Fano(\graph)$: we
select those that are geometric and compute
the respective saturations of~$\graph$.
It is this fact that makes our algorithms converge reasonably fast.

The last stage of the proof (see \autoref{ss.stage.2}) boils down to the
study of but a few extremal configurations. The general approach is outlined
in \autoref{ss.stage.2}, and we give the necessary explanations on a
case-by-case basis.

The general computer aided algorithms are
described
in the appendices (the reason for this decision is explained in \autoref{remark-no-code}), and a
more detailed exposition
and numerous intermediate statements for
the particular case of
octic
\KK-surfaces
are
the  subject of the rest of the paper.
A human readable example of this reasoning ---the case of Kummer octics---
is found in \autoref{S.Kummer}.


\section{A few extreme cases}\label{S.extreme}

Starting from this section, we confine ourselves to
the Fano graphs of octic
\KK-surfaces $(X,\hp)$; thus,
we fix the degree $2d=h^2=8$ and consider
 $m$-admissible $8$-polarized
lattices/graphs, where $m \geq 2$.

\subsection{The star of a pseudo-vertex}\label{s.star}

We start with discussing the local structure,
ruling out a few very simple subgraphs of the Fano graphs of
octics.
We refer to \autoref{def.pseudo} for the definition of a pseudo-vertex
and its \sstar.

\lemma\label{lem.star}
The \sstar~of a \rom(pseudo-\rom)vertex~$e$ of a
$2$-admissible
graph~$\graph$ can be as follows\rom:
\roster
\item\label{star.line}
$\Star(e)\cong\bA_2\oplus a\bA_1$, $a\le3$, or $a\bA_1$, $a\le8$,
if $e\in\graph$ is a line\rom;
\item\label{star.line.admissible}
$\Star(e)\cong a\bA_1$, $a\le7$,
if $e\in\graph$ is a line and $\graph$ is $3$-admissible\rom;
\item\label{star.ex}
$\Star(e)\cong a\bA_1$, $a\le5$,
if $e\in\basis\Weyl$ is an exceptional divisor\rom;
\item\label{star.isotropic}
$\Star(e)\cong\bA_2$ or $a\bA_1$, $a\le9$,
if $e\in\fp$ is a $3$-isotropic vector.
\endroster
\endlemma

\proof
We change the notation and confine ourselves to the
subgraph $\graph:=\Star(e)$ and sublattice $S:=\Z\graph+\Z h+\Z e$,
with $v\cdot e=1$ for each $v\in\graph$. Assuming that $\graph$ contains a
subgraph $\bA_2$, generated by a pair $u_1,u_2$, and/or $a\bA_1$, generated by
$v_1,v_2,\ldots$, we make the following observations:
\roster*
\item
if $e\in\graph$ is a line, then $u_1+u_2+e$ is a $3$-isotropic vector;
\item
if $e\in\graph$ and $a\ge8$, then $3e-h+v_1+\ldots+v_8$ is a $3$-isotropic vector;
\item
if $e\in\graph$ and $a\ge4$, then
$2e-h+u_1+u_2+\ldots+v_4$ is a separating root;
\item
if $e\in\basis\Weyl$ is an exceptional divisor, then $u_1+u_2-e$ is
$2$-isotropic;
\item
if $e\in\fp$ is $3$-isotropic, then
$2e-h+u_1+u_2+v_1$ is a $e$-separating root.
\endroster
All other border cases, \ie, $\graph\supset\tA_2$, or $\graph\supset\bA_3$, or
$\graph\supset b\bA_2\oplus a\bA_1$ with the pair $(a,b)$ exceeding that
announced in the statement, are ruled out by $\Gs_+S\ge2$.
\endproof

\lemma\label{lem.2stars}
Let $v_1,v_2\in\graph$ be two vertices in a
 $2$-admissible
graph~$\graph$. Then the intersection
$\Star(v_1)\cap\Star(v_2)$ is discrete and\rom:
\roster
\item\label{star.u.v=1}
if $v_1\cdot v_2=1$, then $\ls|\Star(v_1)\cap\Star(v_2)|\le1$\rom;
\item\label{star.u.v=0}
if $v_1\cdot v_2=0$, then $\ls|\Star(v_1)\cap\Star(v_2)|\le3$.
\endroster
If, in addition, $\graph$ is $3$-admissible, then
both inequalities above are strict, \ie,
$\graph$ is both \emph{triangle}- and \emph{biquadrangle free}.
\endlemma

\proof
The first inequality is a restatement of
\autoref{lem.star}\iref{star.line}
and~\iref{star.line.admissible}. For the second one, let
$u_1,\ldots,u_k\in\Star(v_1)\cap\Star(v_2)$. If $k\ge4$ or $k\ge2$ and
$u_1\cdot u_2=1$, one has $\Gs_+\graph>1$. Thus, all vertices are
pairwise disjoint and $k\le3$. In the latter case $k=3$,
the vector $h-v_1-v_2-u_1-u_2-u_3$ is $3$-isotropic.
\endproof

As an immediate consequence, we show that, in a large $2$-geometric graph, any pencil has a simple
section.

\lemma\label{lem.Jacobian-moved}
Let $\graph\supset\pencil\supset\fiber$ be a $2$-geometric $\fiber$-graph,
where $\pencil$ is a pencil,
and assume that
$\ls|\graph|>24$.
Then there is at least one simple section of~$\pencil$ in~$\graph$.
\endlemma

\proof
Clearly, $\graph\sminus\pencil\ne\varnothing$, and we only need to show that
\emph{multiple} sections would not suffice to achieve the goal
$\ls|\graph|>24$.

If $\fiber=\tA_4$ or
$\fiber\ge\tA_5$, there are no multiple sections by \autoref{def-type-graph}.

If $\fiber=\tA_2$, multiple sections are ruled out by
\autoref{lem.star}\iref{star.line}.

If $\fiber=\tA_3$,
then, by \autoref{lem.2stars}\iref{star.u.v=0}, there can be at most two (in
fact, at most one) double sections.
However, in
the presence of a double section,
the pencil~$\pencil$ may have at most~$4$
parabolic components, see
\autoref{lem.star}\iref{star.line}; therefore, by~\eqref{eq.Pi.bound}
and~\eqref{eq.pencil}, we have
$\ls|\graph|=\ls|\pencil|+2\le24$.

Finally, if $\fiber=\tD_4$, then, by \autoref{lem.star}\iref{star.line},
there may be $s\le4$ double sections, intersecting the central fiber
$a_1\in\fiber$.
Let $p\ge1$ and $0\le q\le2$ be the numbers of, respectively,
type~$\tD_4$ and type~$\tA_5$ fibers in~$\pencil$. We have $p+q\le4$ and,
hence, $\ls|\pencil|\le22$, see~\eqref{eq.Pi.bound} and~\eqref{eq.pencil}.
Consider the following possibilities:
\roster*
\item
if $s=4$, then $p+q=1$ and $\ls|\graph|=\ls|\pencil|+4\le23$,
see~\eqref{eq.Pi.bound};
\item
if $2\le s\le3$, then $p=1$; hence, $p+q\le3$ and
$\ls|\graph|\le\ls|\pencil|+3\le24$;
\item
if $s\le1$, then $\ls|\graph|\le\ls|\pencil|+1\le23$.
\qedhere
\endroster
\endproof

The following lemma holds for any
polarization $h^2\ge4$. It is obvious for graphs of lines
on  \KK-surfaces --- we
state it merely for future reference.
We formulate it for bi-colored graphs $\tilde{\graph}$, but it
can be applied to
plain graphs as well.
In this case,
we endow  all vertices of a plain $\graph$ with color $1$.

\lemma\label{lem.cycle}
Let $\tilde\graph$ be a bi-colored graph such that the lattice
 $\Fano_{2d}(\tilde{\graph})$ is $2$-admissible and hyperbolic.
Then
\roster*
\item
each triangle in~$\tilde\graph$ consists of lines only, and
\item
a quadrangle in $\tilde\graph$ contains at most one exceptional divisor.
\endroster
More generally, each cycle in~$\tilde\graph$ contains at least three lines
\rom(at least four, if
$\Fano_{2d}(\tilde{\graph})$
is also required to be $3$-admissible\rom).
\endlemma

\proof
Assume the contrary and consider $e:=\sum v$, the summation running over the
cycle in question. We have $e\cdot h=\#(\text{lines in the cycle})$ and
$e^2\ge0$. Hence, either $e$ is $1$- or $2$-isotropic or $\Gs_+(\Z h+\Z e)=2$.
\endproof

According to \autoref{lem.star}\iref{star.line.admissible},
the presence
of a triangle (\ie, an $\tA_2$ subgraph)
in a $2$-geometric graph $\graph$ automatically implies that
$\Fano(\graph)$ has a $3$-isotropic vector.
For such graphs, we have the following lemma (recall \autoref{eq-graphex}).

\lemma[see \autoref{proofs.pencil}]\label{lem.a2p}
Let
$\pencil$ be an $\tA_2$-pencil and $\ls|\pencil|\ge\M_\pencil:=21$. Then,
for any $2$-geometric extension $\graph\supset\pencil$, one has either
\roster*
\item
$\graphex\cong\config{TC.33}$ \rom(see \autoref{tab.main}\rom) and
\smash{$\pencil\cong8 \tA_2$}, or
\item
$\graphex\cong\config{TC.30.3}$ \rom(see \autoref{tab.main}\rom) and
\smash{$\pencil\cong6 \tA_2 \oplus 3\bA_1$},
\endroster
or $\ls|\graph|<30$.
In the case $\graphex\cong\config{TC.30.3}$,
the three $\bA_1$-components and three exceptional
divisors constitute a single \smash{$\tA_5$} type fiber of the elliptic
pencil.
\endlemma


\subsection{Special octics}\label{s.special}
In this section we study line configurations on special octics
(\cf. \autoref{def.triquadric}).
We fix a degree $8$ quasi-polarized \KK-surface $(X, \hp)$ and a 3-isotropic
vector $\pdp \in \Nef(X)$.
Moreover, we assume that
$\NS(X) \ni \hp$ is $2$-admissible.

The linear system $|\pdp|$ endows $X$ with a genus-one fibration
$f_\pdp\:  X \rightarrow \Cp1$.
As in \autoref{s.partitioned}, we obtain a partition of the Fano graph $$\graphl(X,\hp) = C_0  \cup C_1,$$ where $C_1$ (\resp. $C_0$)
are sections (\resp. fiber components) of the fibration $f_\pdp$. 


\lemma \label{lemm.pdp-notriangles}
If $C_0$ contains no triangle, then $\ls|C_0| \le 18$.
\endlemma

\proof
Assume
that $\ls|C_0| > 18$.
Since each singular fiber of $f_\pdp$ contains at most three lines,
at least seven singular fibers contain lines. By the assumption, each
such fiber contains an extra component,
\viz. the one of degree
other than~$1$.
(Recall that,
by the $2$-admissibility, two lines cannot form an $\tA_1$-fiber.) Thus,
the
singular fibers of  $f_\pdp$  contain at least $26$ rational curves.
Contradiction.
\endproof

\lemma\label{lem.6A2}
If \smash{$C_0\supset6\tA_2$}, then any exceptional divisor on~$X$ is
orthogonal to~$\pdp$.
\endlemma

\proof
If $e$ is an exceptional divisor and $e\cdot\pdp>0$, then $e$ intersects each
fiber of~$f_\pdp$.
Hence, $e$
intersects at least six lines, contradicting to
\autoref{lem.star}\iref{star.ex}.
\endproof

\lemma\label{lem.7A2}
If \smash{$C_0\supset7\tA_2$}, then $X$ is smooth.
\endlemma

\proof
Let $e$ be an exceptional divisor in~$X$. By \autoref{lem.6A2}, $e$ is
a fiber component of~$f_\pdp$ and, hence, $e$ is orthogonal to the sublattice
\smash{$\Fano(7\tA_2)$}, \ie,
$e\in T:=\smash{\Fano(7\tA_2)^\perp}\subset H_2(X)$.
On the other hand, it is easily seen that, up to
automorphism, \smash{$7\tA_2$} has a unique geometric kernel and the
lattice $T\cong\bU(3)^2\oplus\bA_2(3)$ is root
free.
\endproof

\subsection{Proof of \autoref{th.special}}\label{proof.special}
We keep the notation of \autoref{s.special} and follow
	the strategy
outlined in \autoref{s.extendedversusplain}.  By
\autoref{lem.star}\iref{star.isotropic},
\begin{equation} \label{eq-c1nine}
\ls|C_1| \le 9.
\end{equation}
If $C_0$ contains no triangles, then \autoref{lemm.pdp-notriangles} yields $\ls|\graphl(X,\hp)|\leq27$.
Thus, we can assume that $\ls|C_0| \ge 21$ and $C_0$ contains a triangle,
in which case,
by \autoref{lem.a2p}, the statement follows from \autoref{th.octics}
\emph{ provided that
$\NS(X)$ is \qgen
by lines}. This completes Stage 1 of the proof (\cf. \autoref{ss.stage.2}).

For Stage~$2$, we need to extend
$\graph:=\spp_1\graphex$, where
$\graphex\cong\config{TC.30.3}$ or
$\config{TC.33}$, by an extra exceptional divisor~$e$. In the notation therein,
for $\graphex\cong\config{TC.30.3}$ we have to take $\bar\graph=\graph$.
By \autoref{lem.6A2}, the
new divisor~$e$ would be a fiber component of~$f_\pdp$;
hence,
by
\autoref{lem.a2p}, 
 it would coincide with one of the three exceptional divisors given by the
 three $0$-vertices of  $\graphex\cong\config{TC.30.3}$.
(Note
that we do \emph{not} assert that the existing exceptional divisors
remain such. However, the \emph{curves} do not disappear: they may
only become
reducible, adding even more fiber components to the elliptic pencil
which already has $24$ components.)

Thus, assume that $\bar\graph\subset\config{TC.33}$ and
$\ls|\bar\graph|\ge30$.
Obviously, the new octic $X$ contains at least five
\smash{$\tA_2$}
connected components of $\pencil$;
hence,
 the new part~$C_0$ must be one of
 \smash{$8\tA_2$}, \smash{$7\tA_2\oplus\bA_1$}, or \smash{$7\tA_2$}
 (see \autoref{lem.a2p} and observe that,
because of the degree, a \emph{single} line cannot be removed from an
\smash{$\tA_2$} component.)
 In
 each case, by \autoref{lem.7A2}, the octic
would have to
remain smooth,
contradicting the fact that an exceptional divisor has been added.
 \qed


\section{Kummer octics}\label{S.Kummer}

In this section,
we treat separately the Kummer and so-called \emph{almost Kummer}
(see \autoref{s.Kummer-} below) octics. We have two reasons to single out
these two classes. On the one hand, the corresponding lattices have very
large geometric kernels and our standard na\"{\i}ve algorithms fail due to
the  lack of simplification (\cf. \autoref{wrn.sat}). On the other hand,
these classes serve as an example where our approach can be explained in
detail in a human readable form.

\subsection{The Golay code\pdfstr{}{ \rm(see \cite[Chapter 11]{Conway.Sloane})}} \label{s.Golay}
Consider the extended binary Golay code~$\CC_{24}$, pick a codeword~$\fo$ of
length~$16$, and denote by~$\CC$ the code $\{o\in\CC_{24}\,|\,o\subset\fo\}$:
it consists of~$\varnothing$, $\fo$, and $30$ octads.
For a subset $s\subset\fo$, we let $[s]:=\sum_{l\in s}l\in\Z\fo$.

We distinguish certain subsets of~$\fo$. Namely, define
\[*
\CS:=\bigl\{o\cap\fo\bigm|o\in\CC_{24}\bigr\}\sminus\CC
=\CS_4\cup\CS_6\cup\CS_8\cup\CS_{10}\cup\CS_{12},
\]
where the last splitting is according to the length of a subset $s\in\CS$.
The involution $s\mapsto\bar s:=\fo\sminus s$ sends $\CS_n$ to $\CS_{16-n}$.
Define also a $\bar\quad$-invariant equivalence relation
\[*
\text{$r\sim s$ if and only if $r\vartriangle s\in\CC$},
\]
where $\vartriangle$ is the symmetric difference. Then,
\roster*
\item
$\CS_4\cup\CS_8\cup\CS_{12}$ splits into $35$ equivalence classes of
size~$4+24+4$;
\item
$\CS_6\cup\CS_{10}$ splits into $28$ equivalence classes of size~$16+16$.
\endroster
The equivalence class of a set $\o\in\CS$ is denoted by $\eclass\o$,
and we let $\eclass\o_n:=\eclass\o\cap\CS_n$.
For a fixed element $\o\in\CS$ and residue $m\in\Z/2$, define
\[*
\CS_n(\o,m):=\bigl\{s\in\CS_n\bigm|\ls|s\cap\o|=m\bmod2\bigr\};
\]
this set depends on the class $\eclass\o$ only.

As an alternative description, for a subset $s\subset\fo$ we have
\[
\text{$\ls|s\cap o|=0\bmod2$ for all $o\in\CC$ iff $s\in\CS\cup\CC$}.
\label{eq.even.sets}
\]

The group $\Aut\CC$ of automorphisms of~$\CC$
is the stabilizer of~$\fo$ in $M_{24}$: it has order $322560$,
preserves $\bar\ $ and $\sim$, and acts transitively
on each~$\CS_n$.

According to Nikulin~\cite{Nikulin:Kummer}, sixteen is the maximal number of
pairwise disjoint smooth rational curves in a \KK-surface,
not necessarily polarized. Identifying the $16$ rational curves with
the elements of~$\fo$, up to isomorphism, the only
finite index extension of~$\Z\fo$ admitting a primitive embedding
to~$\bL$ has kernel
\[
\CK_\fo:=\bigl\{\tfrac12[o]\bmod\Z\fo\bigm|o\in\CC\bigr\};
\label{eq.K.o}
\]
indeed, this is the only subspace $V\subset\discr\Z\fo\cong\F_2^{16}$ of
dimension
$\dim V\ge5$ and minimal Hamming distance $\ge8$.
This statement is easily taken down to $13$ lines.
For each $n=1,2,3$, there is a unique, up to $\Aut\CC$, subset
$\star\subset\fo$, $\ls|\star|=n$; fixing such a subset and identifying
the $(16-n)$ lines with the elements of $\fo^\star:=\fo\sminus{\star}$,
we conclude that the only geometric finite index extension of $\Z\fo^*$ has
kernel
\[
\CK_\fo^\star:=\bigl\{\tfrac12[o]\bmod\Z\fo^\star\bigm|o\in\CC^\star\bigr\},\qquad
\CC^\star:=\bigl\{o\in\CC\bigm|o\subset\fo^\star\bigr\}.
\label{eq.K.star}
\]
Note, though, that if $\ls|\star|>1$, the automorphism group $\Aut\CC^\star$
is larger than the mere restriction of the stabilizer of~$\fo^\star$ in $M_{24}$.

\subsection{Kummer octics}\label{s.Kummer}
In this paper, by
a \emph{Kummer octic} we mean
a degree-$8$ birationally quasi-polarized \KK-surface
$(X,\hp)$ 
with a distinguished
collection of $16$
lines pairwise disjoint in~$X$.
We identify the $16$
lines with
the elements of~$\fo$, regarding the latter as a graph without edges.

There is extensive literature on projective models of Kummer surfaces, but we were
unable to find the following statement, \viz. the existence of
\emph{exactly} two families of Kummer octics as defined above.

\theorem\label{th.Kummer}
There are two disjoint deformation families $\Kummer_d$, $d=64$ or~$256$,
of Kummer octics\rom; they are distinguished by the determinant
\[*
d:=\ls|\det\NS(X)|,\quad \text{where $(X,\hp)\in\Kummer_d$ is a generic member}.
\]
\rom(Alternatively, for \emph{any} $(X,\hp)\in\Kummer_d$, the parameter $d$
is recovered as the determinant of
$\NS(X)\cap(\Q\fo+\Q \hp)$.\rom)
Both families have dimension~$3$ and consist of triquadrics.
\endtheorem

\proof
The lattice $\Fano(\fo,\CK_\fo)$, see~\eqref{eq.K.o},
is not geometric. However,
starting from~$\CK_\fo$, it is easy to list all geometric kernels
$\CK\supset\CK_\fo$ of~$\fo$. There are two:
\begin{alignat}3
&\CK_{64},&\quad&\text{generated over $\CK_\fo$ by
$\textstyle\frac14\hp+\frac18[\fo]+\frac12[o]$, $o\in\CS_6$, or}
\label{eq.Kummer.64}\\
&\CK_{256},&&\text{generated over $\CK_\fo$ by
$\textstyle\frac12\hp+\frac14[\fo]+\frac12[o]$, $o\in\CS_8$}
\label{eq.Kummer.256}
\end{alignat}
(see \autoref{s.Golay} for the notation),
which give rise to the two families in the statement.

To show that all Kummer octics are triquadrics, we try to add to~$\fo$ a
$3$-isotropic vector~$\iso$.
\autoref{lem.star}\iref{star.isotropic} asserts that
$\bigl|\|\iso\|_\fo\bigr|\le9$; however,
in the presence of all $16$ lines,
only $\bigl|\|\iso\|_\fo\bigr|\in\{0,1,2\}$ passes the Sylvester
test (see \autoref{lem.Sylvester}).
By~\eqref{eq.even.sets}, none of these options passes
the kernel test (see \autoref{lem.kernel}).
\endproof

\remark \label{rem.Kummer.ref}
A construction of the family
$\CK_{64}$ and a description of its relation to the
Jacobians of genus-$2$ curves goes back to Klein (see
\cite[$\S$~10.3.3]{Dolgachev:book} and
references therein).
An exposition of the construction of surfaces in $\CK_{64}$
as Kummer surfaces
can be found in
\cite[$\S$.5.1]{GS:Kummer}. Imitating the arguments from
\latin{loc.\ cit.},
one can use the equality
$\ls|\det\NS(X)| = 256$ to show that,
for a generic $X \in \CK_{256}$, we have
$$
T_X \simeq \langle -16 \rangle \oplus U_2 \oplus U_2
$$
and $X$ is the Kummer surface of a $(1,4)$-polarized Abelian surface $A$.

It is well-known that the ample polarization on $A$ induces  a pseudo-ample divisor $\hat{h}$ on $X$ and the sum of the 16 pairewise disjoint rational curves $C_j$, $j=1,\ldots,16$, on $X$ is even (see \cite{Nikulin:Kummer}).
One can check that the divisor on generic $X \in  \CK_{256}$ that defines its Kummer octic model is given as
$h := (\hat{h} - \frac{1}{2} \sum_{j=1}^{16} C_j)$
(\cf. \cite[$\S$.5]{GS:Kummer}).
A more detailed discussion of this construction
is beyond
the scope of this paper.
\endremark

\proposition\label{prop.Kummer.64}
If $(X,\hp)\in\Kummer_{64}$, then either
\roster*
\item
$X_8$ is smooth and $\Fn(X,\hp)\cong\config{QC.32.K}$,
see~\cite{degt:lines} and \autoref{tab.main}, or
\item
$\Sing(X_8)=8\bA_1$ and $\Fn(X,\hp)=\fo$, see \autoref{ex.16-8} below and
\autoref{fig.K16-8}.
\endroster
The \KK-octics in $\Kummer_{64}$ with $\Fn(X,\hp)=\fo$
\rom(\ie, such that $\Sing X_8 \neq \emptyset$\rom)
 constitute a family of
dimension~$2$.
\endproposition

\proof
Fix an element $\o\in\CS_6$ and consider the corresponding kernel $\CK_{64}$,
see~\eqref{eq.Kummer.64}. The stabilizer $G_{64}:=\stab\CK_{64}$ is an order
$11520$ subgroup of $\Aut\CC$.

Using Lemmas~\ref{lem.star}\iref{star.line} and~\ref{lem.kernel},
if an extra line~$v$ can be added to~$\fo$, then
\[*
\|v\|_\fo\in\CS_6(\o,0)
=\eclass\o_6\cup\{\text{another $G_{64}$-orbit}\}.
\]
The former orbit is already present in $\sat(\fo,\CK_{64})$,
resulting in \config{QC.32.K}. If $\|v\|_\fo\notin\eclass\o_6$,
the pair $(\fo\cup v,\CK_{64})$ is not extensible.

Similarly, if an exceptional divisor~$\ex$ can be added to~$\fo$, then
$\|\ex\|_\fo\in\CS_4(\o,1)$
(a single $G_{64}$-orbit).
On this set, we have a coarser equivalence relation
\[*
r\approx s\quad\text{whenever}\quad r\sim s\ \text{or}\ r\sim s\vartriangle\o,
\]
and the lattice
$\Fano(\fo,\CK_{64})+\Z\ex$
contains the whole octuple $\eclass\ex_\approx$.
On the other hand, if
$\|\ex_1\|_\fo\not\approx\|\ex_2\|_\fo$,
the bi-colored graph $\fo\cup\ex_1\cup\ex_2$ violates the Sylvester test
(\autoref{lem.Sylvester}),
no matter whether $\ex_1\cdot\ex_2=0$ or~$1$ (where the latter option
needs to be
considered only if $\|\ex_1\|_\fo\cap\|\ex_2\|_\fo=\varnothing$, see
\autoref{lem.cycle}).
\endproof

\example \label{rem.kum.sing}
The surfaces in the family $\Kummer_{64}$ (see \autoref{prop.Kummer.64})
demonstrate the phenomenon
discussed in
\autoref{wrn.sat}:
extending the geometric lattice  $\Fano_8(\fo,\CK_{64})$ by an exceptional divisor results in
\KK-surfaces
with fewer ($32\mapsto16$) lines!

Note that, \emph{accidentally}, $\fo$ still is geometric, but with a
different kernel, as a member of the other family $\Kummer_{256}$,
see \autoref{prop.Kummer.256} below.
\endexample

\example\label{ex.16-8}
The second configuration in \autoref{prop.Kummer.64} is one of the very few
that can be (and are, upon popular request) drawn, see
\autoref{fig.K16-8}, which depicts lines and exceptional
divisors as columns and rows, respectively.
\figure
\hbox to\hsize{\hss$
	\makeatletter
	\setbox0\hbox{$\m@th\joinrel\relbar\joinrel\relbar\joinrel$}%
	\def\zbox#1{\kern-\wd0\hbox to\wd0{\hss$#1$\hss}}%
	\def\1{\zbox\bullet}\def\.{}\let\\\cr
	\vcenter{\offinterlineskip\halign{&\copy0\zbox|#\cr
			\1&\1&\1&\1&\.&\.&\.&\.&\.&\.&\.&\.&\.&\.&\.&\.\\
			\.&\.&\.&\.&\1&\1&\1&\1&\.&\.&\.&\.&\.&\.&\.&\.\\
			\.&\.&\.&\.&\.&\.&\.&\.&\1&\1&\1&\1&\.&\.&\.&\.\\
			\.&\.&\.&\.&\.&\.&\.&\.&\.&\.&\.&\.&\1&\1&\1&\1\\
			\1&\.&\.&\.&\1&\.&\.&\.&\1&\.&\.&\.&\1&\.&\.&\.\\
			\.&\1&\.&\.&\.&\1&\.&\.&\.&\1&\.&\.&\.&\1&\.&\.\\
			\.&\.&\1&\.&\.&\.&\1&\.&\.&\.&\1&\.&\.&\.&\1&\.\\
			\.&\.&\.&\1&\.&\.&\.&\1&\.&\.&\.&\1&\.&\.&\.&\1\\}}
	$\hss}
\makeatother

\caption{The Kummer octic with $16$ lines and $8$ nodes
	(see \autoref{ex.16-8})}\label{fig.K16-8}
\endfigure
The eight exceptional divisors are pairwise disjoint, so that the
respective surface~$X_8$ has
eight nodes. Each line contains two nodes and each node is a point of
intersection of four lines.
The lattice $\NS(X_8)$ is an index~$32$ extension of
$\Fano(\Fex X_8)$; the material of this section (\autoref{th.Kummer} and
Propositions~\ref{prop.Kummer.64}, \ref{prop.Kummer.256}) asserts that, up to
abstract graph automorphism, $\Fex X_8$
admits a unique geometric kernel.
More explicitly, this kernel is determined by the closure (under $\vartriangle$)
of the set $\CC^*\cup o$, \cf. \eqref{eq.Kummer.64}, and the latter closure
is the collection of all \emph{even} symmetric differences of the supports of
the exceptional divisors.
\endexample

\proposition\label{prop.Kummer.256}
If $(X,\hp)\in\Kummer_{256}$, then either
\roster*
\item
$X_8$ is smooth and $\ls|\Fn(X,\hp)|=16$, $20$, $24$, or~$28$
\rom(with $\Fn(X,\hp)\cong\config{SC.28'}$\rom),
or
\item
$\Sing(X_8)=4\bA_1$ and $\ls|\Fn(X,\hp)|=16$, $24$, or~$32$
\rom(with $\Fex(X,\hp)\cong\config{QC.32.4}$\rom).
\endroster
The configurations \config{SC.28'} \rom(see~\cite{degt:lines}\rom)
and \config{QC.32.4} are described in Tables~\ref{tab.main}
and~\ref{tab.mainsmall}.
\endproposition

\proof
Fix an element $\o\in\CS_8$ and consider the corresponding kernel $\CK_{256}$,
see~\eqref{eq.Kummer.256}. The stabilizer $G_{256}:=\stab\CK_{256}$ is an order
$9216$ subgroup of $\Aut\CC$.
By Lemmas \ref{lem.star}\iref{star.line} and~\ref{lem.kernel},
if an extra line~$v$ can be added to~$\fo$, then
\[*
\|v\|_\fo\in\CS_4(\o,1)\cup\CS_6(\o,0)
\]
(two $G_{256}$-orbits). In the former case, the saturation
$\sat(\fo\cup v,\CK_{256})$ contains the whole quadruple
$\eclass{v}_4$; in the latter case,
the pair $(\fo\cup v,\CK_{256})$ is not
extensible.
Adding (in the progressive mode,
see \autoref{ss.step.r.general} below) up to three independent lines, we
obtain all but one configurations in the statement (those
\qgen by lines).

Similarly, if an exceptional divisor~$\ex$ can be added to~$(\fo,\CK_{256})$,
then
\[*
\|\ex\|_\fo\in\CS_4(\o,0)
=\eclass\o_4\cup\{\text{another $G_{256}$-orbit}\}.
\]
In the former case, the lattice
$\Fano_\ex:=\Fano(\fo,\CK_{256})+\Z\ex$
has no geometric
extensions; in the latter case, the quadruple $\eclass\ex_4$
gives rise to four nodes in
$\Fano_\ex$.
An
attempt to add an extra exceptional divisor to any of the
configurations obtained above
either fails or produces a configuration that is already on
the list.
\endproof

\subsection{Almost Kummer octics}\label{s.Kummer-}
We define an \emph{almost Kummer} octic $(X,\hp)$
(\resp . $X_8\subset\Cp5$)
 as a degree-$8$ birationally quasi-polarized
\KK-surface (\resp. the image $X_8:=f_{\hp}(X))$
containing at least $15$
pairwise disjoint lines .

\proposition\label{prop.Kummer-}
There is a unique deformation family $\Kummer^\star$ of almost Kummer
octics\rom; they are all triquadrics. For any $(X,\hp)\in\Kummer^\star$, one has
$\ls|\Fn(X,\hp)|\le29$ unless $(X,\hp)$ is Kummer \rom(see \autoref{s.Kummer}\rom) or
$\Fn(X,\hp)\cong\config{QC.33}$, see \cite{degt:lines} and \autoref{tab.main}.
\endproposition

\proof
We fix a one-element set ${\star}\in\fo$ and start with the pair
$(\fo^\star,\CK_\fo^\star)$, see~\eqref{eq.K.star}. It is immediate that
the lattice $\Fano(\fo^\star,\CK_\fo^\star)$
is geometric and it has no further geometric finite index extension
compatible with $\fo^\star$;
thus, there is a single deformation family.

Similar to \autoref{s.Golay}, introduce the sets
\[*
\CS^\star:=\bigl\{o\cap\fo^\star\bigm|o\in\CC_{24}\bigr\}\sminus\CC^\star
=\CS_3^\star\cup\CS_4^\star\cup\ldots
\cup\CS_{11}^\star\cup\CS_{12}^\star\cup\CS_{15}^\star.
\]
For any $3$-isotropic vector~$\iso$ that can be added to~$\fo^\star$,
the Sylvester test (\autoref{lem.Sylvester}) implies that
$\bigl|\|\iso\|_{\fo^\star}\bigr|\le2$, and only
$\|\iso\|_{\fo^\star}=\varnothing$ passes the kernel test
(\autoref{lem.kernel});
however, the partitioned pair
$(\fo^\star\cup\varnothing,\CK^\star)$ is not subgeometric.
Thus, any almost Kummer octic is
a triquadric.
Similarly, by  Lemmas \ref{lem.star}\iref{star.line.admissible}
and~\ref{lem.kernel}, if an extra line~$v$ can be added to~$\fo^\star$, then
\[*
\|v\|_{\fo^\star}\in\CS_3^\star\cup\ldots\cup\CS_7^\star\cup\CC,
\]
and only $\CS_3^\star\cup\ldots\cup\CS_6^\star$ pass the Sylvester test. Now,
adding to~$\fo^\star$ (in the progressive mode,
see \autoref{ss.step.r.general} below) up to four independent lines,
after saturating and sorting the results
we arrive at $107$ configurations
\qgen by lines.
Seven of them are Kummer (see \autoref{s.Kummer}), one is \config{QC.33}, and
the others have at most $29$ lines.
There remains to observe that $\rank\config{QC.33}=20$ and
Stage~$2$ of the proof (see \autoref{ss.stage.2}) is void.
\endproof

\subsection{Digression: Kummer quartics}\label{s.quartics}
Certainly, the arguments used above
apply to other degrees~$h^2$ as well. Without going into detail (to be
published elsewhere), we merely announce a few interesting
findings. First of all,
\roster*
\item
if $2d=2\bmod4$, then $\nodal(X)\le12$ for any birationally quasi-polarized
\KK-surface $(X,h)$ of degree $h^2=2d$;
\endroster
thus, (almost) Kummer \KK-surfaces may exist only in degrees
$0\bmod4$.
Henceforth, we confine ourself to the most interesting case of spatial
quartics.

\theorem\label{th.quartics}
There are but eight equilinear families of Kummer quartics $(X,h)$ with
$\NS(X)$ \qgen by lines. Among them,
there is one with $48$ lines and $4$ nodes.
\endtheorem

Concerning the last statement, recall that the maximal number of lines on a
quartic $X_4\subset\Cp3$ with $\Sing X_4\ne\varnothing$ is still an open
problem. The known upper bound is~$64$ (see~\cite{Veniani}). There is a
single example of a quartic with $52$ lines and two nodes (see
\cite{degt:singular.K3}), whereas all other known examples have at most $40$
lines.

Another open problem is the maximal number of lines in a triangle free
configuration of lines on a smooth quartic: the best bound is $52$
(see~\cite{DIS}), and the best example is $37$ (see~\cite{degt:lines}).
We have discovered a larger example.

\proposition
There exists a smooth almost Kummer quartic $X_4\subset\Cp3$ with a
quadrangular configuration of $39$ lines.
\endproposition

\section{Other triquadrics}\label{S.triquadrics}

In this section, we
consider $8$-polarized $3$-admissible graphs, \ie, we assume that $h^2=8$ and $m=3$.
In other words, we deal with triquadrics.
Here, we state a number of technical lemmas
which are used in the next section to derive \autoref{thm-main};
the \GAP~\cite{GAP4} aided proofs of these
lemmas (consisting mainly in feeding appropriate parameters to the master
algorithm and, therefore, hardly interesting from the mathematical point
of view)
are explained in the appendices.
We maintain the notation
introduced in \autoref{s.extendedversusplain} (\cf. \autoref{eq-graphex}).

\subsection{Locally elliptic configurations}\label{s.elliptic}
From~\eqref{eq.Pi.bound}, \eqref{eq.pencil},
and obvious combinatorial bounds on the number
of sections in a locally elliptic graph (\ie, the maximal number of
vertices that can be added to $\fiber$ without introducing a strictly smaller
parabolic subgraph,
see~\cite[Figure 1]{degt:lines}),
we conclude that, as in the smooth case,
\[
\ls|\graph|\le(18+3)+8=29\quad
\text{for any geometric locally elliptic graph~$\graph$},
\label{eq.locally.elliptic}
\]
\ie, any $\fiber$-graph with $\mu(\fiber)\ge5$.
This bound holds for any degree $h^2\ge4$. In fact, a simple computation with
$\tA_5$- and $\tD_5$-graphs confirms that, in degree~$8$,
\begin{equation} \label{eq.lellipticprecise}
\ls|\graph|\le23\quad
\text{except $\graphex\cong\LE_{25}^{}$, $\LE_{24}\A$, $\LE_{24}^{}$,
$\LE_{24}'$, $\LE_{24}''$, or $\LE_{24}^6$}.
\end{equation}
The first five exceptional configurations (the smooth ones) are introduced
in~\cite{degt:lines}, and
 $\Fano_8(\LE_{24}^6)\supset\Fano_8(\LE_{24}\A)$
is an index~$4$ extension 
with the same set of lines.


\subsection{Astral configurations}\label{s.astral}
A graph
of type $\tD_4$ is called \emph{astral}. In other words,
a graph~$\graph$ is astral if and only if $\girth\graph\ge6$ and $\graph$ has a
vertex~$v$ of $\val v\ge4$.

If $\pencil$ is a parabolic $\tD_4$-graph, then,
by~\eqref{eq.Pi.bound} and~\eqref{eq.pencil}, we
have $\ls|\pencil|\le22$.

\lemma[see \autoref{s.sec.D4}]\label{lem.astral}
Let $\graph\supset\fiber\cong\tD_4$ be an astral geometric configuration,
and
let the central vertex of $\fiber$ be one of the maximal valency
in~$\graph$. Then, one has
\[*
\ls|\sec\fiber|\le12\quad\text{and}\quad
\text{$\ls|\graph|<28$ whenever $\ls|\sec\fiber|\ge\M_\fiber:=11$}.
\]
\endlemma

\lemma[see \autoref{proofs.pencil}]\label{lem.d4p}
Let $\pencil$ be a $\tD_4$-pencil, $\ls|\pencil|\ge\M_\pencil:=18$, and
$\nodal(\pencil)\le14$. Then
\[*
\ls|\graph|<28
\]
for any geometric astral extension $\graph\supset\pencil$.
\endlemma

Thus, we have~\eqref{eq.bounds} and, hence,~\eqref{eq.goal}, \ie,
$\ls|\graph|<\M:=28$ for any astral graph~$\graph$
\emph{containing a $\tD_4$-pencil~$\pencil$ with $\nodal(\pencil)\le14$}.
However, there also is an astral Kummer configuration \config{SC.28'} (see
\autoref{prop.Kummer.256}), containing three pencils of type $4\tD_4$, and,
summarizing, we obtain a sharp bound
\[
\text{$\ls|\graph|\le28$ for any astral
	$3$-geometric $8$-polarized graph $\graph$}.
\label{eq.astral}
\]

\subsection{Pentagonal configurations}\label{s.pent}
A graph of type $\tA_4$ is called \emph{pentagonal}. In other words,
a graph~$\graph$ is pentagonal if and only if $\girth\graph=5$.

\lemma[see \autoref{s.sec.A4}]\label{lem.pent}
Let $\graph\supset\fiber\cong\tA_4$ be a pentagonal geometric configuration.
Then, one has
\[*
\ls|\sec\fiber|\le17\quad\text{and}\quad
\text{$\ls|\graph|<30$ whenever $\ls|\sec\fiber|\ge\M_\fiber:=14$}.
\]
\endlemma

\lemma[see \autoref{proofs.pencil}]\label{lem.a4p}
Let $\pencil$ be an $\tA_4$-pencil, $\ls|\pencil|\ge\M_\pencil:=17$, and
$\nodal(\pencil)\le14$. Then,
\[*
\text{$\graphex\cong\config{QF.30'}$, \config{QF.30''}, or \config{QF.30.5''}
 \rom(see~\cite{degt:lines} or \autoref{tab.mainsmall}\rom)}
\quad\text{or}\quad
\ls|\graph|<30
\]
for any geometric pentagonal extension $\graph\supset\pencil$.
\endlemma

The lattice $\Fano(\config{QF.30.5''})$
is an index~$3$ extension of
$\Fano(\config{QF.30''})$; the two lattices have the same set of lines and
differ by the exceptional divisors only.

Since Kummer or almost Kummer configurations (see~\autoref{S.Kummer}) are
never pentagonal, we conclude that
\[
\text{$\ls|\graph|\le30$ for any pentagonal
 $3$-geometric $8$-polarized graph $\graph$}.
\label{eq.pentagonal}
\]

\subsection{Quadrangular configurations}\label{S.quad}
A graph of type $\tA_3$ is called \emph{quadrangular}. Thus,
a graph~$\graph$ is quadrangular if and only if $\girth\graph=4$.

\lemma[see \autoref{s.sec.A3}]\label{lem.quad}
Let $\graph\supset\fiber\cong\tA_3$ be a quadrangular geometric configuration.
Then, one has
\[*
\ls|\sec\fiber|\le20\quad\text{and}\quad
\text{$\ls|\graph|<31$ whenever $\ls|\sec\fiber|\ge\M_\fiber:=16$}
\]
unless $\graphex$ is one of
the $13$ bi-colored
graphs
listed in \autoref{tab.main}
as the entries with a reference to \autoref{lem.quad}.
\endlemma

\lemma[see \autoref{proofs.pencil}]\label{lem.a3p}
Let $\pencil$ be an $\tA_3$-pencil, $\ls|\pencil|\ge\M_\pencil:=17$ and
$\nodal(\pencil)\le14$, and $\graph\supset\pencil$ a geometric quadrangular
extension. Then
\[*
\ls|\graph|<32
\]
unless $\graphex$ is one of the  $8$
bi-colored graphs listed in \autoref{tab.main}
as the entries with a reference to \autoref{lem.a3p}.
\endlemma

Observing that the Kummer configurations
\config{QC.32.K} (see \autoref{prop.Kummer.64}) and
\config{QC.32.4} (see \autoref{prop.Kummer.256}) and almost Kummer configuration
\config{QC.33} (see \autoref{prop.Kummer-}) appear among the
exceptions in
Lemmas~\ref{lem.quad} and~\ref{lem.a3p}, we conclude that, with the $15$
exceptions listed as $\QC_*$ in
\autoref{tab.main}, one has
\[
\text{$\ls|\graph|\le31$ for any quadrangular
 $3$-geometric $8$-polarized graph $\graph$}.
\label{eq.quadrangular}
\]


\section{Proof of \autoref{thm-main}}\label{S.Proof}

For
special octics the statement of theorem
is given by
\autoref{th.special} (see \autoref{proof.special}):
we obtain
a single configuration (\viz. \config{TC.33}) with $33$ lines on a
smooth \KK-surface.
Thus, we can assume that
$(X, \hp)$ is a triquadric (\ie, the $8$-polarized lattice $\NS(X) \ni \hp$
is $3$-admissible, see \autoref{th.K3.ext}),
and then
\autoref{th.K3.ext} reduces the proof to the classification
of the $8$-polarized bi-colored graphs~$\graph'$ such that
$$
\text{$\graph'$  is $3$-geometric and  $\ls|\spp_1\graph'|\geq32$}.
$$
We follow the general strategy described in \autoref{s.extendedversusplain}.

\subsection{{Stage 1:} the case of $\NS(X)$ \qgen by lines}
In this part of the proof
 we classify the extended saturations $\satex\graph$
assuming that $\graph$ is $3$-geometric and $\ls|\graph|\ge32$; by
\autoref{th.K3}, this corresponds to the triquadrics $(X,\hp)$
such that
the lattice $\NS(X)$ is \qgen by lines.
By~\eqref{eq-fano-must-be-hyperbolic}
and \autoref{lem.2stars}, $\graph$ is  a hyperbolic triangle- and biquadrangle free graph, \ie, it is a
$\fiber$-graph for a
certain affine Dynkin diagram \smash{$\fiber>\tA_2$} (see \autoref{subs.type}).
If $\graph$ is quadrangular (\ie, \smash{$\fiber\cong\tA_3$}),
then,
by \eqref{eq.quadrangular},
$\satex\graph$
is one  of the fifteen
$\QC_*$-graphs
listed in
 \autoref{tab.main}.
The other types are ruled out by
\roster*
\item
\eqref{eq.pentagonal}, if $\graph$ is pentagonal
(\ie, $\fiber\cong\tA_4$),
\item
\eqref{eq.astral}, if $\graph$ is astral
(\ie, $\fiber\cong\tD_4$), or
\item
\eqref{eq.locally.elliptic}, if $\graph$ is locally elliptic
(\ie, $\fiber>\tD_4$).
\endroster

\subsection{{Stage 2:} the general case\pdfstr{}{ \rm(see \autoref{ss.stage.2})}}\label{proof.general}
There
remains to analyze (extend by exceptional divisors) the seven graphs
$\graph:=\spp_1\graphex$ in
\autoref{tab.main} that have rank $\rank\graph<20$.
Since the threshold is $\ls|\bar\graph|\ge32$,
for the six configurations $\Theta^*_{32}$ it is the graph~$\graph$ itself
that is to be extended, and a direct
check shows that
the resulting graph is always one of those listed in \autoref{tab.main}.
(In fact, instead of $\graph$, we use smaller ``natural'' generating sets by
means of which the graphs were constructed in the proof, \cf.
\autoref{rem.generating.set}.)

The configuration
$\graph\cong\config{QC.34'}$ in \autoref{tab.main} needs
more work, as we have $\ls|\graph\sminus\bar\graph|\le2$.
We fix a certain
``natural'' generating set $\Lambda\subset\graph$ (\cf.
\autoref{rem.generating.set}), consisting of an \smash{$\tA_3$}-fiber and $14$
sections,
and run the extension algorithm to
check that there are no corank~$1$ extensions
$\bS\supset\Fano(\Lambda)=\Fano(\graph)$ satisfying
conditions~\iref{stage.2.rank}, \iref{stage.2.lines} in \autoref{ss.stage.2}.
Next, we observe that the action of the group $G:=\Aut\graph$ has three orbits,
$\graph=\bigcup_n\Omega_n$, where $n=\ls|\Omega_n|\in\{2,8,24\}$, and each
orbit has a representative $v\in\graph\sminus\Lambda$, so that the lattice
$\Fano(\graph\sminus v)$ has no extensions other than those of
$\Fano(\Lambda)$.
Finally, the action of~$G$ on the set of unordered pairs of vertices has $19$
orbits, $18$ of which also have representatives in $\graph\sminus\Lambda$.
The remaining orbit consists of a single pair $\{u,v\}=\Omega_2$,
and the graph $\graph\sminus\Omega_2\cong\config{QC.32'}$ has already been
analyzed.

\subsection{Graphs to octics}\label{s.graph2T}
Finally,
there remains to apply Nikulin's theory~\cite{Nikulin:forms} and,
for each of the sixteen Fano graphs~$\graph$ found,
classify the isomorphism
classes of primitive isometric embeddings $\Fano(\graph,\CK)\into\bL$,
$\CK\in\GK(\graph)$.
In each case, we obtain a single connected deformation family.
\qed

\remark\label{rem.uniqueness}
The uniqueness of a family of $K3$-surfaces representing each of the large
extended Fano graphs is a special property of these graphs; it needs to be
confirmed by a straightforward, albeit tedious computation
using \cite{Nikulin:forms}
on a case-by-case basis. For example, the
transcendental lattices corresponding to most maximal locally elliptic Fano
graphs (see \autoref{s.elliptic}) are not unique in their genera.

Occasionally, it may also happen that a graph has several
geometric extensions  (not in our tables) or that a ``deeper''
extensions introduces exceptional divisors while keeping the set of lines
(\cf. $\LE_{24}^6\supset\LE_{24}\A$ in \autoref{s.elliptic}
or $\config{QF.30.5''}\supset\config{QF.30''}$ in \autoref{s.pent}).
\endremark

\appendix

\section{The computation}\label{S.app}

We conclude the paper with a description of the algorithms used in the
proof of \autoref{thm-main}. As
explained in
\autoref{remark-no-code}, we do not present
a pseudocode  that
would be
a direct
translation of  \cite{degt.Rams8.anc}.
Instead, we describe, in purely mathematical
terms and \emph{in
the order convenient for the exposition},
what and how is to be computed, emphasizing but a few key points.
Technical details are left out.
These details/\penalty0
tricks
can be
deduced
from  \cite{degt.Rams8.anc} and
are mostly aimed at overcoming the limitations of the
contemporary hardware and the particular computer algebra
system~\cite{GAP4}
used; as such we find them
inappropriate for a
mathematical paper.
Besides, it is our sincere belief that, once the mathematical ideas are
clear, their implementation is straightforward, although tedious, essentially
boiling down to optimizing the parts that do not work fast enough.

In other words, our main objective in
this appendix
is
to explain the contents of \cite{degt.Rams8.anc}
in a human readable form, \ie, present
the main ideas of the algorithms rather than
their
implementation.
	We start with a more detailed exposition of
the general approach
developed in \autoref{s.approach}.

In a nutshell, we aim at considering \emph{all} graphs and
testing each graph
against the criteria given by Theorems~\ref{th.hyperelliptic}
and~\ref{th.octics} (see \autoref{s.aggraphs} for the definitions); these
tests are outlined in \autoref{s.tests}. (In particular,
this  approach requires
no proof of the correctness/termination of
our algorithm.) Our ultimate goal is
a reasonably small list of
graphs (see the saturation lists in \autoref{ss.Sat}), sufficiently large or
sufficiently ``singular'', that are  worth further investigation ``by hand''
(see \autoref{proof.special} or
\autoref{proof.general} and \autoref{s.graph2T} for an example of such
further analysis).

Clearly, the goal of considering all graphs and ruling out all but a few
	cannot be achieved by brute force,
so we have to use several tricks/ideas to overcome this difficulty. Thus, most
importantly,
\roster*
\item
we build graphs recursively, adding one vertex at a time
(see the graph extension procedure in \autoref{s.ext}) and, at
each step, excluding from the further consideration
graphs failing certain \emph{hereditary} tests (\cf.
\autoref{ss.extensibility} and \autoref{ss.geometric});
that is why heredity plays an important r\^{o}le in the
main text;
\item
we use a number of simple preliminary tests (see \autoref{s.prelim}),
that are quick,  at the beginning of each
step;
\item
still, the computation tends to diverge quite fast; therefore,
we try to start from relatively few sufficiently \emph{large} graphs.
It is this reason
why we use the dichotomy of~\eqref{eq.bounds},
see \autoref{S.Msigma} (in particular,
\autoref{rem.large.standard})
and \autoref{S:largepencils} (where we start from \emph{large} pencils);
\item
the graph extension process of \autoref{s.ext} stops as soon as we get a
graph of rank $20$, upon which only finite index extensions are considered
(\cf. \autoref{ss.Sat});
\item
we employ various \latin{a priori} bounds (derived from our taxonomy
or from the geometry of the problem) and symmetries of graphs involved (\cf.
\autoref{s.apriori}).
\endroster

In \autoref{S.basicalg}, we give a more detailed description of the basic
algorithms, tests,
\etc.\ used repeatedly throughout the computation.
Then, in Appendices \ref{S.Msigma} and
\ref{S:largepencils}, we work, on a case-by-case basis, towards,
respectively, the first and second
inequality in~\eqref{eq.bounds}: we discuss a few changes and tweaks and list
the parameters that were used in the actual computation.
These last two appendices
provide details of a machine aided proof for a number of lemmas stated in
\autoref{S.Proof}.

\section{Basic algorithms}\label{S.basicalg}

In this section, we describe the most basic algorithms applicable to any
degree $h^2=2d$. Essentially, they are those used in~\cite{degt:lines}, with
a few minor modifications adjusting them to the case of non-empty singular
locus. Most notably, the ``smooth'' requirement $\root_0(S,h)=\varnothing$ is
replaced with \autoref{lem.compatible}.

Note that, most of the time, we work with \emph{plain} Fano graphs, even
though computing such a graph starts with computing a Weyl chamber (\cf.
\autoref{ss.extensibility} below). At the very end, when computing the
saturation lists of the large plain graphs (see \autoref{ss.Sat} below), we keep
the record of the exceptional divisors, thus obtaining
the extended Fano graphs that appear in the final statements.

\subsection{The master test}\label{s.tests}
In the heart of all algorithms
is a procedure detecting if a given
graph~$\graph$ is geometric. More precisely, the input consists of a
polarized lattice $S\ni h$ and distinguished subset $\graph\subset\root_1(S,h)$.
Typically, $S$ is of the form $\Fano(\graph,\CK)$; however, occasionally we
take for~$S$ an extension of this lattice by an $m$-isotropic vector
and/\penalty0or a
few ``potential'' exceptional divisors.

\subsubsection{Extensibility and admissibility}\label{ss.extensibility}
We assume~$S$ given by its Gram matrix in a certain basis
$\{h,b_1,b_2,\ldots\}$ containing~$h$. Then, we consider the \emph{rational}
lattice
\[*
h^\perp_\Q:=\sum\Z\biggl(b_k-\frac{b_k\cdot h}{2d}h\biggr),
\]
multiply the form by $-2d$, and use \GAP's~\cite{GAP4} function
\texttt{ShortestVectors} to compute
\[*
\CV_s:=\bigl\{v\in h^\perp_\Q(-2d)\bigm|v^2=s\bigr\}\quad
 \text{for all $1\le s\le4d+1$}.
\]
Given an admissibility level $m\in\{1,2,3\}$, we check if there is
\roster*
\item
$v\in\CV_{m^2}$ such that the $m$-isotropic vector $v+m(2d)\1h$
is in $S$;
\endroster
if such a vector is found,
the algorithm terminates as $(S,h)$ is not $m$-admissible.

Next, we compute the set
\roster*
\item
$\root_0(S,h)=\CV_{4d}\cap S$.
\endroster
If there is a separating (with respect to~$\graph$ or $\graph$ and
a given isotropic vector) root $r\in\root_0(S,h)$, the algorithm terminates
as $\graph$ is not extensible (by Lemma~\ref{lem.ext} or
\ref{lem.ext.partitioned}).
Otherwise, we compute a
Weyl chamber $\Weyl'$ for $S':=(\Z h+\Z\graph)^\perp\subset S$
(\cf. \autoref{lem.compatible}) using any generic functional $\Ga\:S'\to\R$,
and use \autoref{alg.ext} to extend~$\Weyl'$ to the unique Weyl
chamber~$\Weyl$ compatible with~$\graph$.

\remark
On the few occasions where $\root_0S'\ne\varnothing$,
a generic functional is found as follows. Pick any root $r\in\root_0S'$ and
let $\Ga\:v\mapsto-(v\cdot r)$. Whenever we can find a root
$r\in\root_0S'$ such that $\Ga(r)=0$, we change~$\Ga$ to
$v\mapsto2\Ga(v)-(v\cdot r)$, continuing this process until $\Ga$ is generic.
This works since, in an even lattice, we have $\ls|r_1\cdot r_2|\le1$ for any
two roots $r_1\ne\pm r_2$.
\endremark

There remains to compute the set
\roster*
\item
$\root_1(S,h)=\bigl\{v+(2d)\1h\bigm|v\in\CV_{2d+1}\}\cap S$
\endroster
and Fano graphs
$\Fn_\Weyl(S,h)$ and $\Fex_\Weyl(S,h)$
(directly as explained in \autoref{s.polarized}).

\subsubsection{Detecting subgeometric sets}\label{ss.geometric}
To check if a given lattice $S\ni h$ is geometric, we compute the
discriminant group $\discr S$ and apply \cite[\thmcit~1.12.2]{Nikulin:forms}
(see also \cite[\thmcit~3.2]{DIS}). If the result is negative, we proceed as
follows:
\roster*
\item
list all isotropic ($\Ga^2=0\bmod2\Z$) vectors $\Ga\in\discr S$
\emph{of prime order};
\item
for each~$\Ga$, compute the finite index extension $S_\Ga\supset S$ and check
whether it is admissible (\via\ \autoref{ss.extensibility},
using the same polarization~$h$ and subset~$\graph$);
\item
if successful, repeat \autoref{ss.geometric}
(this algorithm) for~$S_\Ga$.
\endroster
The algorithm terminates as soon as a geometric lattice has been found or all
admissible finite index extensions have been tried.

\remark
Primitive as it is, this algorithm serves our needs as the discriminant
groups $\discr S$
are usually reasonably small.
We do use a couple of tricks to speed up the computation:
\roster*
\item
if a large subgroup $G\subset\OG_h(S)$ preserving~$\graph$ is known
(\cf. \autoref{s.extension.algorithm} below), we use a single
representative of each $G$-orbit of isotropic vectors;
\item
at the first step (for~$S$ itself), prime order isotropic vectors~$\Ga$
resulting in non-admissible extensions~$S_\Ga$ are recorded not to be used
again.
\endroster
\endremark

For statistical purpose, whenever we establish that a
graph~$\graph$ is geometric, we automatically record the counts
$(\ls|\graph|,\ls|\Sing\graph|)$ (but not $\graph$ itself).

\subsubsection{The saturation lists}\label{ss.Sat}
The \emph{saturation list} of a pair~$(\graph,\CK)$ is the
set
\[*
\anchor{Sat}
{\Sat_m(\graph,\CK)}:=\bigl\{\satex(\graph,\CK')\bigm|
 \text{$\CK'\in\GK^m(\graph)$ and $\CK'\supset\CK$}\bigr\}.
\]
If $\CK=0$, it is omitted from the notation.
If a certain \emph{global criterion}
(not necessarily hereditary)
$\Num\:\{\text{graphs}\}\to\{\false,\true\}$ is given, we denote
\[*
\Sat_m(\graph,\CK;\Num):=\bigl\{\graph'\in\Sat_m(\graph,\CK)\bigm|
 \Num(\graph')=\true\bigr\}.
\]
Typically, $\Num$ consists of a fixed type $\fiber$ in the taxonomy of graphs
and a certain numeric bound, \eg, $\ls|\graph|\ge M$.
In fact, in order to construct plenty of examples,
we also retain all  $m$-geometric graphs~$\graph$ satisfying
\[*
\ls|\graph|\ge\Mlines\quad\text{or}\quad\ls|\Sing\graph|\ge\Msing:=4,
\]
where $\Mlines$ depends on the kind of the graphs considered.

The saturation lists are computed by the same algorithm as in
\autoref{ss.geometric}
(applied to the lattice $S:=\Fano(\graph,\CK)$ and $\graph$ itself
as the distinguished subset),
except that we do not stop on the first hit, listing all geometric
extensions.
These lists are not intended for further processing (except sorting), and
it is here that we keep track of the exceptional divisors, obtaining extended
Fano graphs for the statements.


\subsection{Preliminary tests}\label{s.prelim}
The algorithms in \autoref{s.tests} are relatively expensive
(mainly, due to the \texttt{ShortestVectors}); for this reason, they are
usually preceded by a few preliminary test ruling out
the vast majority of the possibilities.

As explained in \autoref{s.ext} below, typically we extend a given
graph~$\graph_0$ or pair $(\graph_0,\CK)$ by a one or several
(pseudo-)vertices~$v$ described by means of their supports $\|v\|$. The number of
possibilities for a single extra vertex (essentially, a subset of~$\graph_0$)
is huge, even when restricted by a statement like \autoref{lem.star}. Most of
them are ruled out by the following two obvious tests,
 for which the bulk of
the computation
(\ie~the computation of the inverse Gram matrix)
 depends on $(\graph_0,\CK)$ only and can be done
once.
Namely,
recall that we speak of the support of a (pseudo-)vertex $v$
(\cf. \autoref{def.pseudo})
only within the range of applicability of
Lemmas~\ref{lem.simply-laced} and~\ref{lem.e}, so that
$\|v\|$
determines the projection
\[*
\Fano(\graph_0)+\Z v\to\Fano(\graph_0)\dual,\quad v\mapsto v^*,
\]
and, hence, the one-vector extension
$\Fano(\graph_0)+\Z v$ itself.
(The projection $v^*$ is easily computed
in terms of $\|v\|$ and the inverse Gram matrix of the \emph{original}
graph~$\graph_0$.)

\lemma[the kernel test]\label{lem.kernel}
An extra \rom(pseudo-\rom)vertex~$v$ cannot be added
to a pair $(\graph_0,\CK)$ if $v^*\!\cdot k\notin\Z$ for at least one
element $k\in\CK$.
\done
\endlemma

\lemma[the Sylvester test]\label{lem.Sylvester}
A \rom(pseudo-\rom)vertex~$v$ cannot be added
to a graph $\graph_0$ if $v^2>(v^*)^2$.
If $v^2=(v^*)^2$,
the addition of~$v$ does not increase the rank.
\done
\endlemma

Besides, usually we fix a certain type~$\fiber$ and assume that the
\emph{explicit part} (see \autoref{wrn.sat}) of each new graph is a
$\fiber$-graph (or a similar condition like lack of biquadrangles \etc).
As a rule, we run a few quick tests directly in terms of the supports of the
vertices added (\cf. \autoref{rem.extra} below). Then, the Gram matrix of the
new graph needs to be computed, and, before passing this matrix to
\autoref{s.tests}, we check more thoroughly that it defines a
$\fiber$-graph.

\subsection{Graph extensions}\label{s.ext}
In most algorithms,
we fix a certain \emph{base graph}~$\graph_0$ and
consider its \emph{extensions} $\graph\supset\graph_0$ by a few extra
vertices. In this construction, each extra vertex $v\in\graph\sminus\graph_0$
can be represented by and, henceforth, is \emph{identified} with its support
$\supp_{\graph_0}v$.
Therefore, we
merely
regard $\bv:=\graph\sminus\graph_0$ as a \emph{multiset}
(as we do not assert that the correspondence is injective) of
subsets of~$\graph_0$.
(This approach, \viz. handling extra lines by means of their supports, is
illustrated in the proofs in \autoref{S.Kummer}.)
We fix an ordering
of~$\graph_0$ and assume each extra vertex $v\subset\graph_0$ ordered and
each multiset~$\bv$ of extra vertices
ordered lexicographically, $\bv=\{v_1\le\ldots\le v_r\}$.
For a set $S$ and $n\in\NN$, we denote by
$\anchor{sym}S\sym{n}$ the $n$-th
symmetric power of~$S$ and by $\anchor{Com}C(S,n)$,
the set of all $n$-combinations of~$S$; then, we abbreviate
$C_*(S,n):=\bigcup_{i=0}^n\Com(S,n)$.

Given a multiset~$\bv$ as above, we use the notation
\[
\graph:=\graph_0\gcup\bv\quad\text{or}\quad
\graph:=\graph_0\gcup\bv(\gram)
\label{eq.graph.extension}
\]
for the set theoretic union $\graph_0\cup\bv$ equipped with an extra edge
connecting $u\in\graph_0$ and $v\in\bv$ whenever $u\in v$. In the former case,
$\bv$ itself is regarded as a discrete graph, whereas in the latter case,
the graph structure
on~$\bv$ is an extra piece of data
given by an adjacency matrix $\gram=[m_{ij}]\in\Sym(\ls|\bv|,\F_2)$,
$m_{ij}=v_i\cdot v_j\bmod2$.

A \emph{$\graph_0$-isomorphism} between two graph extensions
$\graph'\supset\graph_0$
and $\graph''\supset\graph_0$ is a graph isomorphism $\graph'\to\graph''$
taking~$\graph_0$ to~$\graph_0$ \emph{as a set}.
The group of $\graph_0$-automorphisms of
$\graph\supset\graph_0$ is denoted by $\Aut(\graph,\graph_0)$.
The \emph{explicit sorting} of a list of graph
extensions is the procedure removing all but one representative of each
$\graph_0$-isomorphism class. In the special case where $\graph_0=\varnothing$, this
procedure is called the \emph{ultimate sorting}.
We use the \texttt{GRAPE} package
\cite{GRAPE:nauty,GRAPE:paper,GRAPE} in \GAP~\cite{GAP4}
(with a few minor performance enhancements),
computing also, as a
by-product, the groups $\Aut(\graph,\graph_0)$.

\subsection{The extension algorithm}\label{s.extension.algorithm}
In this section, we describe a procedure which is the
essential part of all other algorithms considered below.

\subsubsection{The input}\label{ss.input}
Since, usually, the algorithm is part of a larger computation, we always have
a certain \emph{global criterion}
$\Num\:\{\text{graphs}\}\to\{\false,\true\}$ in mind. In addition,
the input consists of the following data:
\roster*
\item
an affine Dynkin diagram~$\fiber$,
\item
an \emph{admissibility level} $1\le m\le3$,
\item
a choice of the \emph{mode},
see \autoref{ss.step.r.discrete} \vs. \autoref{ss.step.r.general} below,
\item
a local \emph{numeric criterion} $\num\:\{\text{graphs}\}\to\{\false,\true\}$,
\item
a base $\fiber$-graph~$\graph_0\supset\fiber$
and a \emph{symmetry group} $G_0\subset\Aut(\graph_0,\fiber)$, and
\item
a $G_0$-invariant initial set~$\CS_1$ of extra vertices (as subsets of~$\graph_0$).
\endroster
The goal is finding all extensions $\graph\supset\graph_0$
as in~\eqref{eq.graph.extension}, with
$\bv:=\graph\sminus\graph_0\in\CS_1\sym{r}$, $r\in\NN$,
satisfying, at least, the following conditions:
\roster
\item\label{ext.fiber}
the graph~$\graph$ \emph{itself} is a hyperbolic $\fiber$-graph, and
\item\label{ext.geometric}
one has $\GK^m(\graph)\ne\varnothing$.
\endroster
Sometimes (\eg, if $\bv$ is not assumed discrete), we also insist that
\roster[\lastitem]
\item\label{ext.rank}
the rank of~$\graph$ is as large as possible:
$\rank\graph=\rank\graph_0+\ls|\graph\sminus\graph_0|$.
\endroster
We proceed step by step, with formal Step~$0$ returning
$\{\graph_0\}$ and $\bCS_0:=\{\varnothing\}$.


\subsubsection{Step~$1$}\label{ss.step.1}
We compute the set
\[*
\bCS_1:=\bigl\{v\in\CS_1\bigm|
 \text{$\graph:=\graph_0\gcup\{v\}$
        satisfies~\autoref{ss.input}\iref{ext.fiber}, \iref{ext.geometric} and
        $\rank\graph<20$}\bigr\}.
\]
(For the last condition $\rank\graph<20$, see \autoref{conv.max.rank} below.)
To this end, we
\roster
\item
run the Sylvester test (\autoref{lem.Sylvester}),
reducing $\CS_1$ to a subset $\CS_1'$;
\item
pick a representative~$v$ of each $G_0$-orbit on~$\CS_1'$ and let
$\bv:=\{v\}$;
\item\label{r.prelim}
for each~$\bv$, run appropriate \emph{preliminary tests} (to be specified
below);
\item\label{r.fiber}
for each~$\bv$ left, check that $\graph:=\graph_0\gcup\bv$ is a $\fiber$-graph;
\item\label{r.main}
for each~$\graph$ left, run the master test (\autoref{s.tests})
to check \autoref{ss.input}\iref{ext.geometric}, upon which we discard all
graphs~$\graph$ of rank~$20$ (see \autoref{conv.max.rank} below).
\endroster

\remark\label{rem.extra}
Typical preliminary tests used (for speed) in item~\iref{r.prelim}
are as follows (where $1\le p<q<r\le\ls|\bv|$):
\roster*
\item
if $\fiber>\tA_2$, then $u_1\cdot u_2=0$ for $u_1,u_2\in v_p$, $u_1\ne u_2$
(no triangles);
\item
if $\fiber>\tA_3$, then $\ls|v_p\cap v_q|\le1$ (no quadrangles);
\item
if $m=3$, then $\ls|v_p\cap v_q|\le2$ and
$\ls|v_p\cap v_q\cap v_r|\le1$ (no biquadrangles,
see \autoref{lem.2stars}).
\endroster
If applicable (and not covered by any of the above), we can also
use~\eqref{eq.sec.bound} below,
checking appropriate $\bnd_{**}$-fold intersections of the vertices to be added.
More subtle combinatorial tests are usually ignored as they are
incorporated to item~\iref{r.fiber}, where we thoroughly check the full
Gram matrix of $\Fano(\graph)$.
\endremark

\subsubsection{Step $r\ge2$, the safe mode}\label{ss.step.r.discrete}
In this version of the algorithm, we assume that the graph~$\bv$ is discrete.
The goal is computing the set
\[*
\bCS_r:=\bigl\{\bv\in\bCS_1\sym{r}\bigm|
 \text{$\graph:=\graph_0\gcup\bv$
        satisfies~\autoref{ss.input}\iref{ext.fiber}, \iref{ext.geometric} and
        $\rank\graph<20$}\bigr\}.
\]
We use the result of Step~$(r-1)$ and start with the set of
ordered $r$-tuples
\[
\CS_r:=\bigl\{\bv\in\bCS_{r-1}\times\bCS_1
\bigm|v_1\le\ldots\le v_{r-1}\le v_r\bigr\}.
\label{eq.Sr}
\]
Then, if $r=2$, we let $\CS_2':=\CS_2$; otherwise, we reduce $\CS_r$ to
\[*
\CS_r':=\bigl\{\bv\in\CS_r\bigm|
 \text{$(\ldots,\hat v_i,\ldots)\in\bCS_{r-1}$ for each $i<r$}\bigr\};
\]
apart from reducing the overcounting,
this \emph{must} be done to ensure a well-defined action of~$G_0$.
Finally, we pick a representative~$\bv$ of each $G_0$-orbit on~$\CS_r'$
and repeat operations~\iref{r.prelim}--\iref{r.main} in \autoref{ss.step.1}.

\subsubsection{Step $r\ge2$, the progressive mode}\label{ss.step.r.general}
If $\bv$ is not assumed discrete, we insist that each step of
the algorithm should increase the rank. Thus,
we modify the output of Step~$1$ to
\[*
\bCS_1':=\bigl\{v\in\bCS_1\bigm|\rank(\graph_0\gcup\{v\})>\rank\graph_0\bigr\},
\]
compute the set
\begin{multline*}
\bCS_r^*:=\bigl\{(\bv,\gram)\in\bCS_1'\sym{r}\times\Sym(r,\F_2)\bigm|\\
 \text{$\graph(\gram):=\graph_0\gcup\bv(\gram)$
        satisfies~\autoref{ss.input}\iref{ext.fiber}--\iref{ext.rank} and
        $\rank\graph(\gram)<20$}\bigr\},
\end{multline*}
and define $\bCS_r$ as the projection of $\bCS_r^*$ to the first
factor.
We proceed as in \autoref{ss.step.r.discrete}, with a few alterations:
\roster*
\item
starting from item~\iref{r.prelim}, we deal with pairs $(\bv,\gram)$ rather
than sets~$\bv$
(in the hope that items~\iref{r.prelim} and~\iref{r.fiber} would rule out
most matrices~$\gram$), and
\item
in item~\iref{r.main}, we check, in addition, that
$\rank\graph(\gram)=\rank\graph_0+r$.
\endroster

In both modes,
the algorithm terminates as soon as $\bCS_r=\varnothing$.
Sometimes, we also terminate it after a preset number~$r\smax$ of steps
given as part of the input.

\convention[graphs of the maximal rank]\label{conv.max.rank}
Since any geometric extension of any lattice $\Fano(\graph)$ of
rank~$20$
is of finite index, \emph{we systematically discard any
graph~$\graph$ of the maximal rank $\rank\graph=20$};
instead, we just store the set $\Sat_m(\graph;\Num)$ in a separate
list~$\maxlist$.
(Clearly,
it suffices to compute the saturation lists only for the extensions
$\graph\supset\graph_0$ satisfying
$\rank\graph=\rank\graph_0+\ls|\graph\sminus\graph_0|$.)
\endconvention

\subsubsection{Plain \vs. saturated output}\label{ss.output}
The \emph{plain output} of the algorithm (usually, it is intended for further
processing) is just the union of the outputs, \ie, graphs
\[
\graph:=\graph_0\gcup\bv,\ \bv\in\bCS_r,\quad\text{or}\quad
\graph:=\graph_0\gcup\bv(\gram),\ (\bv,\gram)\in\bCS_r^*,
\label{eq.output}
\]
of all steps, filtered \via~$\num$:
only the graphs~$\graph$ with $\num(\graph)=\true$ are retained.
The content of the list~$\maxlist$ (see \autoref{conv.max.rank}) is \emph{not}
part of the output, although it \emph{is} taken into account when drawing the
final conclusion.

The \emph{saturated output} (usually intended as the final result) is the
union of the list~$\maxlist$ and
the saturations $\Sat_m(\graph;\Num)$ over all graphs~$\graph$
as in~\eqref{eq.output} and such that $\rank\graph=\rank\graph_0+\ls|\bv|$.
(In the actual implementation, searching for examples,
we compute the saturated output in any case,
even if only plain output is needed for further processing.)

\subsubsection{Two-phase computation}\label{ss.2phase}
In a few cases, when the group~$G_0$ and, hence, intermediate lists $\CS_r$
are too large, we have to break the algorithm into two separate phases.
Namely, we run up to a certain preset number $r_\text{break}$ of steps, upon
which start over and apply the same algorithm to each graph~$\graph$
in the output. To reduce the overcounting, we modify the initial set for
phase~$2$ as follows:
\roster*
\item
start with the set $\bCS_1$ computed in phase~$1$, and
\item
assuming that shorter extra vertices are added first, reduce this set (for
each graph $\graph\supset\graph_0$) to
$\bigl\{u\in\bCS_1\bigm|
 \text{$\ls|u|\ge\ls|v|$ for all $v\in\graph\sminus\graph_0$}\bigr\}$.
\endroster
If plain output is required, the combined output of
phase~$2$ is subject to the explicit sorting (see \autoref{s.ext}).

\subsection{\latin{A priori} bounds and patterns}\label{s.apriori}
We fix an order $\fiber=\{a_1,\ldots,a_N\}$ of each affine Dynkin
diagram~$\fiber$ to be considered. Usually, for $\fiber$-graphs, we have
certain \latin{a priori} bounds
$b_i$, $b_{ij}$, $b_{i0}$, $b_{0i}$, $1\le i,j\le N$, so that
\[
\alignedat2
\#\bigl\{l\in\sec a_i\bigm|l\cdot s=1\bigr\}&\le b_{ij}&\quad&
 \text{for any $s\in\sec a_j$},\\
\#\bigl\{l\in\sec a_i\bigm|l\cdot s=1\bigr\}&\le b_{i0}&\quad&
 \text{for any $s\in\sec\varnothing$},\\
\#\bigl\{l\in\sec\varnothing\bigm|l\cdot s=1\bigr\}&\le b_{0i}&\quad&
 \text{for any $s\in\sec a_i$},\\
\ls|\sec a_i|&\le b_i,
\endalignedat
\label{eq.sec.bound}
\]
where $\sec\varnothing:=\sec_\graph\varnothing$
stands for the set of lines $l\in\graph\supset\fiber$
disjoint from~$\fiber$.
(For $b_i$, we usually take $b_i=v\smax(\fiber)-\val_\fiber a_i$, where
$v\smax(\fiber)$ is a bound for the maximal valency of a vertex in an
$m$-admissible $\fiber$-graph.)

A \emph{$\fiber$-pattern}, or just \emph{pattern},
is a function $\pi\:\fiber\to\NN$ such that
$\pi(a_i)\le\bnd_i$ for
each $i\le N$. The \emph{size} of a pattern~$\pi$ is
$\ls|\pi|:=\sum\pi(a)$, $a\in\fiber$.
The group $\Aut\fiber$ acts on the set of patterns, and we
denote by $\anchor{pat}{\pat(\fiber)}$ the set of the \emph{lexicographically maximal}
representatives of the $(\Aut\fiber)$-orbits.
Clearly, we can use the symmetry and confine ourselves to the
$\fiber$-graphs $\graph\supset\fiber$ such that the
\emph{section count}
$\anchor{scount}\pi_\graph\:a\mapsto\ls|\sec a|$, $a\in\fiber$ is an element of
$\pat(\fiber)$.

Given a function $\rho\:D\to\NN$,
$D:=\operatorname{domain}(\rho)\subset\fiber$, a subset
$S\subset\fiber\sminus D$, and a bound $M\in\NN$, we define the \emph{range}
associated with these data as
\[*
\anchor{range}
{\range_S(\rho,M)}:=\bigl\{{\textstyle\sum_{a\in S}\pi(a)}\bigm|
 \pi\in\pat(\fiber),\ \pi|_D=\rho,\ \ls|\pi|\ge M\bigr\}.
\]

\section{Sections at a single fiber~\pdfstr{Sigma}{$\fiber$}}\label{S.Msigma}

In this section we present the proofs of several versions of the first inequality in
\eqref{eq.bounds}, where parameters
 $M$ and $M_\fiber$ are set on a case-by-case
basis, depending on the type  $\fiber$ of the graph $\graph$.
We rely upon \autoref{rem.large.standard} below that lets us start the
computation from relatively few sufficiently large ``standard'' graphs.

Thus, we start with a \emph{fiber} (affine Dynkin diagram)
$\fiber:=\{a_1,\ldots,a_N\}$ and set a \emph{goal}
$\ls|\sec\fiber|\ge M_\fiber$,
taking for the global criterion in \autoref{ss.Sat}
\[*
\Num(\graph):=\bigl(\text{$\graph$ is a $\fiber$-graph}\bigr)\mathbin\&
 \bigl(\ls|\sec_\graph\fiber|\ge M_\fiber\bigr).
\]
We make the following assumptions
(proved separately in the main text):
\roster
\item\label{sec.simple}
each section to be added intersects exactly one of $a_1,\ldots,a_N$;
\item\label{sec.disjoint}
for each $1\le i\le N$, all sections $s\in\sec a_i$ are pairwise disjoint.
\endroster
The algorithm
computing large sets of sections at~$\fiber$
has $N$ levels, one for each line $a_i\in\fiber$,
starting from the graph $\graph_0:=\Sigma$.
Each subsequent level~$k$ produces a list $\{\graph_{k}\}$ from each
graph~$\graph_{k-1}$ obtained at level~$(k-1)$, and
the full output is the union of these lists.

\subsection{Level \pdfstr{k>=1}{$k\ge1$} of the algorithm}\label{s.fiber.k}
For the input, we fix a graph
\[*
\graph_{k-1}:=\fiber\cup\sec a_1\cup\ldots\cup\sec a_{k-1}
\]
with a known automorphism group~$G_{k-1}$ fixing $\fiber$ \emph{pointwise}.
Consider the function
$\rho\:\{a_1,\ldots,a_{k-1}\}\to\NN$, $a\mapsto\ls|\sec a|$,
and let $R_k:=\range_k(\rho,M_\fiber)$.
We run up to $r\smax:=\max R_k$ steps of the algorithm
in~\autoref{s.extension.algorithm}. In addition to the given fiber~$\fiber$
and admissibility level~$m$, we choose the safe mode of the algorithm
(see \autoref{ss.step.r.discrete}),
take
$\num(\graph):=\bigl(\ls|\sec_\graph a_{k}|\in R_k\bigr)$
for the numeric
criterion, $\graph_{k-1}$ and~$G_{k-1}$ for the base graph and symmetry group,
respectively, and
\[
\CS_1:=\{a_{k}\}\times\Sec_k(\graph_{k-1}),\qquad
\Sec_k(\graph):=\prod_{i=1}^{N}\Coms(\sec_\graph a_i,\bnd_{ik}),
\label{eq.Sec}
\]
see~\eqref{eq.sec.bound}, for the initial set of vertices.
The result is the plain output, see~\autoref{ss.output}.

\remark\label{rem.large.standard}
Since we assume that each graph $\sec a_i$ is discrete, level~$1$ of the
algorithm is ``trivial'': we merely test the $1$-parameter family
$\graph(n):=\fiber\cup\sec a_1$, where
$\ls|\sec a_1|=n\le\bnd_1$, see \autoref{s.apriori}.

Likewise, if
$\bnd_{12}=\bnd_{21}\le1$, level~$2$ is reduced to testing the $3$-parameter family of
graphs
$\graph(n_1,n_2;r):=\fiber\cup\sec a_1\cup\sec a_2$, where we let
$\ls|\sec a_i|=n_i$, $i=1,2$, and
$r\le\min\{n_1,n_2\}$ is the number of
sections $s_{1i}\in\sec a_1$ intersecting a section $s_{2j}\in\sec a_2$;
when constructing the Gram matrix, we can assume that
$s_{1i}\cdot s_{2i}=1$ for $1\le i\le r$ whereas all other pairs
of sections are disjoint.

We can take this observation two steps further:
if $\fiber=\tA_m$, $m=4,5$,
levels $1$ to $(m-1)$ reduce to testing an $m$-parameter family
$\graph(n_1,\ldots,n_{m-1};r)$. Thus, in fact, we always use at most two
``essential'' levels.
\endremark

The output of the $N$-th level
is a complete list of the
subgeometric
graphs of the form $\graph_N:=\fiber\cup\sec\fiber$ satisfying
$\ls|\sec\fiber|\ge M_\fiber$. (Note that, at this point, we do \emph{not}
require that $\graph_N$ should be a $\fiber$-graph.)
In particular, we obtain an upper bound on
$\ls|\sec\fiber|$. To eliminate the overcounting and compute the
groups $\Aut(\graph_N,\fiber)$, we apply the explicit sorting (see
\autoref{s.ext}).


\subsection{Level~$0$: fiber components}\label{s.level.0}
We fix another threshold $M>M_\fiber+\ls|\fiber|$ and try to list all
$m$-geometric  $\fiber$-graphs $\graph\supset\graph_N$ satisfying,
in addition, the inequality
$\ls|\graph|\ge M$.
For this, we run another instance of the algorithm in
\autoref{s.extension.algorithm}, setting
\[*
\num(\graph):=\bigl(\ls|\graph|\ge M\bigr)\quad\text{and}\quad
\Num(\graph):=\bigl(\text{$\graph$ is a $\fiber$-graph}\bigr)\mathbin\&
 \num(\graph)
\]
We take one of the graphs $\graph_N$ and groups $G_N:=\Aut(\graph_N,\fiber)$
for the base graph and symmetry group, respectively,
and
\[*
\CS_1:=\Coms(\sec a_1,\bnd_{10})\times\ldots\times
 \Coms(\sec a_N,\bnd_{N0}),
\]
see~\eqref{eq.sec.bound}, for the initial set of vertices.
The preliminary tests are as in Remarks~\ref{rem.extra}, and
we are interested in the saturated output,
see~\autoref{ss.output}.

The choice of the mode
of the algorithm depends on the type
of~$\fiber$. However, in the most common case $\rank\graph_N=19$,
we choose the \emph{progressive mode}
(see \autoref{ss.step.r.general}), adding exactly one extra line to obtain
(and discard) a graph of rank~$20$.

\remark[the validity test]\label{rem.validity.0}
For most fibers~$\fiber$, the safe mode of the algorithm
(see \autoref{ss.step.r.discrete}) is chosen.
It runs faster, but
the validity of this choice, \ie, the
fact that any \emph{sufficiently large}
extension $\graph\supset\graph_N$ is spanned over~$\graph_N$ by
a collection of \emph{pairwise disjoint}
fiber components, needs justification on a case-by-case basis. This is done
automatically, using the known list of large combinatorial
$\fiber$-pencils (see \autoref{S:largepencils} below).
Namely, for each $\delta\in\NN$, we compute
\[*
r\smax(\delta):=\min\bigl\{\nodal(\pencil\sminus\fiber)\bigm|
 \text{$\pencil$ is a $\fiber$-pencil,
        $\ls|\pencil\sminus\fiber|+\delta\ge M$}\bigr\}.
\]
Then, if the
algorithm
starting from a graph $\graph_N:=\fiber\cup\sec\fiber$
does not terminate in $r\smax(\ls|\graph_N|)$ steps, an error is
signaled and $\graph_N$ is returned as unsettled.
It is this criterion that dictates our choice of the
thresholds $M$ and $M_\fiber$.
\endremark

\subsection{Sections at \pdfstr{Sigma = D\_4}{$\fiber=\tD_4$}}\label{s.sec.D4}
We order~$\fiber$ so that $a_1$ is the ``central'' $4$-valent vertex.
Since any $\tD_4$-graph is pentagon free, \emph{all} sections are pairwise
disjoint; thus, in~\eqref{eq.sec.bound} we have, for $1\le i\ne j\le5$,
\[
\bnd_{i0}=1,\qquad \bnd_{ij}=0.
\label{eq.bound.D4}
\]
The other bounds  $\bnd_{01}=6$ and $\bnd_{0k}=5$ for $k\ge2$ are not
used.

Since $\Aut\tD_4=\Bbb{S}(a_2,\ldots,a_5)$
(the full symmetric group), and assuming
that $a_1$ is chosen to have the absolute maximal valency in the graph, we
have the following description
(assuming $m=3$, so that the maximal valency of a vertex
is~$7$, see \autoref{lem.star}\iref{star.line.admissible})
\[
\pi\in\pat(\fiber)\quad\text{iff}\quad
6\ge\pi(a_1)+3\ge\pi(a_2)\ge\ldots\ge\pi(a_5).
\label{eq.pat.D4}
\]

\proof[Proof of \autoref{lem.astral}]
We use the thresholds $M_\fiber=11$, $M=28$, $\Mlines=26$.
At level~$0$,
we choose the safe mode and have to add up to four
pairwise disjoint extra
lines disjoint from~$\fiber$, justifying the validity as in
\autoref{rem.validity.0}.
\endproof

\subsection{Sections at \pdfstr{Sigma = A\_4}{$\fiber=\tA_4$}}\label{s.sec.A4}
We order $\fiber$ cyclically.
Since an $\tA_4$-graph is quadrangle free,
in~\eqref{eq.sec.bound} we have, for $1\le i\ne j\le5$,
\[
\bnd_{i0}=1,\quad \bnd_{ij}=1\ \text{if $i=j\pm2\bmod5$},\quad
\bnd_{ij}=0\text{ otherwise}.
\label{eq.bound.A4}
\]
The other bound $\bnd_{0i}=6$ (assuming the graph $3$-admissible) is not
used.

We have $\Aut\tA_4=\DG{10}$ (the dihedral group),
and it is not very easy to describe the set
$\pat(\fiber)$; we merely compute it using \GAP~\cite{GAP4}.

\proof[Proof of \autoref{lem.pent}]
We use the thresholds $M_\fiber=14$, $M=30$, $\Mlines=28$.
At level~$0$, we add up to three pairwise disjoint extra lines in the safe
mode.
\endproof

\subsection{Sections at \pdfstr{Sigma = A\_3}{$\fiber=\tA_3$}}\label{s.sec.A3}
We order $\fiber$ cyclically.
Any $\tA_3$-graph is triangle free; assuming, in addition, that it is
$3$-admissible and, hence, biquadrangle free, see \autoref{lem.2stars},
in~\eqref{eq.sec.bound} we have, for $1\le i\ne j\le4$,
\[
\bnd_{i0}=2,\quad \bnd_{ij}=1\ \text{if $i=j\pm1\bmod4$},\quad
\bnd_{ij}=2\ \text{if $i=j+2\bmod4$}.
\label{eq.bound.A3}
\]
The other bound $\bnd_{0i}=5$ is not used.

Since $\Aut\tA_3=\DG8$, we have the following relatively simple description
of the maximal representatives $\pi\in\pat(\fiber)$: a pattern
$\pi\:a_i\mapsto n_i\in\NN$, $1\le i\le4$, belongs to $\pat(\fiber)$
if and only if
\[
n_i\le n_1\le5,\quad
n_4\le n_2,\quad\text{and}\quad
n_4\le n_3\text{ whenever }n_2=n_1.
\label{eq.pat.A3}
\]

\proof[Proof of \autoref{lem.quad}]
We use the thresholds $M_\fiber=16$,
$\Mlines=30$.

Since many graphs $\graph_4$ fail the test in \autoref{rem.validity.0},
we choose (for all $\tA_3$-graphs) the progressive mode of the algorithm
(see \autoref{ss.step.r.general}); since this version reliably lists \emph{all}
extensions of~$\graph_4$, the constant~$M$ is redundant.
In addition
to Remarks~\ref{rem.extra},
for the preliminary tests in~\autoref{ss.step.1}\iref{r.prelim} we use the
fact that the graph must be triangle free to reduce the number of adjacency
matrices $\gram=[m_{ij}]$. Most notably, we assert that
\[*
m_{ij}=0\quad\text{whenever $v_i\cap v_j\ne\varnothing$};
\]
besides, for any triple $1\le i<j<k\le\ls|\bv|$, at
least one of the three entries $m_{ij}$, $m_{ik}$, $m_{jk}$ must be~$0$.
With this reduction, the computation remains feasible, in spite of the fact
that we may have to add up to five extra vertices.
\endproof

\section{Large pencils} \label{S:largepencils}

In this section, we explain the proofs of several versions of the second inequality in
\eqref{eq.bounds}, where the target bounds $M$ and $M_\pencil$ are set
on a case-by-case basis, \cf. \autoref{rem.large.pencil} below.

We start with a \emph{fiber}
$\fiber:=\{a_1,\ldots,a_N\}$ and a
combinatorial
$\fiber$-pencil $\pencil\supset\fiber$ (see \autoref{s.pencilK3})
and set a \emph{goal}
$\ls|\graph|\ge M$ for a $\fiber$-extension $\graph\supset\pencil$.
Without loss of generality, we can assume that $\pencil$ is the
\emph{maximal} (with respect to inclusion)
pencil in~$\graph$ containing~$\fiber$; then we can state the global
criterion in terms of sections:
\[*
\Num(\graph):=\bigl(\text{$\graph$ is a $\fiber$-graph}\bigr)\mathbin\&
 \bigl(\sec_\graph\varnothing=\pencil\sminus\fiber\bigr)\mathbin\&
 \bigl(\ls|\sec_\graph\fiber|\ge M-\ls|\pencil|\bigr).
\]
We still assume that~\iref{sec.simple} and~\iref{sec.disjoint}
in \autoref{S.Msigma} hold;
besides, we consider \emph{simple} sections only: $v\in\sec a_i$, where
$a_i\in\fiber$ has multiplicity~$1$, \ie, $n_i=1$ in~\eqref{eq.kappa}.

The algorithm has up to~$N$ levels, processing one line $a\in\fiber$ at a
time. The order may differ from that fixed in \autoref{S.Msigma}; it is
controlled
by a certain permutation $\sigma\:\{1,\ldots,N\}\to\fiber$ fixed in advance.
We abbreviate
$\anchor{perm}\perm(k):=\sigma(k)$ and $\perm(a_k):=a_{\sigma(k)}$.

\remark\label{rem.large.pencil}
For $\tA_2$- and $\tA_3$-type graphs,
where we cannot assert that \emph{all} sections are
pairwise disjoint, only the first level of the algorithm works reasonably fast.
Therefore, in order to avoid the subsequent levels, we try to choose
$M_\pencil$ as large as possible, in full agreement with our general paradigm
of starting from large graphs. However, since $M_\fiber+M_\pencil\le M+1$
in~\eqref{eq.bounds}, a balance should be kept to make both this computation
and that of \autoref{S.Msigma} equally feasible. The particular values of
$M_\fiber$, $M_\pencil$ used, respectively, in \autoref{S.Msigma} and below in
this section have been found experimentally, after a few test runs of
both algorithms.
\endremark

\subsection{Level~$1$ of the algorithm}\label{s.level.1}
Let $R_1:=\range_{\perm(1)}(\varnothing\into\NN,M-\ls|\pencil|)$.
We run up to $r\smax:=\max R_1$ steps of the algorithm
in~\autoref{s.extension.algorithm}. In addition to the given fiber~$\fiber$
and admissibility level~$m$, we choose the safe mode of the algorithm
(see \autoref{ss.step.r.discrete}),
take $\num(\graph):=\bigl(\ls|\sec_\graph\perm(a_1)|\in R_1\bigr)$
for the local numeric criterion, $\graph_0:=\pencil$ for the base graph,
and the stabilizer $\stab\perm(a_1)\subset\Aut(\pencil,\fiber)$
for the group $G_0$.
The initial set is
\[*
\CS_1:=\{\perm(a_1)\}\times
 \prod_{\Delta\in\pencil_p}\Com(\Delta,1)\times
 \prod_{\Delta\in\pencil_e}\Coms(\Delta,1),
\]
where $\pencil_p$ and $\pencil_e$ are the sets of, respectively, parabolic
and elliptic connected components of~$\pencil$ other than~$\fiber$.
(Technically, when computing $\CS_1$, we remove from each component
$\Delta\in\pencil_p$ all lines of multiplicity greater than~$1$. The same
should be done for each $\Delta\in\pencil_e$, but the multiplicities are not
always known here. In any case, all ``wrong'' sections are immediately ruled out
by the Sylvester test.)

The expected result of level~$1$ is the plain output,
see~\autoref{ss.output}.
Besides, we record
\roster*
\item
all intermediate sets $\bCS_r(\perm(a_1)):=\bCS_r$ and
\item
for each graph~$\graph_1$ of rank $\rank\graph_1=19$ in the output, the set
\[
\Pat(\graph_1,G_0):=\bigl\{\max(\gpat_{\graph}\cdot G_0)\bigm|
 \graph\in\Sat_m(\graph_1)\cup\{\graph_1\}\bigr\}
\label{eq.pat.Gamma}
\]
of the maximal elements of the $G_0$-orbits of the section counts
(see \autoref{s.apriori}) of $\graph_1$ itself and all its $m$-geometric
saturations.
\endroster

\subsection{Level \pdfstr{k>=2}{$k\ge2$} of the algorithm}\label{s.level.k}
Let $\graph_{k-1}$ be one of the graphs in the output of level $(k-1)$.
As in \autoref{s.fiber.k},
denote by~$\rho$ the restriction of $\gpat_{\graph_{k-1}}$ to
$\{\perm(a_1),\ldots,\perm(a_{k-1})\}$
and let $R_k:=\range_{\perm(k)}(\rho,M_\fiber)$.
We apply to $\graph_{k-1}$ the algorithm in~\autoref{s.extension.algorithm},
with a few minor modifications explained below. We run up to
$r\smax:=\max R_k$ steps, taking
\[*
\num(\graph):=\bigl(\ls|\sec_\graph\perm(a_k)|\in R_k\bigr),\qquad
G_{k-1}:=\stab\{\perm(a_1),\ldots,\perm(a_k)\}\subset\Aut(\graph_{k-1},\fiber)
\]
(pointwise stabilizer)
for the numeric criterion and
symmetry group, respectively.

By default, we choose the safe mode (see \autoref{ss.step.r.discrete})
and plain output unless either
\roster*
\item
$k=N$, so that $\perm(a_k)$ is the last point of~$\fiber$, or
\item
$\rank\graph_{k-1}=19$ and
\[
s\smin:=\min R_k
-\max\bigl\{\pi(\perm(a_k))\bigm|\pi\in\Pat(\graph_{k-1},G_{k-2})\bigr\}>0,
\label{eq.smin}
\]
\endroster
in which case the progressive mode (see \autoref{ss.step.r.general})
and saturated output are used.
(Special arrangements can be made for particular types of~$\fiber$.)
In the case of plain output, we also record the sets
$\Pat(\graph_k,G_{k-1})$, see~\eqref{eq.pat.Gamma}.

To describe the modifications,
recall that all simple vertices of~$\fiber$ constitute a single
$(\Aut\fiber)$-orbit; hence, we can pick an element $g\in\Aut\sigma$ such
that $\perm(a_1)\cdot g=\perm(a_k)$ and consider the images
$\CS(\perm(a_k)):=\bCS(\perm(a_1))\cdot g$ (see \autoref{s.level.1};
effectively, we merely change the first entry $\perm(a_1)$ of each vertex~$v$
to~$\perm(a_k)$).
Then, we take
\[*
\CS_1:=\CS_1(\perm(a_k))\times
\Sec_{\perm(k)}(\graph_{k-1}),
\]
see~\eqref{eq.Sec},
for the initial set of vertices and, upon completion of Step~$1$, define
the set $\bCS_1(\perm(a_k))$ as the projection of $\bCS_1$ to the first factor.
At each subsequent step~$r$, instead of~\eqref{eq.Sr}, we start from
\[*
\CS_r:=\bigl(\CS_r(\perm(a_k))\cap\bCS_1(\perm(a_k))[r]\bigr)\times
\Sec_{\perm(k)}(\graph_{k-1}),
\]
thus reusing the results of level~$1$.

\remark\label{rem.virtual}
As yet another performance enhancement, we do not store (just recording their
number) the new sections if the algorithm does not improve the rank, \ie, if
either
$\rank\graph_k=\rank\graph_{k-1}$ or $\rank\graph_k=20$ for \emph{each}
graph~$\graph_k$ obtained from a given graph~$\graph_{k-1}$. This convention
simplifies the computation on all subsequent levels, as we have a smaller set
$\sec_{\graph_k}\fiber$.
\endremark

\subsection{Configurations with large pencils}\label{proofs.pencil}

\proof[Proof of \autoref{lem.d4p}]
We use the thresholds $M=28$, $M_\pencil=18$, $\Mlines=26$ and order~$\fiber$
so that the ``central'' $4$-valent fiber is~$\perm(a_5)$, the last one, so that it
is never used: we merely assume that it may have the maximal possible
valency~$7$.
\endproof

\proof[Proof of \autoref{lem.a4p}]
We use the thresholds $M=30$, $M_\pencil=17$, $\Mlines=28$.
\endproof

\proof[Proof of \autoref{lem.a3p}]
We use the thresholds $M=32$, $M_\pencil=17$, $\Mlines=30$.
For level $r\ge2$, we choose the progressive mode whenever
$\rank\graph_{k-1}=19$; however, if $r\le3$ and
$s\smin\le0$ in~\eqref{eq.smin},
the graph~$\graph_{k-1}$ is carried over to level~$4$.
\endproof

\proof[Proof of \autoref{lem.a2p}]
We
let $m=2$ and choose the thresholds $M=30$, $M_\pencil=21$, $\Mlines=29$.
Besides, we make use of the symmetry $a_2\leftrightarrow a_3$:
\roster*
\item
we do not compute the sets $\Pat(\graph_{k-1},G_{k-2})$ in~\eqref{eq.pat.Gamma},
switching to the progressive mode whenever $\rank\graph_{k-1}=19$, and
\item
we abort the computation if level~$2$ of the algorithm does not improve rank (see
\autoref{rem.virtual}).\qedhere
\endroster
\endproof

{
\let\.\DOTaccent
\def\cprime{$'$}
\bibliographystyle{amsplain}
\bibliography{degt-h5v7}
}

\end{document}